\documentclass[opre,nonblindrev]{informs3} 

\DoubleSpacedXI 

\usepackage{endnotes}
\let\footnote=\endnote

\usepackage{amsfonts}
\usepackage{geometry}
\usepackage{enumerate}
\usepackage{float}
\newcommand{\blue}[1]{{#1}}

\newcommand{\magenta}[1]{{#1}}
\usepackage[final]{showkeys}
\usepackage[colorlinks=true,linkcolor=cyan,citecolor=blue,anchorcolor=blue]{hyperref}
\newcommand{\EE}{\mathbb{E}}
\newcommand{\E}{\mathbb{E}}
\newcommand{\R}{\mathbb{R}}
\newcommand{\Z}{\mathbb{Z}}
\newcommand{\w}{(\underline{w}, \overline{w})}
\newcommand{\Prob}{\mathbb{P}}
\renewcommand{\P}{\mathbb{P}}
\newcommand{\eq}{\eqref}
\newcommand{\lipone}{\text{\rm Lip(1)}}
\newcommand{\startproof}{\proof}
\newcommand{\finishproof}{\hfill $\square$\endproof}
\providecommand{\abs}[1]{\left\lvert#1\right\rvert}
\usepackage{subcaption}

\def\be#1{\begin{equation*}#1\end{equation*}}
\def\ben#1{\begin{equation}#1\end{equation}}
\def\bes#1{\begin{equation*}\begin{split}#1\end{split}\end{equation*}}
\def\besn#1{\begin{equation}\begin{split}#1\end{split}\end{equation}}

\usepackage{natbib}
 \bibpunct[, ]{(}{)}{,}{a}{}{,}%
 \def\bibfont{\small}%
\TheoremsNumberedThrough     
\ECRepeatTheorems

\EquationsNumberedThrough    

\begin{document}


\RUNAUTHOR{Braverman, Dai, and Fang}

\RUNTITLE{High-Order Steady-State Diffusion Approximations}

\TITLE{High-Order Steady-State Diffusion Approximations}

\ARTICLEAUTHORS{%
\AUTHOR{Anton Braverman}
\AFF{Kellogg School of Management, Northwestern University, Evanston, IL 60201, \EMAIL{anton.braverman@kellogg.northwestern.edu}}  
\AUTHOR{J. G. Dai}
\AFF{School of Operations Research and Information Engineering,  Cornell University; School of Data Science and Shenzhen Research Institute of Big Data,  The Chinese University of Hong, Shenzhen, \EMAIL{jd694@cornell.edu}}
\AUTHOR{Xiao Fang}
\AFF{Department of Statistics, The Chinese University of Hong Kong, Hong Kong, \EMAIL{xfang@sta.cuhk.edu.hk}}
} 

\ABSTRACT{%
We derive and analyze new diffusion approximations of stationary distributions of Markov chains that are based on second- and higher-order terms in the expansion of the Markov chain generator. Our approximations achieve a higher degree of accuracy compared to  diffusion approximations widely used for the past fifty years, while retaining a similar computational complexity. To support our approximations, we present a combination of theoretical and numerical results across three different models. Our approximations are derived recursively through Stein/Poisson equations, and the theoretical results are proved using Stein's method.
}%


 \HISTORY{This paper was
first submitted on 12/7/2020.}

\maketitle

%



\section{Introduction}
\label{fse1}
We propose a new class of approximations for stationary
distributions of Markov chains. The new approximations will be
numerically demonstrated  to be accurate in three models:  the $M/M/n$ queue  known as the  Erlang-C model, the hospital model proposed in \cite{DaiShi2017}, and
the autoregressive (AR(1)) model studied in \cite{BlanGlyn2018}.
In addition to numerical results, for the Erlang-C model we   provide theoretical guarantees that our approximation  achieves higher-order accuracy.

 Consider a one-dimensional, positive-recurrent, discrete-time Markov chain (DTMC)  $X=\{X(n), n\geq 0\}$   taking values on a subset of $\R$. We introduce our approach in the DTMC setting, but continuous-time Markov chains (CTMC) can be treated analogously; see Section~\ref{fse3} where we treat the Erlang-C model. Call $\E \big( X(1) - X(0) | X(0) = x \big)$ the drift of the DTMC. We center and scale our DTMC by defining 
$\tilde X=\{\tilde X(n), n\geq 0\}$, where
$\tilde X(n)=\delta(X(n)-R)$ for some constants $\delta>0$ and
$R\in \R$. We typically take $R$ to be the point where the drift of $X$ equals zero, which also happens to be the equilibrium of the corresponding fluid model; c.f.,  \cite{Stol2015} or \cite{Ying2016}. The scaling parameter $\delta$ is related to stochastic fluctuations around $R$.

  Let $\tilde X(0)$ have the stationary distribution of $\tilde
X$, let $W=\tilde X(0)$, and let $W'=\tilde X(1)$.  Stationarity implies that 
\ben{\label{f0} \E f(W')-\E f(W)=0 } for all
$f:\mathbb{R}\to \mathbb{R}$ such that the expectations exist. \blue{Setting $\Delta = W' - W$, for   sufficiently smooth $f(x)$ we can expand} the left-hand side   to get
\begin{align}
0 =&\ \E f(W') - \E f(W) = \E \bigg[\sum_{i=1}^{n} \frac{1}{i!}  \Delta^i f^{(i)}(W) + \frac{1}{(n+1)!} \Delta^{n+1} f^{(n+1)}(\xi)\bigg], \quad n \geq 0, \label{eq:taylorgeneric}
\end{align}
 where  $\xi = \xi^{(n)}$ lies between $W$ and $W'$. Note that   $\Delta$ equals our scaling term $\delta$ multiplied by the displacement of the unscaled DTMC. Informally, the DTMCs we consider are those where the moments of the (unscaled) displacement are bounded by a constant independent of $\delta$, while $\delta$ itself is close to zero. In this setting, the right-hand side of \eqref{eq:taylorgeneric} is governed by its lower-order terms when $\delta$ is small. This motivates our approximations of   $W$.
%

Letting $\mathcal{W}$ be the state space  of $\tilde{X}$, for each $x \in \mathcal{W}$ let $ b(x)= \E(\Delta |W=x) $
be the drift of the DTMC at state $x$.  
Let $\w$ be the   smallest interval containing $\mathcal{W}$, and assume $b(x)$ is extended to be defined on all of $\w$. The precise form of the extension is unimportant for the time being  and we will see in our examples that this extension often has a natural form. 
We approximate $W$ by a continuous random variable $Y \in \w$ with density 
\ben{\label{f10}
\frac{\kappa}{v(x)} \exp\Big(\int_0^x \frac{b(y)}{v(y)} dy \Big),\quad  x\in \w,
}
 where $\kappa$ is the normalizing constant 
and $v(x):\w\to\R_+$ is some function to be specified. 
We note that the distribution of $Y$ is determined for a given \emph{fixed} set of system parameters of the Markov chain.
 In particular, $Y$ is well defined even when no limit is studied,  so the stationary distribution of the unscaled DTMC $X$ would then be approximated by $Y/\delta + R$. 
 
  To discuss how to choose $v(x)$, suppose $\w=\R$ and consider the   diffusion process $\{Y(t), t\ge 0\}$ given by 
\begin{align}
  Y(t) = Y(0)+ \int_0^t b(Y(s))ds + \int_0^t \sqrt{2 v(Y(s))} dB(s), \label{eq:sde}
\end{align} 
where $\{B(t), t\ge 0\}$ is the standard Brownian motion.
Under mild regularity conditions on $b(x)$ and $v(x)$, the above diffusion process is well defined and  has a unique stationary distribution whose density is given by (\ref{f10}); \blue{for a proof, see Chapter 15.5 of \cite{KarlTayl1981}}. Furthermore, the stationary density in (\ref{f10}) is characterized by
\begin{align}\label{eq:bar}
  \E b(Y)f'(Y)+ \E v(Y)f''(Y)=0 \quad \text{ for all suitable } f:\R\to \R.
\end{align}
When one or both of the endpoints of $\w$ are finite, we would account for this by adding suitable boundary reflection terms.

In this paper we   think of $Y$ as being a \emph{diffusion approximation} of $W$.
 The characterization equation (\ref{eq:bar}) is well known for Markov processes; c.f., \cite{EthiKurt1986}.  A related version  is called the \emph{basic adjoint relationship} in the context of multidimensional reflecting Brownian motions by \cite{HarrWill1987}.
Equation (\ref{eq:bar}) is known in the Stein research community as the \emph{Stein equation}; see, for example, \cite{ChenGoldShao2011}.  

Ensuring that $Y$ is a good approximation of $W$ requires a careful choice of $v(x)$. If we consider \eqref{eq:taylorgeneric} with $n = 2$  and ignore the third-order error term,  then  a natural choice is to use $v(x)=v_1(x)$, where $v_1(x)$ is an extension of $\frac{1}{2}\E(\Delta^2|W=x)$ to all of $\w$.  Choosing a diffusion approximation in such a way was done  in \cite{MandMassReim1998} and \cite{ WardGlyn2003}, as well as more recently in \cite{DaiShi2017}.

Despite this natural choice, most of the literature in the last fifty years did not use $v_1 $ to develop diffusion approximations. Instead, the typical choice is $v(x) = v_0$, where
\begin{align}
  \label{eq:v0}
  v_0=   v_1(0) = \frac{1}{2}\E(\Delta^2|W=0);
\end{align}
i.e., $v_0$ is $v_1(x)$ evaluated at the fluid equilibrium $x = 0$. For examples, see \cite{HalfWhit1981,HarrNguy1993,Gurv2014, Ward2012}.  It is usually the case that $v_0$ and $v_1(W)$ are asymptotically close,  so using $v_0$ is enough to prove a limit theorem, which is the focus of most of the diffusion approximation literature. We however, show that using  $v_0$ instead of $v_1$ can lead to significant excess error. 

\blue{One such case is the Erlang-C model. It was shown in \cite{BravDaiFeng2016} that for a large class of performance measures, the $v_0$ approximation error is at most  $C/\sqrt{R}$, where $R$ is a parameter known as the offered load and $C>0$ is a constant. In Section~\ref{fse3} we  prove this upper bound is tight. On the other hand,  the $v_1$ error vanishes at a faster rate of $1/R$. Moreover,  the $v_1$ error is much smaller than the $v_0$ error, even in cases when $R$ is small.  }  

 Given the performance of the $v_1$ approximation in the Erlang-C model, it is natural to wonder whether the $v_1$ error vanishes at a faster rate  (compared to the $v_0$ error) for other models as well. The answer is mixed; e.g., it is not true for the model in Section~\ref{fse4}.

 In this paper we provide other options for $v(x)$ beyond $v_0$ and $v_1(x)$. For $n \geq 1$, we define a $v_n$ approximation to be one that uses information from the first $n+1$ terms of the Taylor expansion in \eqref{eq:taylorgeneric}; $v_n$ approximations are not unique. We adopt the convention that   $v_n$ can refer to either the function $v_n(x)$, or the $v_n$ approximation itself.  As a preview,   we can use third-order information from the Taylor expansion is by setting
\begin{align}
  \label{eq:v2}
v(x) =   v_2(x) = v_2^{(\eta)}(x)=\max \bigg\{\frac{a(x)}{2} -\frac{b(x)c(x)}{3a(x)} -\frac{a(x)}{6}\Big(\frac{c(x)}{a(x)}\Big)', \eta \bigg\}  \quad \text{ for } x\in \w,
\end{align}
where $a(x)$ and $c(x)$ are extensions of $ \E(\Delta^2|W=x)$ and $ \E(\Delta^3|W=x)$ to $\w$, respectively, and  
$\eta>0$ is a tuneable parameter selected to keep $v_2(x)$ positive.

We formally motivate and derive \eqref{eq:v2} in Section~\ref{sec:v2def}, where we also elaborate on the need for $\eta$ and how to choose it. Going beyond $v_2$, we derive $v_3$ for the hospital model of Section~\ref{fse4} and the AR(1) model of Section~\ref{fse5}.   In both cases, numerical work suggests that finding an approximation that achieves either  a faster convergence rate of the  error to zero,  or   a significantly lower approximation error than  $v_0$,  requires us to go as far as  $v_3$.  For a discussion on how to determine which $v_n$ to use, see Section~\ref{sec:hospcompare}.  
 
  This paper is limited to the setting where the Markov chain is one-dimensional  because the derivation of $v_n$ for $n \geq 2$ exploits the one-dimensional nature of the Poisson equation; for an example, see Section~\ref{sec:v2def}. At present we do not know how to generalize this to the multidimensional setting.

The theoretical framework underpinning our work is Stein's method, which was pioneered by \cite{Stei1972}. Specifically, we use the generator comparison framework of Stein's method, which dates back to \cite{Barb1988}  and was popularized recently in queueing theory by \cite{Gurv2014}. We remark that deriving the $v_n$ approximations requires only algebra, which is handy from a practical standpoint as one can derive and implement the approximations numerically without worrying about justifying them theoretically.

In addition to deriving the $v_n$ approximations, we also use Stein's method to provide theoretical guarantees. For the Erlang-C model, Theorem~\ref{thm:md-high} establishes Cram\'er-type moderate-deviations error bounds. If $Y$ is an approximation of $W$, then moderate-deviations bounds refer to bounds on the relative error
\begin{align*}
\bigg| \frac{\Prob(Y \geq z)}{\Prob(W \geq z)} - 1 \bigg| \quad \text{ and } \quad \bigg| \frac{\Prob(Y \leq z)}{\Prob(W\leq z)} - 1 \bigg| .
\end{align*}
Compared to the Kolmogorov distance $\sup_{z \in \R} \big| \Prob(W \geq z) - \Prob(Y \geq z) \big|$, the relative error is a much more informative measure  when the value being approximated is small,
as is the case in the approximation of small tail probabilities. For
many stochastic systems modeling service operations such as customer
call centers and hospital operations, these small probabilities
represent important performance metrics; e.g.,\ at most $1\%$ of
customers waiting more than 10 minutes before getting into
service.  

To summarize, our main contribution is to present a new family of $v_n$ approximations for Markov chains. Using a combination of theoretical and numerical results, we show that the $v_n$ approximations perform significantly better than the traditional $v_0$ approximation across three separate models. Our results suggest that $v_1, v_2, v_3, \ldots$   can, and should, be used whenever possible to achieve much greater approximation accuracy. Before moving on to the main body of the paper, we first provide a brief review of related literature.

\subsection{Literature Review}
\label{sec:literature}

\paragraph{Steady-state diffusion approximations.}
In the last fifty years, diffusion approximations have been a major
research theme in the applied probability community for approximate
steady-state analysis of many stochastic systems; c.f., 
\cite{King1961a,HalfWhit1981,HarrNguy1993,MandZelt2009}.  Some of
these approximations were initially motivated by \emph{process-level
limit theorems} that establish functional central limits in certain
asymptotic parameter regions; e.g., \cite{Reim1984,Bram1998,Will1998a}. The pioneering paper of
\cite{GamaZeev2006} initiated a wave of research providing
\emph{steady-state limit theorems}, justifying steady-state
approximations on top of process-level convergence. For some examples of these, see \cite{Tezc2008, ZhanZwar2008, BudhLee2009, Kats2010, GamaStol2012, DaiDiekGao2014, Gurv2014a, YeYao2016}.

Steady-state limit theorems do \emph{not} provide a rate of
convergence or an error bound.  Recently, building on earlier work by \cite{GurvHuanMand2014},  \cite{Gurv2014} developed a general
approach to proving the rate of convergence for steady-state
performance measures of many stochastic systems.
In the setting of the $M/Ph/n+M$ queue with phase-type service time distributions,  \cite{BravDai2017} refined  the approach in
\cite{Gurv2014}, casting it into the Stein framework that has been
extensively studied in the last fifty years. The Stein framework allows
one to obtain an error bound, not just a limit theorem, for approximate
steady-state analysis of a stochastic system with a \emph{fixed} set
of system parameters. Readers are referred to \cite{BravDaiFeng2016}
for a tutorial introduction to using Stein's method for steady-state
diffusion approximations of Erlang-A and Erlang-C models, where
error bounds were established under a variety of metrics, including the Wasserstein distance, Kolmogorov distance, and moment difference.

%

\paragraph{Stein's method and moderate deviations.}

Stein's method was first introduced by \cite{Stei1972}. We refer the reader to the book by \cite{ChenGoldShao2011} for an introduction to Stein's method. 
Moderate deviations date back to
 \cite{Cram1938}, who obtained expansions for tail probabilities for sums of independent
 random variables about the normal distribution. 
 Stein's method for moderate deviations for general dependent random variables was first studied in \cite{ChFaSh13a}. See \cite{ChFaSh13b}, \cite{ShZhZh18}, \cite{Zh19}, \cite{FaLuSh19} for further developments.

\paragraph{Refined mean-field approximations.} First-order approximations, such as mean-field, or fluid model approximations capture the deterministic flow of the Markov chain while ignoring the stochastic effects.  A recent series of papers, \cite{GastHoud2017}, \cite{GastLateMass2018}, \cite{GastBortTrib2019}, explored refined mean-field approximations for computing moments of the Markov chain stationary distribution. In those papers, the authors were able to explicitly compute correction terms to the mean-field approximation, which significantly improves the accuracy of the approximation and speeds up the rate at which the approximation error converges to zero. However, the computation of these correction terms rests on  assuming that the mean-field model is globally exponentially stable and that the drift of the Markov chain is differentiable. These assumptions   fail to hold even for some basic queueing models; e.g., the Erlang-C model. 

\subsection{Notation and Organization of the Paper}
\label{sec:notation}
For $a, b \in \mathbb{R}$, we use $a^+, a^-, a \wedge b$, and $a \vee b$ to denote $\max(a,0)$, $\max(-a,0)$,  $\min(a,b)$, and $\max(a, b)$, respectively.  We adopt the convention that $\sum_{l=k_1}^{k_2}=0$ if $k_2<k_1$.   In Section~\ref{sec:v2def}, we derive several versions of $v_2$ and discuss how to analyze the approximation error using Stein's method.  In Sections~\ref{fse3}--\ref{fse5}, we study the performance of various $v_n$ approximations for three different Markov chains.   To keep the main paper  a reasonable length, some details of the proofs are left to the Appendix.

\section{Deriving the Diffusion Approximations}\label{sec:v2def}
\blue{In the previous section, we said that for $n \geq 1$, a $v_n$ approximation is one that uses information from the first $n+1$ terms of the Taylor expansion in \eqref{eq:taylorgeneric}. In this section, we justify  $v_2(x)$   proposed in \eqref{eq:v2} by tapping into the third-order terms in \eqref{eq:taylorgeneric}. For examples of accessing fourth-order terms, we refer the reader to the derivations of $v_3$ for the models in Sections~\ref{fse4} and \ref{fse5}.} What follows can be repeated for continuous-time Markov chains (CTMC), with the identity $\E G f(W)=0$ replacing  $ \E f(W') - \E f(W) = 0$, where $G$ is the generator of the CTMC.   \blue{As our starting point, we  recall from \eqref{eq:taylorgeneric} that 
\begin{align*}
0 =&\ \E f(W') - \E f(W) = \E \bigg[\sum_{i=1}^{n} \frac{1}{i!}  \Delta^i f^{(i)}(W) + \frac{1}{(n+1)!} \Delta^{n+1} f^{(n+1)}(\xi)\bigg],  
\end{align*}
where $\Delta = W' - W$, and  that $b(x), a(x)$, and $c(x)$ are   extensions of  $\E(\Delta|W=x)$, $\E(\Delta^2|W=x)$, and $ \E(\Delta^3|W=x)$ to $\w$, respectively. Let  $d(x)$ be an extension of $ \E(\Delta^4|W=x)$ to $\w$. Setting $n = 3$ in the expansion above yields }
\begin{align}
\E b(W) f'(W)+\frac{1}{2}\E a(W)f''(W) + \frac{1}{6} \E c(W) f'''(W) = -\frac{1}{24} \E d(W) f^{(4)}(\xi_1) \label{eq:taylorthird}
\end{align} 
where $\xi_1$ lies between $W$ and $W'$. We implicitly assume $f(x)$ is sufficiently differentiable and the  expectations above exist. Since $\Delta$ is small, we treat the right-hand side as error  and use the left-hand side to derive a diffusion approximation. The challenge to overcome is that the stationary density of the diffusion is characterized by \eqref{eq:bar}, which considers only the first two derivatives of a function $f(x)$, whereas the left-hand side of \eqref{eq:taylorthird} contains three derivatives. \blue{We therefore convert $f'''(W)$ into an expression involving $f''(W)$ plus some error.} Consider    \eqref{eq:taylorgeneric} again, but with $ n = 2$:
\begin{align}
\E b(W) f'(W)+\frac{1}{2}\E a(W)f''(W)    = -\frac{1}{6} \E c(W) f'''(\xi_2), \label{eq:taylorsecond}
\end{align} 
for some $\xi_2$ between $W$ and $W'$.   Fix $f(x)$  and let $g(x) =  \int_{0}^{x} \frac{c(y)}{a(y)} f''(y) dy$. Note that
\begin{align*}
g''(x)=&\ \Big( \frac{c(x)}{a(x)} f''(x)\Big)'  =\Big(\frac{c(x)}{a(x)}\Big)'f''(x)+\frac{c(x)}{a(x)}f'''(x).
\end{align*}
Evaluating \eqref{eq:taylorsecond} with $g(x)$ in place of $f(x)$ there yields 
\begin{align*}
\E \frac{b(W)c(W)}{a(W)}f''(W)+\E \frac{a(W)}{2}\Big(\frac{c(W)}{a(W)}\Big)'f''(W)+\frac{1}{2}\E c(W)f'''(W) =  -\frac{1}{6} \E c(W) g'''(\xi_2).
\end{align*} 
Rearranging terms, we have
\ben{\label{f3}
\frac{1}{6}\E c(W)f'''(W) =  -\E\Big(  \frac{b(W)c(W)}{3a(W)} +\frac{a(W)}{6}\Big(\frac{c(W)}{a(W)}\Big)'  \Big)f''(W) -\frac{1}{18} \E c(W) g'''(\xi_2).
}
Substituting \eq{f3} into \eqref{eq:taylorthird}, we obtain
\begin{align}
& \E b(W)f'(W)+\E\Big(\frac{a(W)}{2}-\frac{b(W)c(W)}{3a(W)}-\frac{a(W)}{6}\Big(\frac{c(W)}{a(W)}\Big)'\Big)f''(W) \notag \\
=&\  \frac{1}{18} \E c(W) g'''(\xi_2) -\frac{1}{24} \E d(W) f^{(4)}(\xi_1). \label{f4}
\end{align}
The left-hand side resembles the generator of a diffusion process. Define
\begin{align}
\underline{v}_2(x) =  \frac{a(x)}{2}-\frac{b(x)c(x)}{3a(x)}-\frac{a(x)}{6}\Big(\frac{c(x)}{a(x)}\Big)' , \quad x \in \w,  \label{eq:underlinev2}
\end{align}
and let $v_2(x) = (\underline v_2(x) \vee \eta)$  for some  $\eta > 0$  to recover the $v_2(x)$ in \eqref{eq:v2}. The value of $\eta$ should be chosen close to zero, and if  $\inf_{x \in \w} \underline v_2(x) > 0$, then we can pick $v_2(x) = \underline v_2(x)$.

We enforce $v(x) > 0$ because there may be issues with the integrability of the density in \eqref{f10} if $v(x)$ is allowed to be negative. For instance, in all three examples considered in this paper, $b(x) > 0$ when $x$ is to the left of the fluid equilibrium of $W$, and $b(x) < 0$ when $x$ is to the right of the fluid equilibrium; i.e.,\ the DTMC drifts back toward its equilibrium.  This drift toward the equilibrium is intimately tied to the positive recurrence of the DTMC  and can therefore be thought of as a reasonable assumption even if we go beyond this paper's three examples. Now, if $v(x)$ is allowed to be negative,   it may be that $\kappa = \infty$ in \eqref{f10}; e.g.,\ if $v(x) < 0$ for $x > K$ for some threshold $K$. Conversely, $\inf_{x \in \w} v(x) > 0$  is sufficient to ensure that $\kappa < \infty$ in all three of our examples. Another, more intuitive, reason that $v(x) > 0$  is that   a diffusion coefficient cannot be negative.

\subsection{The $v_2$ Approximation Error}\label{sec:theoryv2}
Let us discuss the error of our $v_2$ approximation. For simplicity, let us assume that $\w = \R$ and that $\inf_{x \in \R} \underline v_2(x) >0$, i.e., $v_2(x)$ equals the untruncated version $\underline v_2(x)$. We discuss in Section~\ref{sec:hybrid} what happens when the latter assumption does not hold.  Suppose $Y$ is a random variable with density as in \eqref{f10} and with $v(x)$ there equal to $v_2(x)$, \blue{i.e., 
\begin{align*}
\frac{\kappa}{v_2(x)} \exp\Big(\int_0^x \frac{b(y)}{v_2(y)} dy \Big),\quad  x\in \w,
\end{align*}}
  and assume for simplicity that $\w = \R$. Fix a test function $h: \R \to \R$ with $\E \abs{h(Y)} < \infty$,  and let $f_h(x)$ be the solution to the Poisson equation 
\begin{align}
b(x) f_h'(x) + v_2(x) f_h''(x) = \E h(Y) - h(x), \quad x \in \R. \label{eq:poissonde}
\end{align}
Assume that $\E \abs{f_h (W)} < \infty$, which is typically true in practice, and take expected values with respect to $W$ to get
\begin{align*}
\E h(Y) - \E h(W) =&\  \E b(W) f_h'(W) + \E v_2(W) f_h''(W) =   \frac{1}{18} \E c(W) g_{h}'''(\xi_2) -\frac{1}{24} \E d(W) f_{h}^{(4)}(\xi_1).
\end{align*}
The last equality follows from \eqref{f4}, and $g_h(x)  =  \int_{0}^{x} \frac{c(y)}{a(y)} f_h''(y) dy$. We have again made an implicit assumption that   $f_h(x)$ is sufficiently regular. The regularity of $f_h(x)$ is entirely determined  by the regularity of $b(x)$, $v(x)$, and $h(x)$.  
The right-hand side equals
 \begin{align}
&\frac{1}{18} \E \Big[c(W) \Big(\frac{c(x)}{a(x)}\Big)''\Big|_{x = \xi_2}f_h''(\xi_2)\Big]  + \frac{2}{18} \E \Big[ c(W) \Big(\frac{c(x)}{a(x)}\Big)'\Big|_{x = \xi_2} f_h'''(\xi_2)\Big]    \notag \\
&+  \frac{1}{18} \E \Big[c(W)  \frac{c(\xi_2)}{a(\xi_2)} f_h^{(4)}(\xi_2)\Big]  -\frac{1}{24} \E d(W) f_h^{(4)}(\xi_1) \label{eq:v2taylorerror}
\end{align}
because 
\begin{align*}
g_h'''(x) =&\ \Big( \frac{c(x)}{a(x)} f_h''(x)\Big)'' = \Big(\frac{c(x)}{a(x)}\Big)''f_h''(x) + 2\Big(\frac{c(x)}{a(x)}\Big)'f_h'''(x)+\frac{c(x)}{a(x)}f_h^{(4)}(x).
\end{align*}
Note that \eqref{eq:v2taylorerror} contains a term involving $f_h''(x)$ that is not captured by $\underline v_2(x)$. To capture that term, we can consider
\begin{align*}
\overline v_2(x) =  \frac{a(x)}{2}-\frac{b(x)c(x)}{3a(x)}-\frac{a(x)}{6}\Big(\frac{c(x)}{a(x)}\Big)' - \frac{1}{18} c(x) \Big(\frac{c(x)}{a(x)}\Big)'' , \quad x \in \R.
\end{align*}
Truncating $\overline v_2(x)$ produces yet another $v_2$ approximation with error
\begin{align}
&\frac{1}{18} \E \Big[c(W) \Big( \Big(\frac{c(x)}{a(x)}\Big)''\Big|_{x = \xi_2}f_h''(\xi_2) - \Big(\frac{c(W)}{a(W)}\Big)'' f_h''(W)\Big)\Big]     \notag \\
&+ \frac{2}{18} \E \Big[ c(W) \Big(\frac{c(x)}{a(x)}\Big)'\Big|_{x = \xi_2} f_h'''(\xi_2)\Big] +  \frac{1}{18} \E \Big[c(W)  \frac{c(\xi_2)}{a(\xi_2)} f_h^{(4)}(\xi_2)\Big]  -\frac{1}{24} \E d(W) f_h^{(4)}(\xi_1) \label{eq:v2taylorerroralt} 
\end{align}
in place of \eqref{eq:v2taylorerror}. In order to decide between $\underline v_2(x)$ and $\overline v_2(x)$, let us compare the two error terms in \eqref{eq:v2taylorerror} and \eqref{eq:v2taylorerroralt}. We stress that the following is   an informal discussion meant to develop  intuition. Theoretical guarantees for $v_2$ must be established on a case-by-case basis and fall outside the scope of this paper. 

 Consider first the error term \eqref{eq:v2taylorerror}. Recall that $a(x),c(x)$, and $d(x)$ equal $\E(\Delta^k|W=x)$ for $k = 2,3,4$, respectively. Now $\Delta = W' - W$ equals $\delta$ times the one-step displacement of the Markov chain. Let us  assume that the displacement is bounded by a constant independent of $\delta$, in which case $\E(\Delta^k|W=x)$ shrinks at the rate of at least $\delta^k$ as $\delta \to 0$. In particular, $d(x)$ shrinks at least as fast as $\delta^4$. Since $a(x)$ is the extension of the strictly positive function $\E (\Delta^{2} | W=x)$, we assume that this extension is also strictly positive. Furthermore, we assume that $a(x)$ is of order $\delta^2$, as opposed to merely shrinking at a rate of at least $\delta^2$. Formally, we assume that $\inf\{\delta^{-2} \abs{a(x)} : \delta \in (0,1),\ x \in \w\} > 0$, which implies that, provided the derivatives exist, $c(x) \frac{c(x)}{a(x)}$, $c(x) \Big(\frac{c(x)}{a(x)}\Big)'$, and $c(x) \Big(\frac{c(x)}{a(x)}\Big)''$ all shrink at a rate of at least $\delta^4$ as $\delta \to 0$, making them comparable to $d(x)$. 

Now consider the error term \eqref{eq:v2taylorerroralt}, focusing on the first line there. Provided $a(x),c(x),$ and $f_h(x)$ are sufficiently differentiable, the mean value theorem implies 
\begin{align*}
&\E \Big[c(W) \Big( \Big(\frac{c(x)}{a(x)}\Big)''\Big|_{x = \xi_2}f_h''(\xi_2) - \Big(\frac{c(W)}{a(W)}\Big)'' f_h''(W)\Big)\Big]  \\
=&\ \E \Big[c(W) (\xi_2 - W) \Big( \Big(\frac{c(x)}{a(x)}\Big)'' f_h''(x)\Big)'\Big|_{x = \xi_3}\Big].
\end{align*}
Under the two assumptions from before, the terms in front of the derivatives of $f_h(x)$ above   shrink at a rate of at least $\delta^5$. If the rest of the terms in  \eqref{eq:v2taylorerroralt} shrink at the rate of $\delta^4$, then using $\overline v_2(x)$ instead of $\underline v_2(x)$ as the $v_2$ approximation would not make the  error converge to zero faster. For this reason and also because $\underline v_2(x)$ is simpler than $\overline v_2(x)$, we work with $\underline v_2(x)$ in the  models we consider.

\section{Erlang-C Model}\label{fse3}
\blue{In this section we consider the Erlang-C model. We prove that the $v_1$ error converges to zero at a faster rate than the $v_0$ error. We also conduct  numerical experiments where we observe that the $v_1$ error is much smaller, often by a factor of 10, than the $v_0$ error. After defining the model, we introduce the approximations  in Section~\ref{sec:ecv0v1} and then present theoretical and numerical results in Sections~\ref{sec:erlangctheory} and \ref{sec:erlangcnumerical}, respectively. } 

\blue{The Erlang-C, or $M/M/n$, system has a single buffer served by $n$ homogeneous servers working in} a first-come-first-served manner. Customers arrive according to a Poisson process with rate $\lambda$, and   service times are  i.i.d., exponentially distributed with mean $1/\mu$. We let $R = \lambda/\mu$ and   $\rho = \frac{\lambda}{n\mu} = R/n$ be the  offered load and utilization, respectively.

Let $X(t)$ be the number of customers in the system at time $t$. We assume that $\rho < 1$, implying that $X = \{X(t), t \geq 0\}$ a  positive recurrent CTMC. Set $\delta = 1/\sqrt{R}$, $\tilde X = \{ \tilde X(t) = \delta(X(t) - R),\ t \geq 0\}$, and let $W$ be the random variable having the stationary distribution of $\tilde X$. The support of $W$  is $\mathcal{W} =\{\delta(k-R): k\in \Z_+\}$, so we let
\begin{align}
\w = (-\delta R, \infty) =  (-\sqrt{R}, \infty). \label{eq:smallestinterval}
\end{align}
The generator of $\tilde X$ satisfies (cf. Eq. (3.6) of \cite{BravDaiFeng2016})
\ben{ \label{eq:gx}
G_{\tilde X}f(x)=\lambda (f(x+\delta)-f(x))+\mu\big[ (x/\delta+R) \wedge n  \big] (f(x-\delta)-f(x)), 
}
where   $x = \delta (k-R)$  for some integer $k \geq 0$. Proposition 1.1 in \cite{Hend1997} states that 
\ben{ \label{eq:gxbar}
\E G_{\tilde X} f(W) = \E \Big[ \lambda (f(W+\delta)-f(W))+\mu\big[  (W/\delta+ R) \wedge n  \big] (f(W-\delta)-f(W))   \Big]=0
}
for all $f(x)$ such that $\E \abs{f(W)} < \infty$. 
\subsection{The $v_0$ and $v_1$ Approximations}
\label{sec:ecv0v1}
Let us perform Taylor expansion on the left-hand side of \eqref{eq:gxbar}:
\begin{align}
 \E b(W) f'(W)+\E \frac{a(W)}{2}f''(W) =&\  -\frac{1}{6}\big( \delta^3 \lambda f^{'''}(\xi_1) - \delta^3 \mu\big[  (W/\delta+ R) \wedge n  \big]f^{'''}(\xi_2)\big), \label{f6}
\end{align}
where $\xi_{1} \in (W,W+\delta)$, $\xi_2 \in (W-\delta, W)$, 
\begin{align}
  &  b(x) = \delta \big(\lambda -\mu \big[(x/\delta+R)\wedge n\big]\big), \quad \text{ and } \label{eq:tb}\\
  & a(x) = \delta^2 \big(\lambda+ \mu \blue{[(x/\delta+R)\wedge n]} \big) = 2\mu - \delta b(x), \quad x\in \mathcal{W}. \label{eq:ta} 
\end{align}
The second equality in \eqref{eq:ta} holds because $\delta^2 = 1/R = \mu/\lambda$. Let
\begin{align}
   \beta= \delta(n-R)>0,  \quad \text{ or } \quad n = R +\beta \sqrt{R}. \label{eq:staffing}
\end{align}  
   When $\beta$ is fixed and $R,n \to \infty$, the asymptotic regime is known as the Halfin-Whitt regime; see \cite{HalfWhit1981}. It is also
known as the \emph{quality  and efficiency--driven regime}  because in
this parameter region, the system simultaneously achieves short
average waiting time (quality) and high server utilization
(efficiency); \cite{GansKoolMand2003}. Some of our results assume that $\beta$ is fixed, while others do not.  

By considering the cases when $x \leq \beta$ and $x > \beta$ in \eqref{eq:tb}, we see that $b(x) = -(\mu x \wedge \mu \beta)$  for $x \in \mathcal{W}$, and we extend $b(x)$ to the entire real line via
\begin{align}
b(x) = -(\mu x \wedge \mu \beta), \quad x \in \R. \label{eq:bdef} 
\end{align}
We also want a strictly positive extension of  $a(x)$ to $\R$. Since $\mathcal{W} \subset [ -\sqrt{R}, \infty)$, we define 
\begin{align}
a(x) = 2\mu - \delta b(-\sqrt{R} \vee x), \quad x \in \R, \label{eq:adef}
\end{align}
and since $b(x)$ is nonincreasing and $b(-\sqrt{R}) = \mu \sqrt{R} = \mu/\delta$, we have $a(x) \geq a(-\sqrt{R}) = \mu$. Recall from \eqref{f10} that our diffusion approximations all have density of the form 
\begin{align*}
\frac{\kappa}{v(x)} \exp\Big(\int_0^x \frac{b(y)}{v(y)} dy \Big),\quad  x\in \R,
\end{align*}
for some normalizing constant $\kappa > 0$.  The $v_0$ and $v_1$ approximations are obtained by setting 
\begin{align*}
v(x) = v_0 = \frac{1}{2} a(0) = \mu  \quad \text{ and } \quad v(x) = v_1(x) = \frac{1}{2} a(x), \quad x \in \R. 
\end{align*}
Let $Y_0$ and $Y_1$ be the random variables corresponding to  $v_0$ and $v_1$, respectively.  

\begin{remark}
To better approximate $W$, we can use a diffusion process defined on $[-\sqrt{R},\infty)$ with a reflecting condition at the left boundary of $x = -\sqrt{R}$. However, our theorems in Section~\ref{sec:erlangctheory} are intended for the asymptotic regimes when $R \to \infty$.  Since the probability of an empty system shrinks rapidly as $R$ grows,   the choice between a reflected diffusion on $[-\sqrt{R}, \infty)$ and a diffusion defined on $\R$ is inconsequential.  
\end{remark}

\subsection{Theoretical Guarantees for the Approximations}\label{sec:erlangctheory}
We now present several theoretical results showing that the $v_1$ error vanishes faster than the $v_0$ error. Define the class of all Lipschitz-$1$ functions by
\begin{align*}
\lipone = \big\{h: \R \to \R \ \big|\  \abs{h(x)-h(y)} \leq \abs{x-y} \text{ for all $x,y\in \R$}\big\}.
\end{align*}
It was shown in \cite{BravDaiFeng2016} that 
\begin{align}
\sup_{h \in \lipone} \big| \EE h(W) - \EE h(Y_0) \big| \leq \frac{205}{\sqrt{R}}, \quad \text{ if } R < n. \label{eq:oldmain}
\end{align}
The quantity on the left-hand side above is known as the Wasserstein distance and, as was shown in \cite{GibbSu2002}, convergence in the Wasserstein distance implies convergence in distribution. To add to the result of \cite{BravDaiFeng2016}, we prove the following lower bound in Section~\ref{sec:lowerbound} of the electronic companion.
\begin{proposition}
\label{prop:lowerbound}
Assume $n = R + \beta \sqrt{R}$ for some fixed $\beta > 0$. 
There exists a constant $C(\beta) > 0$ depending only on $\beta$ such that
\begin{align*}
 \big| \EE W - \EE  Y_0  \big| \geq \frac{C(\beta)}{\sqrt{R}}.
\end{align*}
\end{proposition}
An immediate implication of Proposition~\ref{prop:lowerbound} is that the Wasserstein distance between $W$ and $Y_0$ is at least $C(\beta)/\sqrt{R}$. The assumption that $\beta$ is fixed can likely be removed (with additional effort), but that is not the focus of our paper.  We turn to the $v_1$ approximation. Define $W_2 = \{h: \R \to \R\ |\  h(x), h'(x) \in \lipone \}$  and for two random variables $U,V$, define the $W_2$ distance as
\begin{align*}
 d_{W_2}(U,V) = \sup_{h \in W_2} \big| \EE h(U) - \EE h(V) \big|.
\end{align*}
Although $W_2 \subset \lipone$, it still rich enough to imply convergence in distribution. In particular, Lemma 3.5 of \cite{Brav2017} shows that by approximating the indicator function of a half line by Lipschitz functions with bounded second derivatives,  convergence
in the $d_{W_2}$ distance implies convergence in distribution. The following result first appeared as Theorem 3.1 in \cite{Brav2017}.
\begin{theorem}
\label{thm:w2}
There exists a constant $C > 0$ (independent of $\lambda, n$, and $\mu$)  such that for all $n \geq 1, \lambda > 0$, and $\mu > 0$ satisfying $1 \leq R < n $,
\begin{align*}
 \sup_{h \in W_2} \big| \EE h(W) - \EE h(Y_1) \big| \le \frac{C}{R}.
\end{align*}
\end{theorem} 
Note that $h(x) = x$ belongs to $W_2$,  so Theorem~\ref{thm:w2} and Proposition~\ref{prop:lowerbound} tell us that the the $v_1$ approximation error of $\E (W)$ is guaranteed to vanish faster than the $v_0$ error as $R \to \infty$. 

Error bounds of the flavor of Theorem~\ref{thm:w2} were established in \cite{GurvHuanMand2014, Gurv2014, BravDai2017, BravDaiFeng2016}, all of which studied convergence rates for steady-state diffusion approximations of various models. The rate of $1/R$ is an order of magnitude better than the rates in any of the previously mentioned papers, where  the authors  obtained rates  that would be equivalent to $1/\sqrt{R}$ in our model.

Going beyond error bounds for smooth test functions, we now present moderate-deviations bounds for our two approximations. Namely, we are interested in the relative error of approximating the cumulative distribution function (CDF) and complementary CDF (CCDF). We define the relative error of the right tail to be 
\begin{align*}
\abs{\frac{\Prob(Y_i \geq z)}{\Prob(W \geq z)} - 1}, \quad i = 0,1.
\end{align*}
The relative error for the left tail is defined similarly. The first result is for $v_0$. 
\begin{theorem} 
\label{thm:md-std}
Assume that $n = R + \beta \sqrt{R}$ for some fixed $\beta > 0$. \blue{There exist positive constants $c_0$ and $C$ depending only on $\beta$ such that 
  \begin{align}
  & \left|\frac{\Prob(Y_0\geq z)}{\Prob(W\geq z)}-1\right|\leq \frac{C}{\sqrt{R}}\left(1+z\right) \quad \text{ for } 0<z\leq c_0 R^{1/2}\ \text{and} \label{f12} \\
&\left|\frac{\Prob(Y_0\leq -z)}{\Prob(W\leq -z)}-1\right|\leq \frac{C}{\sqrt{R}}\left(1+z^3\right),\ \text{ for } 0<z\leq \min\{c_0 R^{1/6}, R^{1/2}\}. \label{f15}
  \end{align}
}
\end{theorem}
The second result presents analogous bounds for the $v_1$ approximation. 
\begin{theorem}
\label{thm:md-high}
Assume $n = R + \beta \sqrt{R}$ for some fixed $\beta > 0$. 
\blue{There exist positive constants $c_1$ and $C$ depending only on $\beta$ such that 
  \begin{align}
  & \left|\frac{\Prob(Y_1\geq z)}{\Prob(W\geq z)}-1\right|\leq \frac{C}{\sqrt{R}}\left(1+\frac{z}{\sqrt{R}}\right) \quad \text{ for } 0< z\leq c_1 R\ \text{and} \label{f13}\\
&\left|\frac{\Prob(Y_1\leq -z)}{\Prob(W\leq -z)}-1\right|\leq \frac{C}{\sqrt{R}}\left(1+z+\frac{z^4}{\sqrt{R}}\right),\ \text{ for } 0<z\leq \min\{c_1 R^{1/4}, R^{1/2}\} \label{f16}.
  \end{align}
} 
\end{theorem}
Inequality \eqref{f13} follows from Theorem~4.1 of \cite{Brav2017}.  We prove \eqref{f16} in   Section \ref{fse8};
Theorem \ref{thm:md-std} follows from a similar and simpler proof in Section~\ref{fap3}. 

These are called moderate deviations bounds because they cover the case when $z$ is ``moderately'' far from the origin, with ``moderately'' being quantified by intervals of the form $z \in [0,c_0R^{1/2}]$, $z \in [0,c_1R]$, etc. In contrast,  large-deviations results focus on understanding the behavior of $\Prob(W \geq z)$ as $z \to \infty$. To compare the two theorems, suppose $z = c_0 \sqrt{R}$ and consider the upper bounds in \eqref{f12} and \eqref{f13}. The $v_0$ error is guaranteed to be bounded as $R$ grows, while the $v_1$ error shrinks at a rate of at least $1/\sqrt{R}$.

\subsection{Numerical Results}\label{sec:erlangcnumerical}
\begin{figure}
\centering
\includegraphics[scale=0.24]{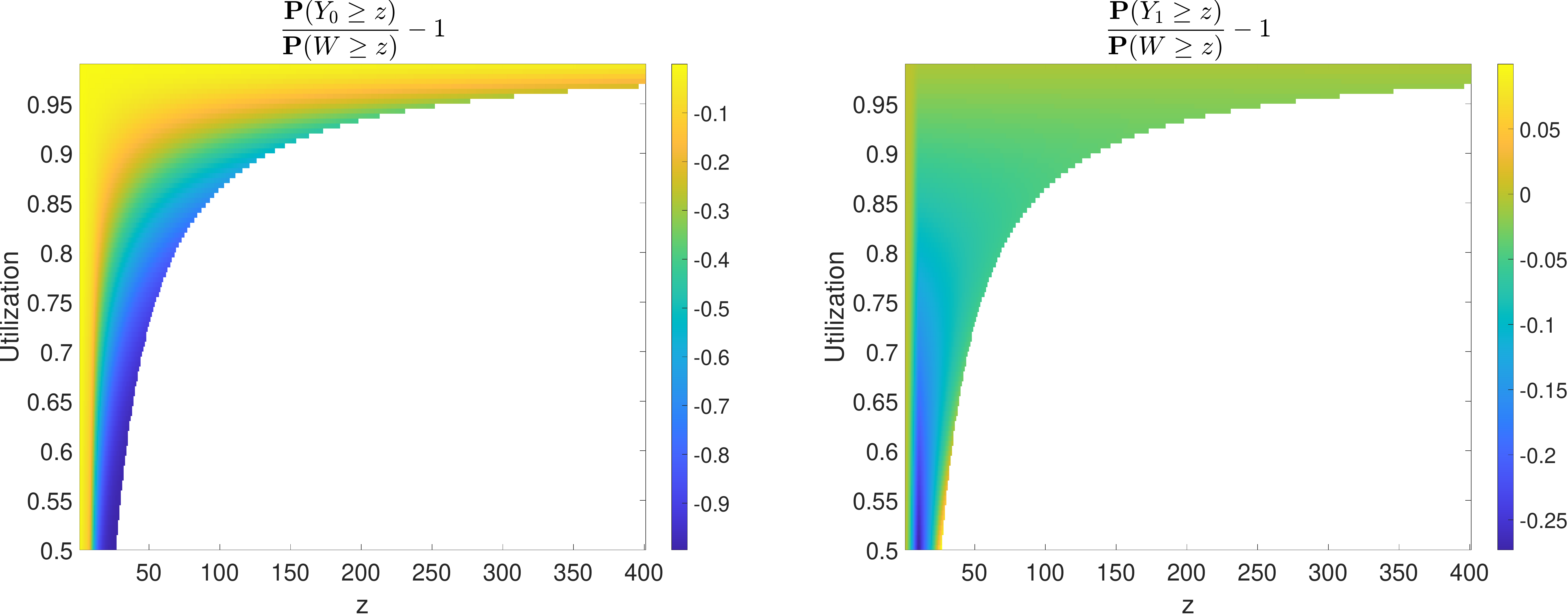}
\caption{The $v_1$ approximation is much more accurate. }
 \label{fig:ecmd1}
\end{figure}

Although the Erlang-C system depends on three parameters $\lambda$, $\mu$, and $n$,  the stationary distribution  depends on only  $\rho$ and $n$; see Appendix C in \cite{Alle1990}. Figure~\ref{fig:ecmoment} displays the relative errors of the $v_0$ and $v_1$ approximations of $\E (W)$  when $5 \leq n \leq 100$ and $0.5 \leq \rho \leq 0.99$. Note that the $v_1$ error is   about ten times smaller. 
We also compare how well $v_0$ and $v_1$ approximate the CCDF of $W$. Figure~\ref{fig:ecmd1}  plots the relative error of approximating $\Prob(W/\delta +R \geq z)$ for various values of $z$ when $n = 10$ and   $\rho \in [0.5,0.99]$, and shows that the  $v_1$ error is again much smaller. In results not reported in the paper, we observed that the $v_1$ error remains much smaller even as we vary $n$. 
\begin{figure}[h!]
\centering
	\includegraphics[scale=0.55]{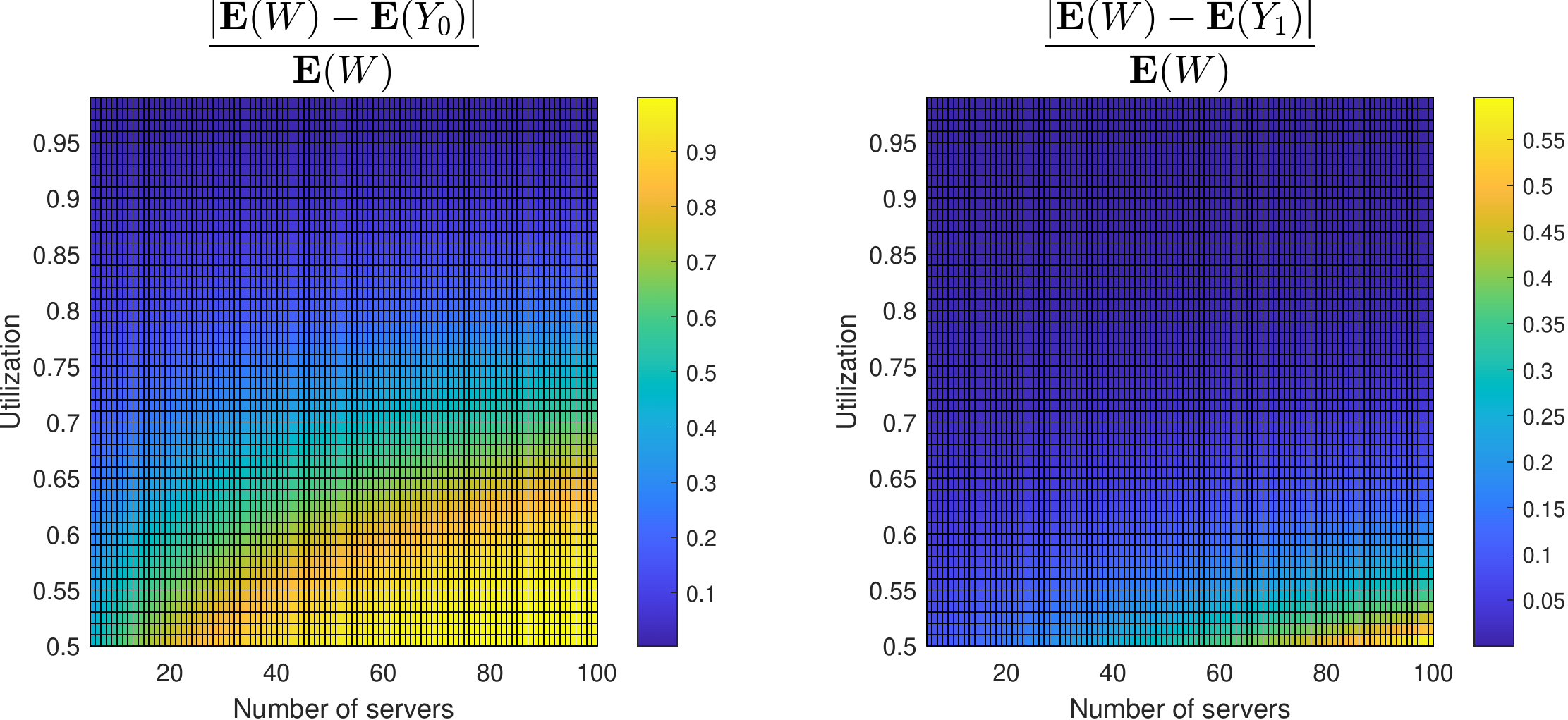}
	\caption{The errors increase towards the bottom-right corner of each plot. This is due to the fact that $\E (W)$ is very close to zero in that region and not because  approximations perform poorly. \label{fig:ecmoment}}
\end{figure}

%
%
%
%

\section{Hospital Model}\label{fse4}
In this section we consider the discrete-time model for hospital inpatient flow proposed by \cite{DaiShi2017}. Numerical experiments presented later  suggest that both the $v_1$ and  $v_2$ errors vanish at the same rate as the $v_0$ approximation error. To observe a faster convergence rate, we have to resort to the $v_3$ approximation.

Consider a discrete-time queueing model with $N$ identical servers. Let $X(n)$ be the number of customers in the system at the end
of time unit $n$. Given $X(0)$, we define 
\begin{align*}
X(n) = X(n-1) + A(n) - D(n), \quad n \geq 1,
\end{align*}
where $A(n) \sim$ Poisson$(\Lambda)$ represents the number of new arrivals in the time period $[n-1,n)$. At the end of each time period, every customer in service flips a coin and, with probability $\mu \in (0,1)$, departs the system at the start of the next time period. Thus, conditioned on $X(n-1) = k$, we have $D(n) \sim $ Binomial$(k\wedge N,\mu)$. Assuming $ \Lambda < N \mu$,   \cite{DaiShi2017} showed that
  $X=\{X(n): n=1, 2, \ldots\}$ is a positive-recurrent DTMC.

\blue{This DTMC is similar to the Erlang-C model, but unlike the Erlang-C model where the customer count only changes by one at a time, the jump size $X(1) - X(0)$ is not bounded because $A(n)$ is unbounded. As a result, computing the stationary distribution takes a long time when  $\Lambda$  is large and the utilization $\rho = \Lambda/(N\mu)$ is near one because the state-space truncation has to be large to account for potential arrivals. }

 We are interested in the scaled DTMC  $\tilde X = \{ \tilde X(n) = \delta (X(n) - N)\}$.  To stay consistent with \cite{DaiShi2017}, we  center $X(n)$ around $N$. We consider the  parameter ranges  studied in \cite{DaiShi2017}, which, given some constant  $\beta>0$,   are
\begin{align}\label{eq:hosqed}
\Lambda=\sqrt{N}-\beta, \quad \mu=\delta=1/\sqrt{N}   
\end{align}

\subsection{Motivating the Need for a $v_3$ Approximation}
\label{sec:hospcompare}
This section contains an informal discussion aimed at explaining why the  $v_0$, $v_1$, and  $v_2$ errors vanish at the same rate of $\delta = 1/\sqrt{N}$, and why we need the $v_3$ approximation to observe a convergence rate of $\delta^2 = 1/N$.

Initialize $\tilde X(0)$ according to the stationary distribution of $\tilde X$, let $W = \tilde X(0)$, $W' = \tilde X(1)$, and set $\Delta = W' - W$. The support of $W$ is $\mathcal{W}=\{ \delta(k-N): k \in \mathbb{N} \} \subset [-\delta N, \infty)$. As we are accustomed to doing  by this point, we let $b(x) =  \E(\Delta  | W = x)$. We know from (37) and (38) of \cite{DaiShi2017} that for $x \in \mathcal{W}$, 
\begin{align}
 b(x) = \E(\Delta| W= x) =&\  \delta( x^{-} - \beta ), \quad \text{ and }  \label{eq:1}\\ 
   \E(\Delta^2|W=x)  =&\ 2\delta + \big(b^{2}(x)-\delta b(x) -\delta^2 - 2\delta^2\beta\big) + \delta^3 x^{-}.  \label{eq:2}
\end{align}
For the higher moments of $\Delta$, let us use $\epsilon(x)$ to represent a generic function that may change from line to line  but always satisfies the property 
\begin{align}
\abs{\epsilon(x)} \leq C (1 + \abs{x})^{5} \label{eq:epsdef}
\end{align}
for some constant $C>0$ that depends only on $\beta$. We show in Section~\ref{sec:hosmecproof} that 
\begin{align}
\E(\Delta^3|W=x) = \delta^2 \epsilon(x), \quad \E(\Delta^4|W=x) = \delta^2 \epsilon(x), \quad \E(\Delta^5 |W=x) = \delta^3 \epsilon(x), \quad x \in \mathcal{W}. \label{eq:345}
\end{align}
%
 All of our $v_n$  approximations share the same drift $b(x)$, but the diffusion coefficients vary. As always,  $v_1(x) = \frac{1}{2} \E(\Delta^2|x  )$. Since $b(x) = 0$ at $x = -\beta$, \eqref{eq:2} implies that
\begin{align}
v_0 =&\ v_1(-\beta) = \frac{1}{2}\big( 2\delta - \delta^2 (1+2\beta) + \delta^3 \beta\big). \label{eq:hospv0}
\end{align} 
The following informal discussion assumes that all functions are sufficiently differentiable and that all expectations exist. Our starting point, as always, is the Poisson equation 
\begin{align}
b(x) f_h'(x) + v(x) f_h''(x) = \E h(Y) - h(x), \quad x \in \R, \label{eq:hosppoisson}
\end{align} 
where  $v(x)$ is a temporary placeholder and $Y$ has density given by \eqref{f10}. The Taylor expansion in \eqref{eq:taylorgeneric}  tells us that 
\begin{align*}
\E b(W) f_h'(W) + \E \bigg[\sum_{i=2}^{n} \frac{1}{i!}  \Delta^i f_h^{(i)}(W) + \frac{1}{(n+1)!} \Delta^{n+1} f_h^{(n+1)}(\xi)\bigg] = 0,
\end{align*} 
where $\xi = \xi^{(n)}$ lies between $W$ and $W'$. Subtracting this equation from \eqref{eq:hosppoisson} and taking expected values there with respect to $W$ we see that for $n \geq 1$, 
\begin{align}
\E h(W) - \E h(Y) =&\ -\E v(W) f_h''(W)  + \E \bigg[\sum_{i=2}^{n} \frac{1}{i!}  \Delta^i f_h^{(i)}(W) + \frac{1}{(n+1)!} \Delta^{n+1} f_h^{(n+1)}(\xi)\bigg]. \label{eq:hosperr}
\end{align}
Consider the $v_0$ error. When $v(x) = v_0$, equation \eqref{eq:hosperr} with $n = 2$ there becomes
\begin{align}
\E h(W) - \E h(Y) =&\    \E  \Big(\frac{1}{2}\Delta^2 - v_0 \Big) f_h''(W) + \frac{1}{6} \E \Delta^3 f_h'''(\xi). \label{eq:v0errorhosp}
\end{align}
Note that $f_h(x)$ depends on $v(x)$ because it solves \eqref{eq:hosppoisson}.
In Lemma~3 of  \cite{DaiShi2017}, the authors proved that   $\abs{f_h''(x)} \leq C/\delta$ and $\abs{f_h'''(x)} \leq C/\delta$ for some constant $C > 0$ dependent only on $\beta$. Assuming that $f_h''(x)$ and $f_h'''(x)$ indeed grow at the rate of $C/\delta$ as $\delta \to 0$, the forms of $\E (\Delta^2 | W = x)$ and $v_0$ in \eqref{eq:2} and  \eqref{eq:hospv0} yield
\begin{align*}
  \frac{1}{2}\E (\Delta^2 | W = x) - v_0 =&\  \frac{1}{2}  \Big(b^{2}(x)-\delta b(x)   + \delta^3 (x^{-}-\beta)\Big) =   \frac{1}{2}  \Big(b^{2}(x)-\delta b(x)   + \delta^2 b(x)\Big).
\end{align*}
Since $b(x) = \delta(x^{-}-\beta)$, this quantity is of order $\delta^2$. Now \eqref{eq:345} says that  $\E (\Delta^3| W =x)$ is also of order $\delta^2$. Therefore, we expect both $ \E  \big(\Delta^2/2 - v_0 \big) f_h''(W)$ and $\E \Delta^3 f_h'''(\xi)$ to be of order $\delta$, so even if $ \E  \big(\Delta^2/2 - v_0 \big) f_h''(W)$ were not present in \eqref{eq:v0errorhosp}, the approximation error would still be of order $\delta$ due to $\E \Delta^3 f_h'''(\xi)$. We believe this is why the $v_0$ and $v_1$ errors appear to vanish at the same rate  despite $v_1(x)$ capturing the entire second order term of
\begin{align}
 \E \bigg[\sum_{i=2}^{n} \frac{1}{i!}  \Delta^i f_h^{(i)}(W) + \frac{1}{(n+1)!} \Delta^{n+1} f_h^{(n+1)}(\xi)\bigg]. \label{eq:hosptaylorerror}
\end{align}

Going beyond $v_0$ and $v_1$, we see that for the error to be of order $\delta^2$, the diffusion approximation must  capture all the terms in \eqref{eq:hosptaylorerror}  that are of order $\delta$. If we assume for the moment that the derivatives of $f_h(x)$ are all of order $1/\delta$, we see that our approximation has to capture all terms of order $\delta^2$ or larger in the functions $\{ \E(\Delta^i|x)\}_{i=1}^{\infty}$.  
 
From \eqref{eq:1}, \eqref{eq:2}, and  \eqref{eq:345} we see that  $\E (\Delta^{1} | x)$ through  $\E (\Delta^{4} | x)$ are all of order $\delta$ or $\delta^2$, while $\E(\Delta^{5} | x)$ is of order $\delta^3$. Thus, if we want the error to be of order $\delta^2$, our approximation must capture the terms of order $\delta$ and $\delta^2$ in $\E (\Delta^{1} | x)$ through  $\E (\Delta^{4} | x)$  and can ignore terms of order $\delta^3$ like $\E(\Delta^{5} | x)$.   We remark that $v_2(x)$ in \eqref{eq:v2} depends only on  $\E (\Delta^{1} | x)$ through  $\E (\Delta^{3} | x)$ and not  on $\E (\Delta^{4} | x)$. We suspect this is  why we observe the $v_2$ error to be of order $\delta$. In Section~\ref{app:hospital_proofs} we derive a $v_3$ approximation of the form 
\begin{align}
  \label{eq:hosv3}
  v_3(x)=\max\Big\{ \delta +\frac{1}{2} \Big(\delta^{2} 1(x<0)- \delta b(x)-\delta^2-2\delta^2\beta\Big),  \delta/2  \Big\}
\end{align}
and in the following section we present numerical results that suggest that the error of this approximation converges to zero at a rate of $\delta^2$.
 
\subsection{Numerical Results} \label{sec:hospital_numeric}
\begin{figure}[H]
\centering
\begin{subfigure}[t]{0.45\textwidth}
\centering
    \includegraphics[width=.8\linewidth]{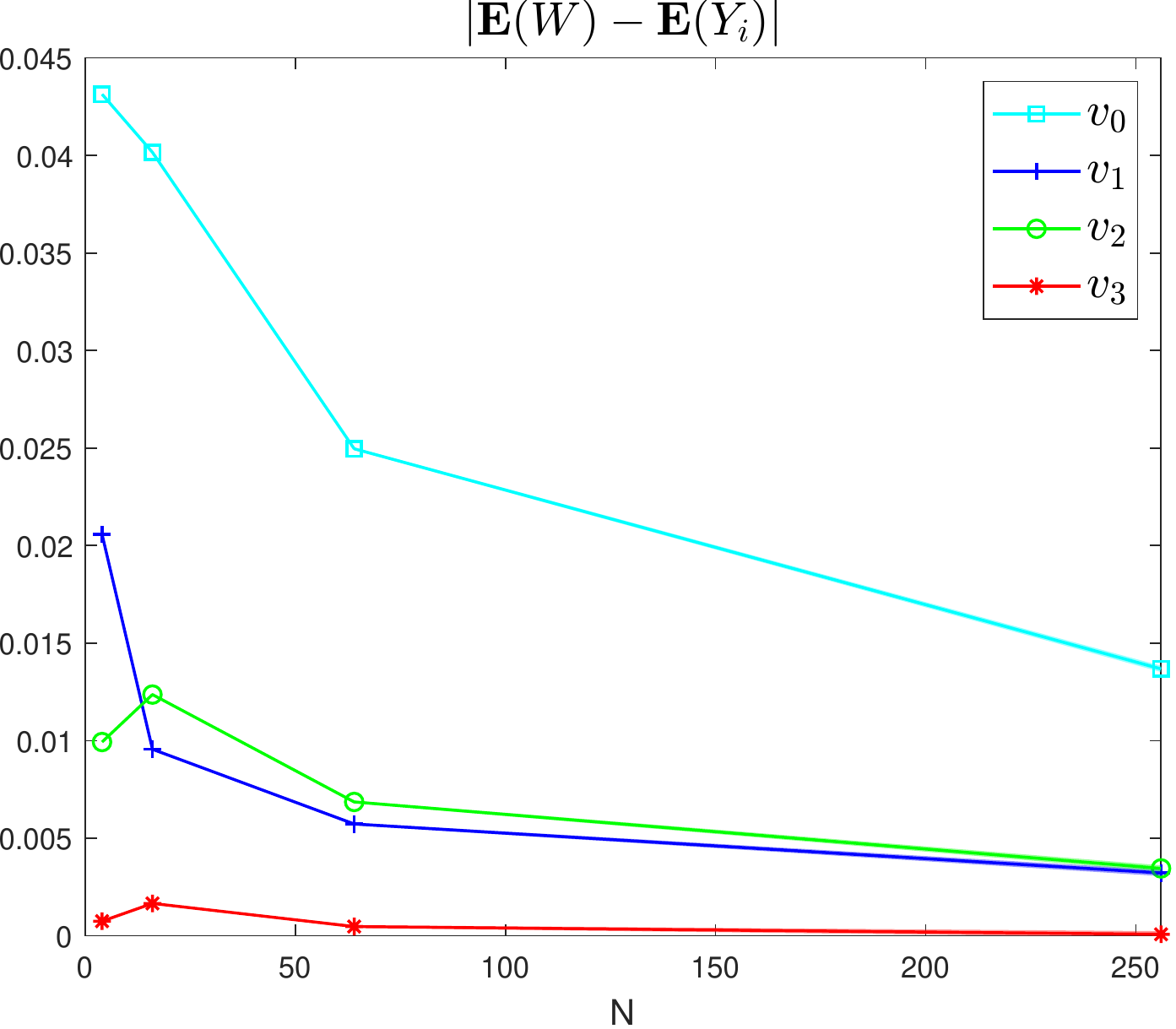}
%
\end{subfigure}%
\begin{subfigure}[t]{0.45\textwidth}
\centering
\setlength\tabcolsep{3pt}
\renewcommand{\arraystretch}{0.5}
\vspace{-3.5cm}
\begin{tabular}{rc | c | c  }
 $N$ & $\beta$ & $\E W$ &$95\% \text{CI for } \big| \E (W) - \E (Y_3) \big|$  \\
\hline
4    & 1  & -0.933      & $[0.0007, 0.0008] $ \\
16   & 1 & -0.865   & $[0.0016, 0.0017]$  \\
64  & 1 &  -0.823 &$[0.0004, 0.0005]$    \\
256 &  1 & -0.801 & $[0,0.0002]$  \\~\\
\end{tabular}
\end{subfigure}
\caption{$\beta = 1$. $Y_n$ corresponds to the $v_n$ approximation. }
 \label{fig:hospmoments}
\end{figure}

\begin{figure}
\centering
	\includegraphics[width=3.5in]{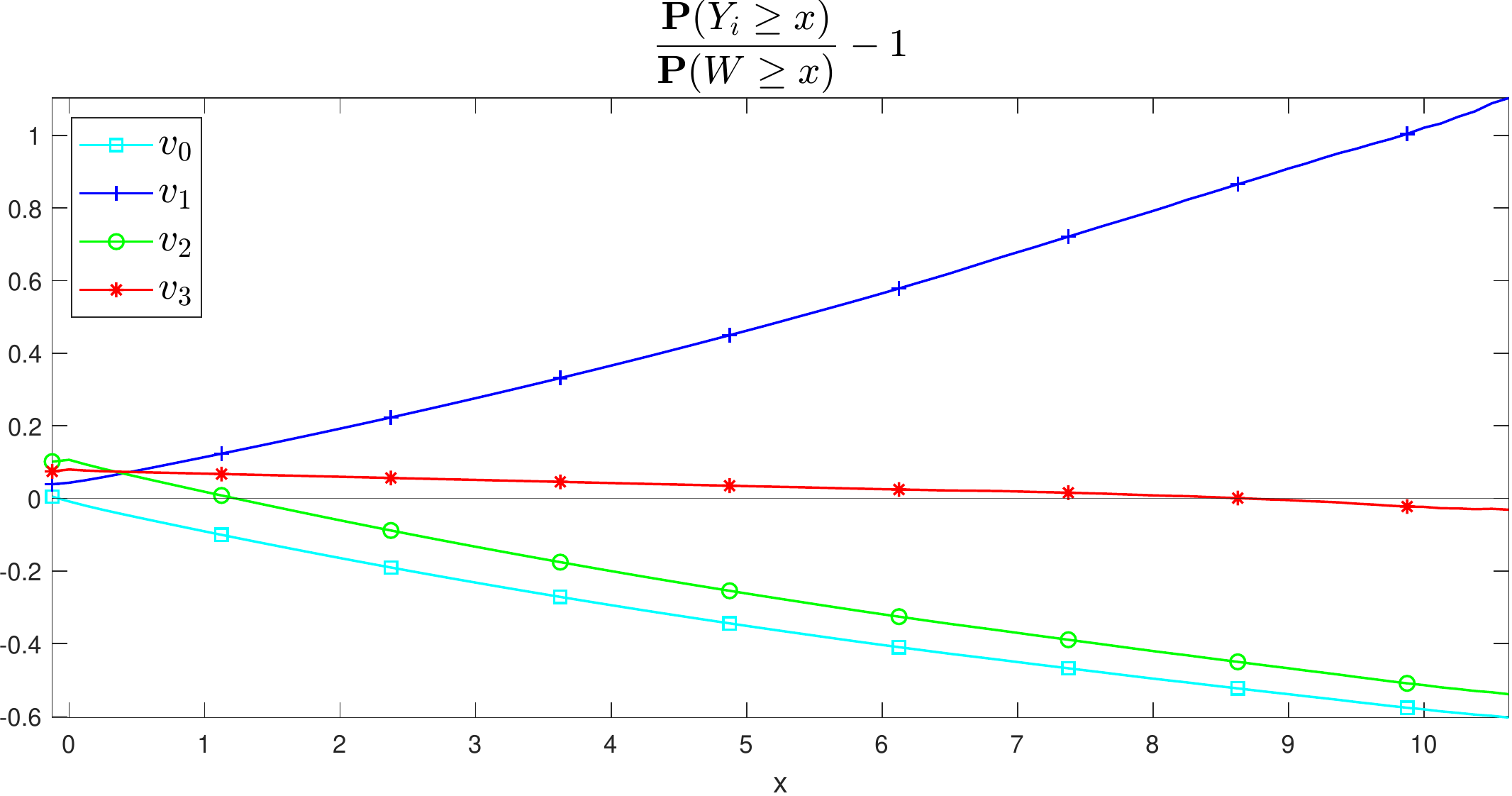}
	\caption{$N = 64$ and $\beta = 1$. $\Prob(W \geq x) \approx 10^{-6}$ for $x$  at the right endpoint of the $x$-axis.  }
	\label{fig:hosp64}
\end{figure}

In Figure~\ref{fig:hospmoments} we compare the   $v_n$ approximations for $\E (W)$ when $\beta = 1$ and $N \in  \{4,16,64,256\}$. The values of $\E (W)$ are estimated using a simulation; the width of the $95\%$ confidence intervals (CIs) is on the order of $10^{-4}$. Though we do not report them, the $v_0, v_1$, and $v_2$ approximation errors appear to decay at the rate of $1/\sqrt{N}$, but the $v_3$ error in the table appears to decay linearly in $N$.  When it comes to approximating the CCDF,  $v_3$ also outperforms the other approximations; see Figure~\ref{fig:hosp64} for an example when $N = 64$ and $\beta = 1$. Our findings were consistent for other values of $\beta$ and $N$. 
\vspace{-1cm}

%

\section{AR(1) Model}\label{fse5}
In this section we consider the first-order autoregressive model with random coefficient and general error distribution, which we refer to as the AR$(1)$ model. \cite{BlanGlyn2018} studied this model and   used an Edgeworth expansion to approximate its stationary distribution. Following the notation of \cite{BlanGlyn2018}, we compare the performance their expansion to our diffusion approximations.  

We encounter the issue that for large values of $x$, the untruncated $v_3(x)$ becomes sufficiently negative to make our usual solution of truncation from below perform poorly when approximating the tail of the distribution. We resolve this via a hybrid approximation that uses both $v_3(x)$ and $v_2(x)$ to construct a diffusion coefficient $\hat v_3(x)$ that combines the extra accuracy of $v_3(x)$ with the positivity of $v_2(x)$. 

To introduce the model, let $\{X_n, n\geq 1\}$, $\{Z_n, n\geq 1\}$ be two independent sequences of i.i.d.\ random variables. We assume that both $X_1$ and $Z_1$ are exponentially distributed with unit mean so that we can compare our approximations   to those of \cite{BlanGlyn2018}. Given $D_{0} \in \R$ and  $\alpha>0$, consider the DTMC $D=\{D_n, n\geq 0\}$ defined as 
\begin{align}
D_{n+1} = e^{-\alpha Z_{n+1}} D_n + X_{n+1} \label{eq:defar1}
\end{align}
 and let $D_\infty > 0$ denote the random variable having the stationary distribution of $D$. Using \eqref{eq:defar1}  we can see that $D_{\infty}$ is equal in distribution to $\sum_{k=0}^\infty X_k e^{-\alpha\sum_{j=0}^{k-1}Z_j}$, which is   the random variable studied in Section 4 of  \cite{BlanGlyn2018}. 

 We consider  $\tilde D = \{ \tilde D_n = \delta (D_n - R),\ n \geq 0\}$, where $\delta = \sqrt{\alpha}$ and,  to be consistent with \cite{BlanGlyn2018}, we choose  $R = 1/\alpha$. The asymptotic regime we consider is   $\alpha \to 0$, so   going forward we assume that $\alpha \in (0,1)$. 
It follows from \eqref{eq:defar1} that 
\[
\tilde D_{n+1} = e^{-\alpha Z_{n+1}} \tilde D_n + \delta\big(X_{n+1} + R(e^{-\alpha Z_{n+1}} - 1)\big).
\]
Let $W = \delta(D_{\infty}-R)$   and $W' = e^{-\alpha Z}W + \delta\big(X +  R(e^{-\alpha Z} - 1)\big)$, where $(X, Z)$ is an independent copy of $(X_1, Z_1)$, which is also independent of $W$. Since $D_{\infty} > 0$, the support of $W$ is $\mathcal{W} = (-1/\sqrt{\alpha}, \infty)$, which grows as $\alpha \to 0$. 
Stationarity implies that $\E f(W')-\E f(W)=0$ provided $\E \abs{f(W)} < \infty$. 
Note that the one-step jump size 
\begin{align*}
\Delta = W' - W = W\big(e^{-\alpha Z}-1\big) + \delta\big(X +  R(e^{-\alpha Z} - 1)\big)   =   \delta \big(D_{\infty}(e^{-\alpha Z} - 1) + X\big)
\end{align*}
does not depend on the choice of $R$. To present our diffusion approximations, we need expressions for $\E(\Delta^{k} | W = x)$. The following lemma  is proved in Section~\ref{app:ar1proof}.
 \begin{lemma}\label{lem:ar1}
 Recall that $\delta = \sqrt{\alpha}$. For any $k \geq 1$, 
 \begin{align*}
 \E(\Delta^{k} | D_{\infty} = d) =  \delta^{k} k! \bigg(1 + \sum_{i=1}^{k} (-1)^{i}d^{i} \prod_{j=1}^{i} \frac{ \alpha }{1 + j \alpha}  \bigg), \quad d > 0.
\end{align*}
\end{lemma} 
The relationship between $D_{\infty}$ and $W$ implies that  
\begin{align*}
  \E(\Delta^{k} | W = x) =  \E(\Delta^{k} | D_{\infty} = x/\delta + R) =&\  \delta^{k} k! \bigg(1 + \sum_{i=1}^{k} (-1)^{i}\Big( x\sqrt{\alpha}  +  1 \Big)^{i} \prod_{j=1}^{i} \frac{1}{1 + j \alpha}  \bigg),
\end{align*}
for $x \in \mathcal{W}$, where we used the facts that $\delta = \sqrt{\alpha}$ and  $R = 1/\alpha$ in the second equality. Extending $\E(\Delta^{k} | W = x)$ to all $x\in \R$ in the obvious way, we now state the $v_0$, $v_1$, and $v_2$ approximations, whose forms are all standard. Namely, $v_{2}(x)$ follows from \eqref{eq:v2},  $v_1(x) =   \E(\Delta^{2} | W = x)/2$, and $v_0 = \E(\Delta^{2} | W = x^{\ast})$, where $x^{\ast} = 1+1/\alpha$ solves $\E(\Delta | W = x) = 0$. To present $v_3(x)$, let us note that $\E(\Delta^{k} | W = x)$ takes the form  $\delta^k p_{k}(x)$ for some  degree-$k$ polynomial $p_{k}(x)$; we omit the dependence on $\alpha$ to ease notation. Given a truncation level $\eta > 0$, we let $v_3(x) = (\underline{v}_3(x) \vee \eta)$, where
\begin{align*}
\underline{v}_3(x) =&\ \delta^2 \Big(\frac{p_{2}(x)}{2}-\frac{p_{1}(x)\bar p_3(x)}{ \underline{p}_2(x)}-\delta \underline{p}_2(x)\Big(\frac{\bar p_3(x)}{\underline{p}_2(x)}\Big)'\Big), \\ 
\bar p_3(x) =&\   \frac{1}{6}  \Big( p_3(x) -   \frac{p_1(x)p_4(x)}{ 2 p_2(x)} - \frac{1}{4} \delta  p_2(x)\Big(\frac{p_4(x)}{p_2(x)}\Big)'  \Big),\\
\underline p_2(x) =&\   \Big(\frac{p_2(x)}{2}-\frac{p_1(x)p_3(x)}{3p_2(x)}-\frac{p_2(x)}{6}\Big(\frac{p_3(x)}{p_2(x)}\Big)'\Big).
\end{align*}
We derive $v_3(x)$ by successively applying the ``trick'' we used in Section~\ref{sec:v2def} to derive the $v_2(x)$ approximation in order to gain access to higher order terms in the Taylor expansion. The details are left to  Section~\ref{app:ar1proof} of the electronic companion. Before presenting our numerical results, we discuss how to modify the $v_3$ approximation to overcome the issue that $\underline{v}_3(x)$ becomes negative for large values of $x$.

%

\subsection{Hybrid Approximation} \label{sec:hybrid}
\begin{figure}
	\includegraphics[width=2.5in]{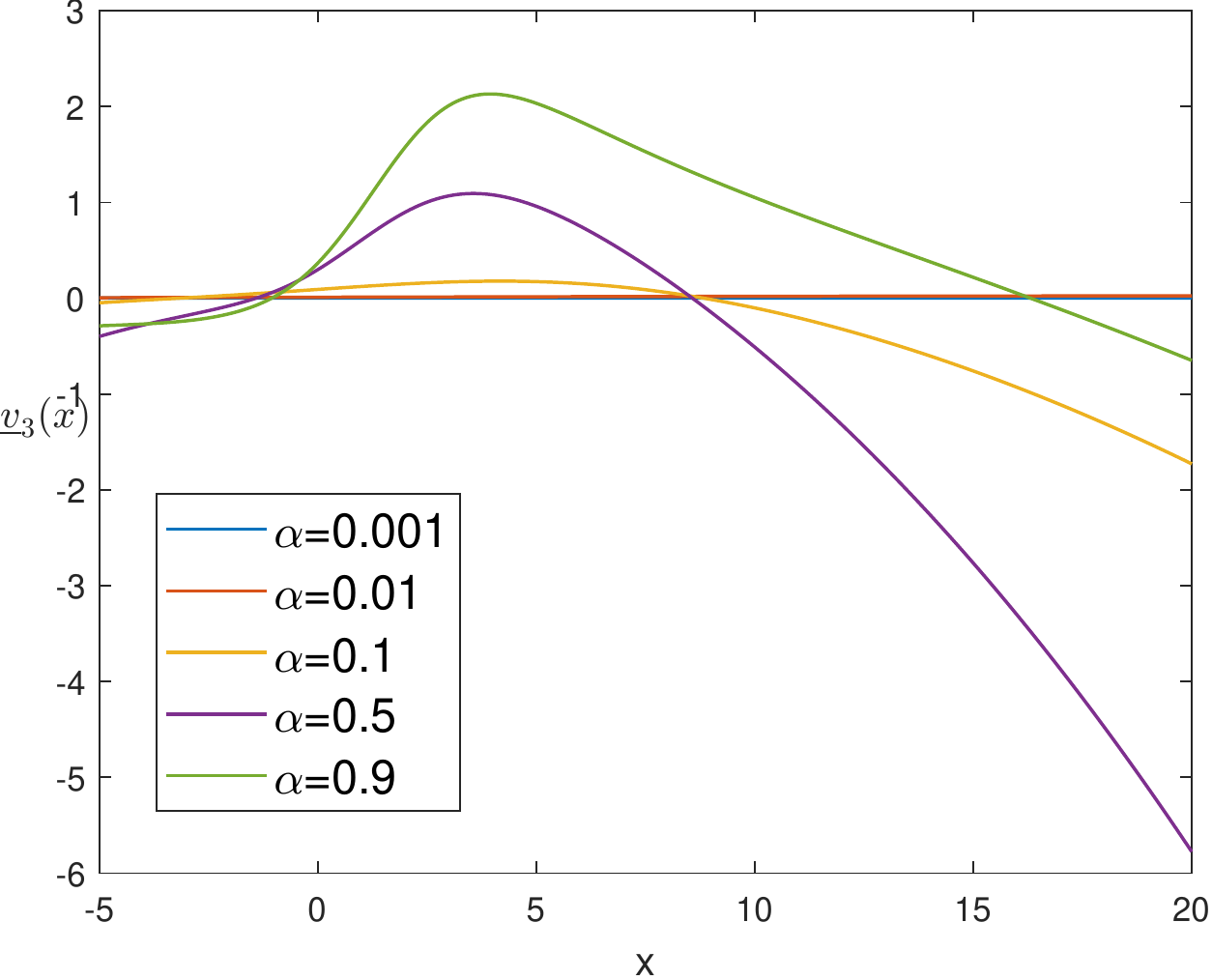}
	\centering
	\caption{ When $\alpha = 0.001$ and $0.01$, $\underline{v}_3(x)$  is nonnegative at all points plotted.  \label{fig:v3plots}}
\end{figure}

 Figure~\ref{fig:v3plots} displays  $\underline{v}_3(x)$ for several values of $\alpha$. When $\alpha = 0.001$ or $0.01$, the plots of $\underline{v}_3(x)$  are very close to zero but remain nonnegative. However, $\underline{v}_3(x)$ is negative when $\alpha = 0.1, 0.5$, or $0.9$. The behavior of $\underline{v}_3(x)$ in the left tail is not as important  because the left boundary of the support of $W = \sqrt{\alpha} (D - 1/\alpha)$ is $-1/\sqrt{\alpha}$; we therefore ignore the negativity in the left part of $\underline v_3(x)$ for $x < 0$. The farther $\underline{v}_3(x)$ drops below zero, the worse we expect the truncated  $\underline{v}_3(x)$ to perform. For example, the plot with $\alpha = 0.9$ in Figure~\ref{fig:arbig} of  Section~\ref{sec:ar1numerical} shows that while $v_3$ performs well in regions  where $\underline{v}_3(x)>0$, it does not perform as well  when estimating $\Prob(W > x)$ for large $x$. 

To improve upon the $v_3$ approximation, we propose the hybrid approximation $\hat v_3(x) = \underline{v}_3(x) 1(x \leq K) + \underline{v}_2(x) 1(x > K)$.  The threshold $K$ is numerically chosen to equal the right-most point of intersection of $\underline{v}_2(x)$ and $\underline{v}_3(x)$. The idea is for $\hat v_3(x)$ to enjoy the increased accuracy of $\underline{v}_3(x)$ in the center with the performance of $\underline{v}_2(x)$ far in the tail. 
We expect $\hat v_3(x)$ to outperform a truncated $\underline{v}_3(x)$  when $\underline{v}_3(x)$ drops far below zero; e.g., when $\alpha = 0.9$. If $\underline{v}_3(x)$ is nonnegative,  we expect little  benefit from $\hat v_3(x)$; e.g., when $\alpha = 0.001$. Our expectations are consistent with our numerical findings in Section~\ref{sec:ar1numerical}. 
Lastly, we remark on what can be done in the case when both $\underline{v}_2(x)$ and  $\underline{v}_3(x)$ are negative in the same region: instead of falling back on $\underline{v}_2(x)$, we can combine $\underline{v}_3(x)$  with $v_1(x)$, which is always positive because  $v_1(x) = \frac{1}{2} \E (\Delta^{2} | W=x) > 0 $ for all $x \in \mathcal{W}$.

\subsection{Numerical Results} \label{sec:ar1numerical}
\begin{figure}
\centering
\begin{subfigure}[t]{0.45\textwidth}
\centering
    \includegraphics[width=.8\linewidth]{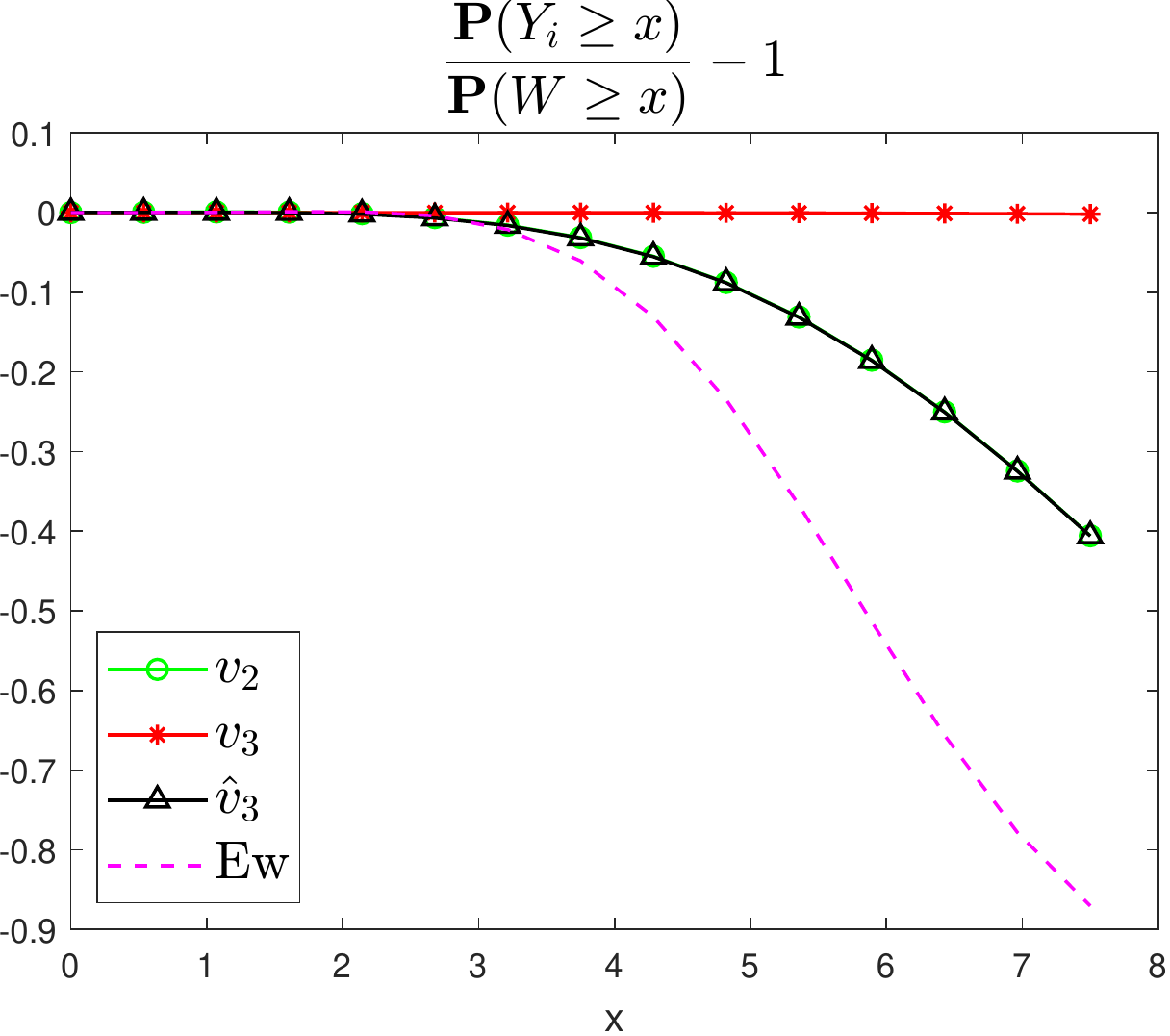}
    \vspace{-.3cm}
    \caption{\footnotesize{$\alpha = 0.001$}}
\end{subfigure}%
\begin{subfigure}[t]{0.45\textwidth}
\centering
    \includegraphics[width=.8\linewidth]{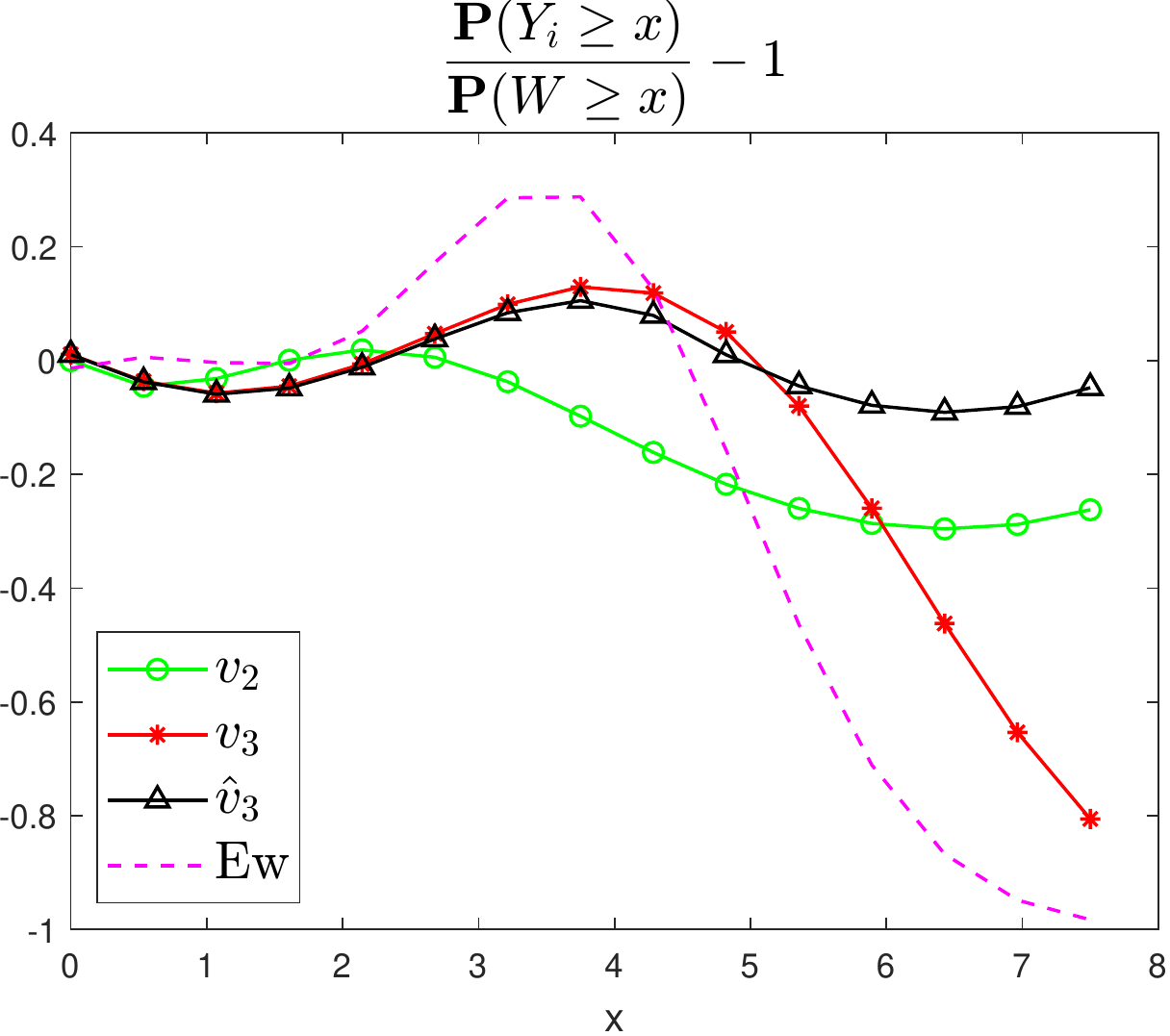}
    \vspace{-.3cm}
    \caption{\footnotesize{$\alpha = 0.9$}}
\end{subfigure}%
\caption{The plots exclude $v_0$ and $v_1$, the worst performing approximations.}
  \label{fig:arbig}
\end{figure}

It is well known that Edgeworth expansions, obtained for the probability distribution at a particular point, can suffer from two issues: (1) they may not be a proper probability distribution function, and (2) they may not be sufficiently accurate in the tails. 
Our diffusions approximate the entire distribution.
We  compare the quality of the $v_n$ and $\hat v_3$ approximations to the Edgeworth expansion of \cite{BlanGlyn2018}. Figure~\ref{fig:arbig} displays the relative error of approximating the CCDF of $W$ for different values of $\alpha$ and contains two plots: one where $\alpha$ is close to zero and one where $\alpha$ is far from zero. In the latter plot, the hybrid approximation is the best performer because $\alpha = 0.9$ and $\underline{v}_3(x)$ is negative in Figure~\ref{fig:v3plots},  whereas in the former plot $v_3$ is the best performer because $\alpha = 0.001$ and $\underline{v}_3(x)$ is nonnegative in Figure~\ref{fig:v3plots}. In addition to estimating the CCDF of $W$,  Table~\ref{tab:ar1log} compares the performance of our approximations when estimating the expectation of a smooth test function like $\E \log(W + \delta R) = \E \log(\alpha D_{\infty})$.

\begin{table}
  \begin{center}
\renewcommand{\arraystretch}{0.5}
   \begin{tabular}{|c|c|c|c|c|c|  }
   \hline 
 & $\alpha = 0.64$ & $\alpha = 0.32$  &$\alpha = 0.16$ & $\alpha = 0.08$& $\alpha = 0.04$\\ 
\hline
 $\abs{\E f(Y_0) - \E f(W)}$ &   0.095  & 0.039  &  0.011  & 0.002  & $3.6 \times 10^{-4}$  \\ \hline
 $\abs{\E f(Y_1) - \E f(W)}$      &   0.104  & 0.037  &  0.008  &  $9.4 \times 10^{-4}$ &  $1.0\times 10^{-4}$  \\ \hline
 $\abs{\E f(Y_2) - \E f(W)}$     &   0.034  & 0.011  &  0.003 & $4.7 \times 10^{-4}$  & $7.4 \times 10^{-5}$    \\ \hline
$\abs{\E f(Y_3) - \E f(W)}$    &   0.0201  & 0.005  & $7.9 \times 10^{-4}$  & $8.9 \times 10^{-5}$ &$7.8\times 10^{-6}$  \\ \hline
 $\abs{\E f(\hat{Y}_3) - \E f(W)}$      &   0.0194  & 0.005  &  $7.8 \times 10^{-4}$  & $8.9 \times 10^{-5}$ & $7.8 \times 10^{-6}$ \\ \hline
 $\abs{\E f( {Y}_e) - \E f(W)}$     &   0.148  & 0.053  &  0.009  & $7.3 \times 10^{-4}$  & $6.3 \times 10^{-4}$  \\  \hline
 $ \E f(W) $     &   0.510  & 0.721  &  0.994  & 1.302  &  1.629 \\  \hline
  \end{tabular} 
  \end{center}
  \caption{ $f(W) = \log(W + \delta R)$. The random variable $Y_n$ corresponds to the $v_n$-approximation, $\hat Y_3$ corresponds to the $\hat v_3$-approximation, and $Y_e$ corresponds to the Edgeworth expansion estimate.  \label{tab:ar1log}}
\end{table} 
%

\section{Conclusion}
We have outlined a general procedure to derive $v_n$ approximations for one-dimensional Markov chains. Although the expressions for $v_n(x)$  get more complicated as $n$ increases, the diffusion approximations remain computationally tractable. A natural question is how to extend this work to the multi-dimensional setting.

Another direction worth exploring relates to establishing theoretical guarantees for the approximations. The only results we have are for the $v_1$ error in Section~\ref{fse3}, a key ingredient of which are bounds on the derivatives to the solution of the Poisson equation, also known as Stein factor bounds; c.f.,  Lemma~\ref{lem:higherders} of the electronic companion. Since the Poisson equation depends on the diffusion coefficient,  we have to reestablish Stein factor bounds for each new $v_n(x)$; the difficulty of this grows with the complexity of the expression for $v_n(x)$. The  prelimit generator approach, recently proposed by \cite{Brav2022}, may offer a simpler avenue for theoretical guarantees because it uses Stein factor for the Markov chain instead of the diffusion.


\ECSwitch


\ECHead{Accompanying Proofs}

This e-companion contains the proofs of certain theoretical results in the paper. It is divided into three main sections. The first section is about the Erlang-C model,  and contains the proofs of Proposition~\ref{prop:lowerbound}, Theorem~\ref{thm:md-high}, and Theorem~\ref{thm:md-std}. The second  and third sections derive the $v_3$ approximation for the hospital model and AR(1) model, respectively.

\section{Companion for the Erlang-C Model}
\blue{To prepare for the arguments to come, let us recall the notation related to the Erlang-C model.  The Erlang-C model is defined by the customer arrival rate $\lambda >0$, the service rate $\mu > 0$, and the number of servers $n > 0$. Additional important quantities include 
\begin{align*}
R=\frac{\lambda}{\mu} <n,\ \beta=\frac{n-R}{\sqrt{R}}>0,  \text{ and }  \delta=\frac{1}{\sqrt{R}}.
\end{align*}
We study $W$, which has the stationary distribution of the CTMC $\{\tilde X(t) = \delta (X(t) - R)\}$, where $X(t)$ is the number of customers in the system at time $t \geq 0$. Equation \eqref{eq:gxbar} states that
\ben{ \label{eq:gxbarf}
\E G_{\tilde X} f(W) = \E \Big[ \lambda (f(W+\delta)-f(W))+\mu\big[  (W/\delta+ R) \wedge n  \big] (f(W-\delta)-f(W))   \Big]=0
}
for all $f(x)$ satisfying $\E \abs{f(W)} < \infty$. We also note that the support of $W$ is
\begin{align}
 \mathcal{W} =  \{-\sqrt{R}, -\sqrt{R}+\delta,  -\sqrt{R}+2\delta,\dots\}. \label{f29}
\end{align}
To define $v_0(x)$ and $v_1(x)$, we recall from \eqref{eq:tb} and \eqref{eq:ta} that for $x \in \mathcal{W}$, 
\begin{align}
b(x) =&\  \delta \big(\lambda -\mu \big[(x/\delta+R)\wedge n\big]\big) \quad \text{ and } \quad a(x)=  \delta^2 \big(\lambda+ \mu \blue{[(x/\delta+R)\wedge n]} \big),  \label{eq:tabeff}
\end{align}
and from \eqref{eq:bdef}  and \eqref{eq:adef} that the extensions of these to $\R$ are
\begin{align}
b(x) =&\  -(\mu x \wedge \mu \beta) \quad \text{ and } \quad a(x)=    2\mu - \delta b(  -\sqrt{R}\vee x), \quad x \in \R. \label{eq:adeff}
\end{align}
We define  $v_1(x) = \frac{1}{2} a(x)$ and $v_0(x) = v_0 = v_1(0) = \mu$,  and for $n \in \{0,1\}$ we define the $v_n$ approximation to be the random variable $Y_n$ with density
\begin{equation}
  \label{eq:stddenf}
   \frac{\kappa}{v_n(x)} \exp\Big({\int_0^x \frac{b(y)}{v_n(y)}dy}\Big), \quad x \in \R,
\end{equation}
where $\kappa > 0$ is a normalization constant that depends on $n$. Lastly, assuming that  $-z$ belongs to $ \mathcal{W}$ and setting   $f(x)=1(x\geq -z)$ in \eq{eq:gxbarf}, we get
\ben{\label{f26}
\lambda \P(W=-z-\delta)=\mu[(-z/\delta+R)\wedge n] \P(W=-z),
}
which are the flow-balance equations for the CTMC. 
}

\subsection{Proving Proposition~\ref{prop:lowerbound}}
\label{sec:lowerbound}
We repeat the statement of Proposition~\ref{prop:lowerbound} for convenience. 
\begin{proposition}
\label{prop:lowerboundec}
Assume that $n = R + \beta \sqrt{R}$ for some fixed $\beta > 0$. 
There exists a constant $C(\beta) > 0$ depending only on $\beta$ such that
\begin{align*}
 \big| \EE W - \EE  Y_0  \big| \geq \frac{C(\beta)}{\sqrt{R}}.
\end{align*}
\end{proposition}
We prove the proposition with the help of four auxiliary lemmas. The lemmas are proved at the end of this section after we prove Proposition~\ref{prop:lowerboundec}.  We use $C = C(\beta)>0$ to denote a constant that may change from line to line, but does not depend on anything other than $\beta$. 
\begin{lemma} \label{lem:momentequiv}
For any $\beta > 0$, 
\begin{align}
\beta \Prob(W \geq \beta) = -\E \big( W 1(W < \beta)\big)  \quad \text{ and } \quad \beta \Prob(Y_0 \geq \beta) = -\E \big( Y_0 1(Y_0 < \beta)\big). \label{eq:bizero}
\end{align} 
Consequently,
\begin{align*}
\E W = \E (W-\beta)^{+} \quad \text{ and } \quad \E Y_0 = \E (Y_0-\beta)^{+}.
\end{align*}
\end{lemma}
Lemma~\ref{lem:momentequiv} implies that
\begin{align*}
\abs{ \E W - \E Y_0} = \abs{ \E (W-\beta)^{+} - \E (Y_0-\beta)^{+}}.
\end{align*}
The next lemma rewrites the right-hand side above using the Poisson equation so that we can bound it from below.  
\begin{lemma}
\label{lem:expliciterror}
Fix $h \in \lipone$  and let $f_h(x)$ be the solution the the Poisson equation 
\begin{align}
b(x) f_h'(x) + v_0 f_h''(x) = \E h(Y_0) - h(x), \quad x \in \R. \label{eq:poissonv0}
\end{align}
Then  $f_h'''(x-) = \lim_{y \uparrow x} f_h'''(y)$ is defined for all $x \in \R$ and 
\begin{align*}
 \E h(W) - \E h(Y_0) =  - \frac{1}{2} \delta \E b(W) f_h''(W)+ \E \varepsilon(W),
\end{align*} 
where 
\begin{align*}
\varepsilon(W) =&\ \frac{1}{6}\delta^2 b(W) f_h'''(W-)+ \lambda(\varepsilon_{+}(W) + \varepsilon_{-}(W)) - \frac{1}{\delta} b(W) \varepsilon_{-}(W), \\
\varepsilon_{+}(W) =&\  \frac{1}{2}\int_{W}^{W+\delta} (W+\delta - y)^{2} (f_h'''(y) - f_h'''(W-)) dy, \\
\varepsilon_{-}(W) =&\  - \frac{1}{2}\int_{W-\delta}^{W} (y- (W-\delta))^{2} (f_h'''(y) - f_h'''(W-)) dy.
\end{align*}
\end{lemma}
Our plan is to show that for any $h \in \lipone$, the term $\E \varepsilon(W)$ vanishes at a rate of at least $\delta^2$. We then fix $h(x) = (x-\beta)^{+}$ and show that $\abs{\E b(W) f_h''(W)}$ can be bounded away from zero by a constant independent of $\delta$, which implies Proposition~\ref{prop:lowerboundec}. The following two lemmas are needed for this. The first one is for the upper bound on $\E \varepsilon(W)$, and the second is for the lower bound on $\abs{\E b(W) f_h''(W)}$.
\begin{lemma}
\label{lem:higherders}
Assume  that $n = R + \beta \sqrt{R}$ for some fixed $\beta > 0$. 
There exists $C = C(\beta) > 0$ depending only on $\beta$ such that for any $h \in \lipone$,
\begin{align*}
\abs{ f_h^{'''}(x-)} \leq  \frac{C}{\mu }  \quad \text{ and } \quad \E \abs{b(W)} \leq \mu \E \abs{W} \leq&\ \mu C , \quad x \in \R.
\end{align*}
Additionally, if $h(x) = (x-\beta)^{+}$, then for all  $x \neq \beta$, $f_h^{(4)}(x)$ exists and $\abs{f_h^{(4)}(x)} \leq \frac{C}{\mu }(1 + \abs{x})$.
\end{lemma}
\begin{lemma}
\label{lem:secondder}
If $h(x) = (x-\beta)^{+}$, then $f_h''(x) = \frac{1}{\mu \beta}$  for $x \geq \beta$ and
\begin{align*}
f_h''(x) =\frac{1}{\mu \beta}  \frac{1 + x e^{x^2/2} \int_{-\infty}^{x} e^{-y^2/2} dy }{1 +  \beta e^{\beta^2/2}  \int_{-\infty}^{\beta} e^{-y^2/2} dy   }	 , \quad x \leq \beta.
\end{align*}
\end{lemma}
Before proving Proposition~\ref{prop:lowerboundec}, let us remark that we restrict ourselves to $h(x) = (x-\beta)^{+}$ to keep the proof simple. Our arguments can likely be extended to work for other   $h(x)$  at the expense of added complexity.  
\startproof{Proof of Proposition~\ref{prop:lowerboundec}}
Lemma~\ref{lem:expliciterror} implies that
\begin{align*}
 \E h(W) - \E h(Y_0) =  - \frac{1}{2} \delta \E b(W) f_h''(W)+ \E \varepsilon(W).
\end{align*}
In the first part of the proof, we show that $\E \abs{\varepsilon(W)} \leq C \delta^2$. For convenience, we recall that
\begin{align*}
\varepsilon(W) =&\ \frac{1}{6}\delta^2 b(W) f_h'''(W-)+ \lambda(\varepsilon_{+}(W) + \varepsilon_{-}(W)) - \frac{1}{\delta} b(W) \varepsilon_{-}(W), \\
\varepsilon_{+}(W) =&\  \frac{1}{2}\int_{W}^{W+\delta} (W+\delta - y)^{2} (f_h'''(y) - f_h'''(W-)) dy, \\
\varepsilon_{-}(W) =&\  - \frac{1}{2}\int_{W-\delta}^{W} (y- (W-\delta))^{2} (f_h'''(y) - f_h'''(W-)) dy.
\end{align*}
Lemma~\ref{lem:higherders} immediately implies that $\E \abs{b(W) f_h'''(W-)} \leq C$. Next, we show that $\E  \abs{\varepsilon_{+}(W)} \leq \frac{C}{\mu} \delta^4 $. By considering the cases when $W = \beta$ and $W \neq \beta$ and using Lemma~\ref{lem:higherders}, we see that
\begin{align*}
\E  \abs{\varepsilon_{+}(W)}  \leq&\ \frac{1}{2}\Prob(W = \beta)\bigg|\int_{\beta}^{\beta+\delta} (\beta+\delta - y)^{2} (\abs{f_h'''(y)} + \abs{f_h'''(\beta-)}) dy\bigg| \\
&+ \frac{1}{2}\E \Bigg[1(W \neq \beta)\bigg|\int_{W}^{W+\delta} (W+\delta - y)^{2}\int_{W}^{y} \abs{f_h^{(4)}(u)} du  dy \bigg|\Bigg]\\
 \leq&\ \frac{C}{\mu}\Prob(W = \beta)\bigg|\int_{\beta}^{\beta+\delta} (\beta+\delta - y)^{2}  dy\bigg|\\
 &+\frac{C}{\mu}\E \Bigg[1(W \neq \beta) \bigg|\int_{W}^{W+\delta} \delta (1+\abs{W}+\delta) (W+\delta - y)^{2}  dy \bigg|\Bigg]\\
 \leq&\ \frac{C}{\mu} \delta^3 \Prob(W = \beta) + \frac{C}{\mu} \delta^4.
\end{align*}
This argument   can be repeated to show $\E \abs{\varepsilon_{-}(W)} \leq \frac{C}{\mu} \delta^4$.
It was shown in (3.29) of \cite{Brav2017} that $\Prob(W = \beta)\leq C \delta$, but for completeness we repeat the argument at the end of the proof.  Combining these results, we arrive at 
\begin{align*}
\E\abs{\varepsilon(W)} \leq C \delta^2 + \lambda \frac{C}{\mu} \delta^4 + \E \abs{b(W)} \frac{C}{\mu} \delta^3 \leq C \delta^2,
\end{align*}
where in the last inequality we use the bound on $\E \abs{b(W)}$ from Lemma~\ref{lem:higherders} and the fact that $\delta^2 = 1/R = \mu/\lambda$. We now show that $\abs{\E b(W) f_h''(W)} \geq C$. Combining the fact that $b(x) =  -(\mu x \wedge \mu \beta)$ with the form  of $f_h''(x)$ from Lemma~\ref{lem:secondder}, we have 
\begin{align*}
\E b(W) f_h''(W) =&\ \frac{-\mu \beta}{\mu \beta} \Prob(W \geq \beta) - \frac{\mu}{\mu \beta} \frac{\E\bigg[  W\Big(1+ W e^{W^2/2}\int_{-\infty}^{W} e^{-y^2/2} dy     \Big) 1(W < \beta) \bigg]}{1 +  \beta e^{\beta^2/2}  \int_{-\infty}^{\beta} e^{-y^2/2} dy   } \\
\leq&\ -\Prob(W \geq \beta) - \frac{1}{\beta} \frac{\E\big(  W  1(W < \beta) \big)}{1 +  \beta e^{\beta^2/2}  \int_{-\infty}^{\beta} e^{-y^2/2} dy   }.
\end{align*}
From  Lemma~\ref{lem:momentequiv} we know that  $\E \big( W 1(W < \beta)\big) = -\beta\Prob(W \geq \beta)$, so
\begin{align*}
\E b(W) f_h''(W) =&\ -\Prob(W \geq \beta) + \frac{1}{ \beta} \frac{ \beta \Prob(W \geq \beta) -  \E\bigg[ \Big( W^2 e^{W^2/2}\int_{-\infty}^{W} e^{-y^2/2} dy     \Big) 1(W < \beta) \bigg]}{1 +  \beta e^{\beta^2/2}  \int_{-\infty}^{\beta} e^{-y^2/2} dy   } \\
\leq&\ \Prob(W \geq \beta)\Big( -1 + \frac{1}{1 + \beta e^{\beta^2/2} \int_{-\infty}^{\beta} e^{-y^2/2} dy}  \Big)\\
\leq&\ -C \Prob(W \geq \beta).
\end{align*}
Proposition 1 of \cite{HalfWhit1981} tells us  that $\Prob(W \geq \beta)$ converges to a positive constant (depending on $\beta$) as $R \to \infty$. This implies the lower bound on $\abs{\E b(W) f_h''(W)}$. Lastly, we prove that $\Prob(W = \beta)\leq C \delta$. Let $\phi_0(x)$ be the density of $Y_0$. Since $W$ is grid valued, for any $z \in (0,\delta)$ we have
\begin{align*}
\Prob(W = \beta) =&\ \Prob( \beta - z \leq W \leq \beta + z) \\
\leq&\ \Prob( \beta - z \leq Y_{0} \leq \beta + z) + \abs{\Prob( \beta - z \leq W \leq \beta + z) - \Prob( \beta - z \leq Y_0 \leq \beta + z)}\\
\leq&\ 2z \sup_{x \in \R} \phi_0(x) + 2\sup_{x \in \R} \abs{\Prob(W \leq x) - \Prob(Y_0 \leq x)}.
\end{align*}
To reach the desired conclusion, we use Lemma 7 and Theorem 3 of \cite{BravDaiFeng2016}. The former says that $\phi_0(x) \leq \sqrt{2/\pi}$, while the latter result says that $\sup_{x \in \R} \abs{\Prob(W \leq x) - \Prob(Y_0 \leq x)} \leq C \delta$. 
\finishproof
 

\subsubsection{Proof of Lemma~\ref{lem:momentequiv}.}
\startproof{Proof of Lemma~\ref{lem:momentequiv} }
If $f(x) = x$, then $\E G_{\tilde X} f(W) = 0$ because $\E |W| < \infty$; see Lemma~2 of \cite{BravDaiFeng2016} for a proof of the latter fact. It follows from \eqref{eq:gxbarf} that 
\begin{align*}
\lambda  - \E \Big( \mu\big[  (W/\delta+ R) \wedge n  \big] \Big) = \E b(W) = - \mu \E(W \wedge \beta) = 0,
\end{align*}
implying that  
\begin{align}
-\beta \Prob(W \geq \beta) = \E \big( W 1(W < \beta)\big). \label{eq:leminterm}
\end{align} 
Adding $ \E \big( W 1(W \geq \beta)\big)$ to both sides  proves that $\E W = \E (W-\beta)^{+}$.    The claim  about $Y_0$ follows similarly. The density of $Y_0$ is given by \eqref{eq:stddenf}, so 
\begin{align*}
\E b(Y_0) = \kappa \int_{-\infty}^{\infty} \frac{b(x)}{\mu} \exp\Big(\int_0^x \frac{b(y)}{\mu} dy \Big) dx = 0.
\end{align*}
The last equality follows from integration by parts and the fact that $\lim_{x\to \pm \infty} \exp\big(\int_0^x \frac{b(y)}{\mu} dy \big)  = 0$. Therefore
\begin{align*}
\frac{1}{\mu} \E b(Y_0) = -\beta \Prob(Y_0 \geq \beta) - \E \big(Y_0 1(Y_0 < \beta)\big) = 0 
\end{align*}
and $\E Y_0 = \E (Y_0 -\beta)^{+}$.
\finishproof

\subsubsection{Proof of Lemma~\ref{lem:expliciterror}.}
\startproof{Proof of Lemma~\ref{lem:expliciterror}}
First, note that $f_h''(x) = -\frac{b(x)}{v_0} f_h'(x) + \frac{1}{v_0}( \E h(Y_0) - h(x))$ is differentiable almost everywhere because $f_h'(x)$ is continuously differentiable and both $b(x)$ and $h(x)$ are Lipschitz functions. The former statement follows, for example, from (B.1) of \cite{BravDaiFeng2016}. Therefore, $f_h'''(x)$ exists almost everywhere.  Now \eqref{eq:gxbarf} implies that   
\begin{align*}
\E G_{\tilde X} f_h(W) = \E \Big[ \lambda (f_h(W+\delta)-f_h(W))+\mu\big[  (W/\delta+ R) \wedge n  \big] (f_h(W-\delta)-f_h(W))   \Big]=0,
\end{align*}
provided $\E \abs{f_h(W)} < \infty$. The integrability of $f_h(W)$ has already been established in \cite{BravDaiFeng2016}; see Lemma 1 and Remark 2 there. Since $f_h'''(x)$ does not exist everywhere (for instance at $x = \beta$), to perform Taylor expansion we need to use the integral form of the remainder term. We claim that 
\begin{align}
f_h(x+\delta) - f_h(x) =&\ \delta f_h'(x) + \frac{1}{2}\delta^2 f_h''(x) + \frac{1}{6}\delta^{3} f_h'''(x-) + \varepsilon_+(x), \notag \\
f_h(x-\delta) - f_h(x) =&\ -\delta f_h'(x) + \frac{1}{2}\delta^2 f_h''(x) - \frac{1}{6}\delta^{3} f_h'''(x-) + \varepsilon_{-}(x).  \label{eq:vareps}
\end{align}
To verify the claim, note that
\begin{align*}
f_h(x+\delta) - f_h(x) = \int_{x}^{x+\delta} f_h'(y) dy =&\ \delta f_h'(x) + \int_{x}^{x+\delta}   (f_h'(y)-f_h'(x))  dy \\
=&\ \delta f_h'(x) + \int_{x}^{x+\delta} \int_{x}^{y}  f_h''(u) du dy\\
=&\ \delta f_h'(x) + \int_{x}^{x+\delta} (x+\delta - u) f_h''(u) du.
\end{align*}
A similar treatment of $f_h(x+\delta) - f_h(x)$ yields   \eqref{eq:vareps}. Letting $s(W) = \mu\big[  (W/\delta+ R) \wedge n  \big]$, we therefore have
\begin{align*}
 \E G_{\tilde X} f_h(W) =&\ \E \Big[ \delta (\lambda - s(W)) f_h'(W) + \frac{1}{2}\delta^2 (\lambda + s(W))  f_h''(W) \\
&+ \frac{1}{6}\delta^3 (\lambda - s(W)) f_h'''(W-) + \lambda \varepsilon_{+}(W) + s(W) \varepsilon_{-}(W)\Big].
\end{align*}
We know from \eqref{eq:tabeff} that $\delta(\lambda - s(W)) = b(W)$ and $\delta^2(\lambda + s(W)) = a(W)$, and consequently $s(W) = \lambda - b(W) / \delta$. Therefore,  
\begin{align*}
& \E G_{\tilde X} f_h(W) \\
 =&\ \E \Big[b(W)f_h'(W) + \frac{1}{2} a(W) f_h''(W) + \frac{1}{6}\delta^2 b(W) f_h'''(W-) + \lambda (\varepsilon_{+}(W) + \varepsilon_{-}(W)) - \frac{1}{\delta} b(W) \varepsilon_{-}(W)\Big] = 0.
\end{align*}
Taking expected values in the Poisson equation \eqref{eq:poissonv0}, we get
\begin{align*}
\E h(Y_0) - \E h(W) =&\   \E \Big[ b(W) f_h'(W) + v_0 f_h''(W)\Big] - \E G_{\tilde X} f_h(W) \\ 
=&\ - \frac{1}{2}\E \big(a(W) - a(0)\big) f_h''(W) - \E \varepsilon(W).
\end{align*}
We conclude by noting that $a(W) - a(0) = -\delta b(W)$.
\finishproof

\subsubsection{Proof of Lemma~\ref{lem:higherders}.}
\startproof{Proof of Lemma~\ref{lem:higherders}}
 To bound $\abs{f^{(4)}(x)}$, first note that the  Poisson equation \eqref{eq:poissonv0} implies that
\begin{align*}
f_h'''(x) =   -\frac{b(x)}{v_0} f_h''(x)  -\frac{b'(x)}{v_0} f_h'(x)   -\frac{1}{v_0} h'(x).
\end{align*}
Since $b(x)$ and $h(x)$ are piece-wise linear with a kink at $x = \beta$, the derivative above exists for all $x \neq \beta$.  Differentiating again, we get
\begin{align*}
f_h^{(4)}(x) =   -\frac{b(x)}{v_0} f_h'''(x) -\frac{b'(x)}{v_0} f_h''(x)  -\frac{b'(x)}{v_0} f_h''(x)  , \quad x \neq \beta.
\end{align*}
We conclude that $|f_h^{(4)}(x)| \leq (C/\mu)(1+\abs{x})$ for $x \neq \beta$ because $v_0 = \mu$,  $\abs{b(x)} \leq \mu \abs{x}$, $\abs{b'(x)} \leq \mu$, $\abs{f_h'''(x)} \leq C/\mu$, and $\abs{f_h''(x)} \leq C/\mu $, where the last two inequalities follow from Lemma 3 of \cite{BravDaiFeng2016}. To conclude the proof, we note that $\E \abs{b(W)} \leq \mu\E |W| \leq \mu C$, where the last inequality follows from
Lemma 2 of \cite{BravDaiFeng2016}, which tell us  that $\E \abs{W} \leq C $. 
\finishproof

\subsubsection{Proof of Lemma~\ref{lem:secondder}.}

\startproof{Proof of Lemma~\ref{lem:secondder}}
First assume that $x \geq \beta$. In (B.8) of \cite{BravDaiFeng2016} it is shown that 
\begin{align*}
f_h''(x) =&\  e^{-\int_{0}^{x} \frac{b(u)}{v_0} du}\int_{x}^{\infty} \frac{1}{\mu}\big( h'(y) + f_h'(y) b'(y)\big) e^{\int_{0}^{y} \frac{b(u)}{v_0} du} dy = e^{-\int_{\beta}^{x} \frac{b(u)}{v_0} du}\int_{x}^{\infty} \frac{1}{\mu}  e^{\int_{\beta}^{y} \frac{b(u)}{v_0} du} dy,
\end{align*}
where in the second equality we use $b'(x) = 0$ and $h'(x) = 1$ for $x \geq \beta$. Since $b(x)/v_0 = -\beta$ for $x \geq \beta$, 
\begin{align*}
f_h''(x) = e^{\beta(x-\beta)} \int_{x}^{\infty} \frac{1}{\mu } e^{-\beta(y-\beta)} dy = \frac{1}{\mu \beta}.
\end{align*}
Now suppose that $x \leq \beta$. The Poisson equation \eqref{eq:poissonv0} and the fact that $h(x) = 0$ imply that
\begin{align*}
f_h''(x) =&\ -\frac{b(x)}{v_0} f_h'(x) + \frac{1}{v_0} \E h(Y_0)  = x f_h'(x) + \frac{\E h(Y_0)}{\mu }.
\end{align*}
One can verify by differentiating that
\begin{align*}
f_h'(x) = e^{-\int_{0}^{x} \frac{b(u)}{v_0} du} \int_{-\infty}^{x} \frac{1}{v_0} \big(\E h(Y_0) - h(y) \big)  e^{\int_{0}^{y} \frac{b(u)}{v_0} du} dy = \frac{\E h(Y_0)}{\mu }e^{ \frac{1}{2}x^2 } \int_{-\infty}^{x}   e^{- \frac{1}{2} y^2} dy.
\end{align*}
 The first equality  appears as equation (B.1) in \cite{BravDaiFeng2016}. The second equality follows from the form of $b(x)$ in \eqref{eq:adeff} and the fact that $h(x) = 0$ for $x \leq \beta$.  Lastly, since the density of $Y_0$ is given by \eqref{eq:stddenf}, we have 
\begin{align*}
\E h(Y_0) = \frac{\int_{-\infty}^{\infty} h(y) e^{\int_{0}^{y} \frac{b(u)}{v_0} du} dy}{\int_{-\infty}^{\infty}   e^{\int_{0}^{y} \frac{b(u)}{v_0} du}dy} =&\ \frac{\int_{\beta}^{\infty} (y-\beta)^{+} e^{\int_{\beta}^{y} \frac{b(u)}{v_0} du}dy}{\int_{-\infty}^{\beta}   e^{\int_{\beta}^{y} \frac{b(u)}{v_0} du}dy + \int_{\beta}^{\infty}   e^{\int_{\beta}^{y} \frac{b(u)}{v_0} du}dy}\\
=&\  \frac{\int_{\beta}^{\infty} (y-\beta)  e^{-\beta(y-\beta)}dy}{\int_{-\infty}^{\beta}   e^{-\frac{1}{2}(y^2-\beta^2) }dy + \int_{\beta}^{\infty}   e^{-\beta(y-\beta)}dy} =  \frac{1/\beta^2 }{ e^{\frac{1}{2}\beta^2 } \int_{-\infty}^{\beta}   e^{-\frac{1}{2} y^2  }dy + 1/\beta}.
\end{align*}
This verifies the form of $f_h''(x)$ when $x \leq \beta$. 
\finishproof

\subsection{Proving Theorem~\ref{thm:md-high}}\label{fse8}
We first recall Theorem~\ref{thm:md-high}.
\begin{theorem}
\label{thm:md-highec}
Assume $n = R + \beta \sqrt{R}$ for some fixed $\beta > 0$. 
\blue{There exist positive constants $c_1$ and $C$ depending only on $\beta$ such that 
  \begin{align}
  & \left|\frac{\Prob(Y_1\geq z)}{\Prob(W\geq z)}-1\right|\leq \frac{C}{\sqrt{R}}\left(1+\frac{z}{\sqrt{R}}\right) \quad \text{ for } 0< z\leq c_1 R\ \text{and} \label{f13ec}\\
&\left|\frac{\Prob(Y_1\leq -z)}{\Prob(W\leq -z)}-1\right|\leq \frac{C}{\sqrt{R}}\left(1+z+\frac{z^4}{\sqrt{R}}\right),\ \text{ for } 0<z\leq \min\{c_1 R^{1/4}, R^{1/2}\} \label{f16ec}.
  \end{align}
} 
\end{theorem}
We  use \blue{$c_1,\ C,\ C_1, K$} to denote positive constants which may differ from line to line, but will only depend on $\beta$. We first prove \eqref{f16ec}, and then   \eqref{f13ec}.

\blue{ Note that if $R \leq K$ for some $K > 0$, then $n$ is also bounded because $\beta$ is fixed and $n = R + \beta \sqrt{R}$.  We argue that  \eqref{f16ec} holds trivially in such a case. Observe that for any $K > 0$, \magenta{because $\{W = -\sqrt{R}\}$ corresponds to an empty system,}
\begin{align*}
\inf_{\substack{0 < R \leq K \\ 0 < z \leq \sqrt{R}}} \Prob(W \leq -z) \geq \inf_{\substack{0 < R \leq K }} \Prob(W = -\sqrt{R})  \geq L(K,\beta) > 0,
\end{align*}
\magenta{where $L(K,\beta)$ is a positive constant depending only on $K$ and $\beta$. The second-last inequality is true because $\inf_{\substack{0 < R \leq K }} \Prob(W = -\sqrt{R})$ is at least as large as $\Prob(W = -\sqrt{R})$ when $R = K$, because for each $n$,  the probability of an empty system decreases in $R$.}  We can choose a sufficiently large $C$  (that depends on $K$) to ensure \eq{f16ec} trivially holds. Therefore, in the following, we assume $R\geq C_1$ for a sufficiently large $C_1$. Since \eqref{f16ec} requires $ 0 < z \magenta{\leq} c_1 R^{1/4}$,  we can also assume that  
\ben{\label{fEC4}
\delta = \frac{1}{\sqrt{R}} <\min\{1/2,\beta\}, \quad 0<z<\sqrt{R}-2,\quad \delta(z+1)<1/2.
} 
If not, we simply increase the value of $C_1$ and decrease the value of $c_1$ until \eqref{fEC4} holds. Without loss of generality let us therefore assume \eqref{fEC4} going forward. } Given $z \in \R$, we let $f_z(x)$ be the solution \magenta{(cf. \eq{f18})} to the Poisson equation 
\ben{\label{f17}
\frac{1}{2}a(x)f_z''(x)+b(x)f_z'(x)=\P(Y_1 \leq -z) - 1(x\leq -z), \quad x \in \R. 
}
The following object will be of use. Define\magenta{, for $W$ in its support \eq{f29},}
\begin{align}
K_W(y) = 
\begin{cases}
(\lambda - b(W)/\delta) (y+\delta) \geq 0, \quad y \in [-\delta, 0], \\
\lambda(\delta-y) \geq 0, \quad y \in [0,\delta].
\end{cases} \label{eq:kdef}
\end{align}
It can be checked that
\begin{align}
\int_{-\delta}^{0} K_W(y) dy =&\ \frac{1}{2} \delta^2 \lambda -\frac{1}{2}\delta b(W), \quad \int_{0}^{\delta} K_W(y) dy = \frac{1}{2}\delta^2 \lambda, \notag \\
\int_{-\delta}^{\delta} K_W(y) dy =&\ \frac{1}{2}a(W) = \mu - \frac{\delta}{2} b(W), \quad \text{ and } \quad \int_{-\delta}^\delta y  K_W(y) dy=\frac{\delta^2 b(W)}{6}. \label{MD:k3}
\end{align}
Our first result is an expression for $\Prob(Y_1 \leq -z) - \Prob(W \leq -z)$, and is proved in Section~\ref{sec:prooftaylormdhigh}.
\begin{lemma} \label{lem:taylormdhigh}
For any $z \in \R$,  
\begin{align}
&\Prob(Y_1 \leq -z) - \Prob(W \leq -z) \notag \\
  =&\ \E \bigg[\int_{-\delta}^{\delta} \Big(\frac{2b(W+y)}{a(W+y)}f_z'(W+y) - \frac{2b(W)}{a(W)}f_z'(W) \Big) K_W(y) dy  \bigg]  \notag \\
&\magenta{-}\E \bigg[\int_{-\delta}^{\delta} \Big( \frac{2}{a(W)} 1(W\leq -z) - \frac{2}{a(W+y)}1(W+y\leq -z) \Big) K_W(y) dy  \bigg] \notag  \\
&\magenta{-}\Prob(Y_1 \leq -z)\E \bigg[\int_{-\delta}^{\delta} \Big(  \frac{2}{a(W+y)} -\frac{2}{a(W)} \Big) K_W(y) dy  \bigg]. \label{f20}
\end{align} 
\end{lemma}
The bulk of the effort to prove \eqref{f16ec} comes   from the first term. The following lemma is proved in Section~\ref{fap4}.
\begin{lemma}\label{l4}
For $x \in \R$, define $r(x) = 2b(x)/a(x)$. There exists   constants $c_1,C,C_1 > 0$ depending only on $\beta$ such that for any $R\geq C_1$ and $0<z\leq c_1 R^{1/4}$ satisfying \eqref{fEC4},
\besn{\label{f21}
\left| \E  \bigg[ \int_{-\delta}^\delta \big(r(W+y)f_z'(W+y)-r(W)f_z'(W)\big)K_W(y)dy \bigg] \right|\leq&\ C \delta^2 (z\vee 1)^4  \P(Y_1\leq -z).
}
\end{lemma} 
\startproof{Proof of \eqref{f16ec}}
We first bound the  third term in \eqref{f20}. Using the form of $a(x)$ in  \eq{eq:adeff} and the assumption that $\delta < 1/2$ in \eqref{fEC4}, it is not hard to check that 
\ben{\label{fEC6}
\mu \leq a(x) \leq 2\mu+\delta\mu \beta \leq C \mu, \quad \text{ and } \quad \abs{a'(x)} \leq \delta \mu,
} 
from which it follows that 
\begin{align}
\frac{1}{a(x)} \leq 1/\mu  \quad \text{ and } \quad \abs{\frac{1}{a(x)} - \frac{1}{a(y)}} = \frac{\abs{a(y) - a(x)}}{a(y)a(x)} \leq \frac{\delta \abs{y-x}}{\mu}. \label{eq:alipsch}
\end{align}
Therefore, the third term in \eq{f20}  satisfies
\begin{align}
  \P(Y_1\leq -z)\left| \E \bigg[\int_{-\delta}^{\delta} \Big(  \frac{2}{a(W+y)} -\frac{2}{a(W)} \Big) K_W(y) dy \bigg] \right|  \leq&\ \P(Y_1\leq -z)\left| \E \bigg[\frac{2\delta^2}{\mu } \int_{-\delta}^{\delta} K_W(y) dy \bigg] \right| \notag \\
=&\ \P(Y_1\leq -z)\left| \E \bigg[\frac{2\delta^2}{\mu } \frac{1}{2}a(W)  \bigg] \right| \notag \\
 \leq&\ C\delta^2 \P(Y_1\leq -z). \label{f24}
\end{align}
The equality is due to \eqref{MD:k3}. The second term in \eqref{f20} is bounded similarly. Namely, 
\begin{align}
& \left| \E \bigg[\int_{-\delta}^{\delta} \Big( \frac{2}{a(W)} 1(W\leq -z) - \frac{2}{a(W+y)}1(W+y\leq -z) \Big) K_W(y) dy  \bigg] \right| \notag \\
\leq &\ \left| \E \bigg[\int_{-\delta}^{\delta} \Big( \frac{2}{a(W)}  - \frac{2}{a(W+y)}  \Big)1(W\leq -z) K_W(y) dy  \bigg] \right| \notag \\
&+\left| \E \bigg[\int_{-\delta}^{\delta}  \frac{2}{a(W+y)}\Big( 1(W\leq -z) -1(W+y\leq -z) \Big) K_W(y) dy  \bigg] \right| \notag \\
\leq &\ C \delta^2 \Prob(W \leq -z) + \frac{2}{\mu }  \E \bigg[\int_{-\delta}^{\delta}  \Big| 1(W\leq -z) -1(W+y\leq -z) \Big| K_W(y) dy  \bigg]. \label{f22}
\end{align}
Letting $-\tilde{z}$ denote the smallest value in the support of $\{W: W >  -z \}$ \magenta{(cf. \eq{f29})}, we have
\begin{align}
\frac{2}{\mu }  \E \bigg[\int_{-\delta}^{\delta}  \Big| 1(W\leq -z) -1(W+y\leq -z) \Big| K_W(y) dy  \bigg] \leq  C \P(W=-\tilde{z})+ C \P(W=-\tilde{z}-\delta). \label{eq:interm30}
\end{align}
At the end of this proof we argue that
\ben{\label{f23}
\P(W=-\tilde{z})+  \P(W=-\tilde{z}-\delta)\leq C\delta (z\vee 1) \P(W\leq -z).
}
Combining the bounds in \magenta{\eq{f24},} \eqref{f22}, \eqref{eq:interm30} and \eqref{f23} with Lemma~\ref{l4} yields 
\begin{align*}
|\P(Y_1\leq -z)-\P(W\leq -z)|\leq C\delta^2 (z\vee 1)^4 \P(Y_1\leq -z)+C\delta (z\vee 1) \P(W\leq -z).
\end{align*}
Dividing both sides by $P(W\leq -z)$, which is allowed because our assumption that $z < \sqrt{R} -2$  \magenta{in \eqref{fEC4}}  implies that $P(W\leq -z)\geq \P(W=-\sqrt{R})>0$, we arrive at
\ben{\label{fEC5}
\left|\frac{\P(Y_1\leq -z)}{\P(W\leq -z)}-1   \right|\leq  C \bigg( \delta^2 (z\vee 1)^4 \frac{\P(Y_1\leq -z)}{\P(W\leq -z)}+\delta(z\vee 1) \bigg).
}
Since we assumed $z \leq c_1 R^{1/4}$ and $R \geq C_1$, it follows that $\delta^2(z\vee 1)^4 \leq  c_1^4 \vee (1/\magenta{C_1})$, which can be made arbitrarily small (without affecting $C$ above) by decreasing $c_1$ and increasing $C_1$.  Therefore, without loss of generality  we can assume $C\delta^2(z\vee 1)^4<1/2$, so 
\begin{align*}
\frac{1}{2}\frac{\P(Y_1\leq -z)}{\P(W\leq -z)} \leq  (1-C\delta^2 (z\vee 1)^4)\frac{\P(Y_1\leq -z)}{\P(W\leq -z)}\leq 1+C\delta (z\vee 1) \leq C,
\end{align*}
where the second inequality is due to \eqref{fEC5} and the last inequality follows from \eqref{fEC4}. Combining the upper bound above with \eqref{fEC5} proves \eqref{f16ec}.
It remains to verify \eq{f23}. Equation \eqref{f26} implies that for $-y$ in the support of $W$,
\begin{align*}
\lambda \P(W=-y-\delta)=\mu[(-y/\delta+R)\wedge n] \P(W=-y).
\end{align*}
Since $\beta = \delta(n-R)$, the set $\{-y/\delta+R \leq n\}$ equals $\{-y \leq \beta\}$, so dividing both sides above by $\lambda$ and using the fact that $\mu/\lambda = 1/R$, it follows that 
\ben{\label{f27}
\P(W=- y-\delta)=(1-\delta y) \P(W=- y) \quad \text{ for  $-y\leq \beta$}.
} 
Recall that $-\tilde{z}$ denotes the smallest value in the support of $\{W: W > -z \}$, so $-\tilde z - \delta \leq -z < -\tilde  z$ and the set $\{W \leq -z\}$ equals $\{W \leq -\tilde z - \delta\}$. Also recall that we assumed in \eqref{fEC4} that $0 < z < \sqrt{R} -2$ and $\delta < \beta$. Therefore,
\begin{align*}
 \P(W\leq -z) =&\ \P(W=-\tilde z-\delta)+\P(W=-\tilde z-2\delta)+\cdots +\P(W=-\sqrt{R})\\
=&\ \P(W=-\tilde z) \Big( (1-\delta \tilde z) +\big((1-\delta \tilde z)(1-\delta (\tilde z+\delta))\big)\\
&+\big((1-\delta \tilde z)(1-\delta (\tilde z+\delta))(1-\delta (\tilde z+2\delta))\big)+\cdots \Big)\\
\geq&\ \P(W=-\tilde z) \Big( (1-\delta \tilde z) +\big((1-\delta \tilde z)(1-\delta (\tilde z+\delta))\big)+\cdots\\
& +\big((1-\delta \tilde z)\cdots(1-\delta(\tilde z+ \lfloor \frac{1}{\delta}\rfloor \delta))\big) \Big).
\end{align*}
The second equality follows from \eqref{f27} with $\tilde z$ in place of $y$. This requires that  $-\tilde z \leq \beta$, which follows from the fact that $-\tilde z - \delta \leq -z < 0$ and therefore $-\tilde z < \beta$ by our assumption that $\delta < \beta$ in \eqref{fEC4}. In the last inequality, $\lfloor \cdot \rfloor$ denotes the integer part and the inequality itself follows from the fact that $-(\tilde z+ \lfloor \frac{1}{\delta}\rfloor \delta)\geq -(\tilde z+1)>-(z+1)>-\sqrt{R}$. Now for any $0 \leq k \leq   \lfloor \frac{1}{\delta}\rfloor$ we have  $1 - \delta (\tilde z + \delta k) \geq 1 - \delta(\tilde z + 1) > 1 - \delta(z + 1) > 1/2$, where the last inequality follows from our assumption in \eqref{fEC4}. Therefore, the right-hand side above is bounded from below by
\begin{align*}
 & \P(W=-\tilde z) \Big( \big(1-\delta (\tilde z+1)\big) +\big(1-\delta (\tilde z+1)\big)^2+\cdots +\big(1-\delta (\tilde z+1)\big)^{\lfloor \frac{1}{\delta}\rfloor +1}\Big)\\
=&\ \P(W=-\tilde z)\big(1-\delta (\tilde z+1)\big) \frac{1-\big(1-\delta (\tilde z+1)\big)^{\lfloor \frac{1}{\delta}\rfloor+1}}{\delta (\tilde z+1)} \\
\geq & \P(W=-\tilde z)\big(1-\delta (\tilde z+1)\big) \frac{1-(e^{-\delta (\tilde z+1)})^{\lfloor \frac{1}{\delta}\rfloor+1}}{\delta (\tilde z+1)}  \\
\geq& \P(W=-\tilde z)[1-\delta ( z+1)] \frac{1-e^{- 1/2}}{\delta ( z+1)}   \\
\geq&   \P(W=-\tilde z) \frac{1-e^{- 1/2}}{2\delta ( z+1)}. 
\end{align*}
The first inequality is true because $1-x \leq e^{-x}$ and $\delta(\tilde z+1)\geq \delta(z-\delta+1)>0$. The second inequality is true because $0<\delta(\tilde z+1)<\delta (z+1)$ and $\delta (\tilde z+1)(\lfloor 1/\delta\rfloor+1)>(\delta/2) (\lfloor 1/\delta\rfloor+1) > 1/2$   (because $\lfloor 1/\delta\rfloor+1 > 1/\delta$,  $\delta < 1/2$ by \eqref{fEC4}, and $-\tilde z - \delta \leq -z < 0$ from it follows that $\tilde z > -\delta > -1/2$). The last inequality follows from $\delta(z+1)<1/2$ in \eqref{fEC4}. 
This implies that $\P(W=-\tilde z)\leq C\delta(z\vee 1) \P(W\leq \magenta{-z})$.
Using $-\tilde z < \delta < 1/2$ and \eq{f27} again we get
\be{
\P(W=-\tilde z-\delta)=(1-\delta \tilde z) \P(W=-\tilde z)\leq (1+\delta^2) \P(W=-\tilde z)\leq \frac{5}{4} \P(W=-\tilde z),
}
This proves \eq{f23}.  
\finishproof

 \subsubsection{Proof of Lemma~\ref{lem:taylormdhigh}. } \label{sec:prooftaylormdhigh}
\startproof{Proof of Lemma~\ref{lem:taylormdhigh}}
Lemma~\ref{fEC3}, which we  state in Section~\ref{fap4}, implies that $f_z'(x)$ is bounded and absolutely continuous with bounded $f_z''(x)$.   
Lemma~\ref{lem:higherders} implies $\E \abs{W} < \infty$, which when combined with the fact that $f_z'(x)$ is bounded implies $\E \abs{f_z(W)} < \infty$, and in turn $\E G_{\tilde X} f_z(W) = 0$ due to \eqref{eq:gxbarf}.  Letting 
\begin{align*}
G_Y f(x)=\frac{1}{2}a(x)f''(x)+b(x) f'(x),
\end{align*} 
 taking expected values with respect to $W$ in the Poisson equation \eqref{f17} and subtracting $\E G_{\tilde X} f_z(W)=0$ from it, we get 
\ben{ \label{eq:f171}
\E G_Y f_z(W) - \E G_{\tilde X} f_z(W) =\P(Y_1\leq -z) - \P(W\leq -z).
}
To prove the lemma we work on the left-hand side. Similar to \eqref{eq:vareps}, for any function $f: \R \to \R$ with an absolutely continuous derivative,  and any $x,\delta \in \R$, 
\begin{align*}
f(x+\delta) -f(x) =   \delta f'(x)   + \int_x^{x+\delta} (x+\delta -y)f''(y)dy=   \delta f'(x)   + \int_{0}^{\delta} (\delta -y)f''(x+y)dy.  
 \end{align*}
 Applying this expansion to  $G_{\tilde X}$  in  \eqref{eq:gxbarf}, we get 
\begin{align*}
G_{\tilde X} f_z(W) =&\ \delta \big( \lambda -  \mu \big[(W/\delta+R)\wedge n\big]\big) f_z'(W) \\
&+ \lambda \int_{0}^{\delta} (\delta -y) f_z''(W+y) dy  + \mu \big[(W/\delta+R)\wedge n\big] \int_{-\delta}^{0} (y+\delta)f_z''(W+y) dy.
\end{align*}
Recall from \eqref{eq:tabeff} that $b(W) =\delta \big( \lambda -  \mu \big[(W/\delta+R)\wedge n\big]\big) $, and consequently $\mu \big[(W/\delta+R)\wedge n\big] = \lambda - b(W)/\delta$. We recall $K_W(y)$ from \eqref{eq:kdef}, so
\begin{align}
  G_{\tilde X} f_z(W)  =&\ b(W)f_z'(W) +\int_{-\delta}^{\delta} f_z''(W+y) K_{W}(y) dy \notag \\
 =&\ b(W)f_z'(W) + \frac{1}{2} a(W) f_z''(W) +   \int_{-\delta}^{\delta} \big(f_z''(W+y) - f_z''(W)\big) K_{W}(y) dy . \label{eq:steinx}
\end{align}
The last equality follows from $\int_{-\delta}^{\delta} K_W(y) dy = \frac{1}{2}a(W)$ in \eqref{MD:k3}. Thus,
\begin{align}
  \Prob(Y_1 \leq -z) - \Prob(W \leq -z)  =&\ \E G_{Y} f_z(W) - \E G_{\tilde X} f_z(W)  \notag \\ 
  =&\  \E \bigg[ \int_{-\delta}^{\delta} \big(f_z''(W) - f_z''(W+y)\big) K_{W}(y) dy\bigg]. \label{MD:exp1}
\end{align}
In the last equality above we used \eqref{eq:steinx}. The lemma follows from \eqref{f17}; i.e.,   $f_z''(x) = -\frac{2b(x)}{a(x)}f_z'(x) + \frac{2}{a(x)}\big(\Prob(Y_1  \leq -z ) - 1(x \leq -z) \big)$.
\finishproof

\subsubsection{Proof of Lemma~\ref{l4}.} \label{fap4}
 We first present a series of intermediary lemmas that represent the main steps in the proof, and then use them to prove Lemma~\ref{l4}.  We remind the reader that  \eqref{eq:adeff}  implies that
\begin{align}
r(x) = 
\begin{cases}
-2x, \quad x \leq -1/\delta,\\
\frac{-2x}{2 +\delta x}, \quad x \in [-1/\delta, \beta], \\
\frac{-2\beta}{2 + \delta\beta}, \quad x \geq \beta,
\end{cases}
\quad
r'(x) = \begin{cases}
-2, \quad x \leq -1/\delta, \\
\frac{-4}{(2+\delta x)^2}, \quad x \in (-1/\delta, \beta], \\
0, \quad x > \beta,
\end{cases}\label{eq:rform}
\end{align}
where $r'(x)$ is interpreted as the derivative from the left at the points $x = -1/\delta, \beta$. In particular, note that $\abs{r'(x)} \leq 4$. The first lemma decomposes
\begin{align*}
\E  \bigg[ \int_{-\delta}^\delta \big(r(W+y)f_z'(W+y)-r(W)f_z'(W)\big)K_W(y)dy \bigg] 
\end{align*}
into a more convenient form. It is proved in Section~\ref{sec:errorexpproof}.
\begin{lemma} \label{lem:error_expansion}
Let $f_z(x)$ solve \eqref{f17}, then 
\besn{\label{26}
& \int_{-\delta}^\delta \big(r(W+y)f_z'(W+y)-r(W)f_z'(W)\big)K_W(y)dy  \\
=&\int_{-\delta}^\delta  K_W(y)  r(W+y) y f_z''(W)     dy\\
&- \int_{-\delta}^\delta  K_W(y)  r(W+y)\int_0^y \int_0^s r(W+u)f_z''(W+u) du ds     dy \\
&- \int_{-\delta}^\delta  K_W(y)   r(W+y)\int_0^y \int_0^s r'(W+u)f_z'(W+u) du ds     dy \\
&-\int_{-\delta}^\delta  K_W(y) r(W+y)\int_0^y \bigg[\frac{ 1(W+s\leq -z)}{a(W+s)/2}-\frac{ 1(W\leq -z)}{a(W)/2} \bigg]ds      dy \\
&\magenta{+} \P(Y_1\leq -z) \int_{-\delta}^\delta  K_W(y) r(W+y) \int_0^y \bigg[\frac{2}{a(W+s)}-\frac{2}{a(W)}\bigg] ds      dy \\
&+1(W=-1/\delta) f_z'(W) \int_0^\delta K_W(y) \int_0^y r'(W+s) ds dy\\
&+1(W=\beta) f_z'(W) \int_{-\delta}^0 K_W(y) \int_0^y r'(W+s) ds dy\\
&+1(W\in [-1/\delta+\delta, \beta-\delta]) f_z'(W) \int_{-\delta}^\delta K_W(y) \int_0^y \int_0^s r''(W+u) duds dy\\
&+1(W\in [-1/\delta+\delta, \beta-\delta]) f_z'(W) r'(W) \frac{\delta^2 b(W)}{6}.
}
\end{lemma}
Next we  need a bound on $f_z'(x)$ and $f_z''(x)$ as provided by the following lemma.
\begin{lemma}\label{fEC3}
Let $f_z(x)$ solve \eqref{f17}. There exists a constant $C$ that depends only on $\beta$ such that for any $z > 0$, 
\begin{align}
\abs{f_z'(x)} \leq&\ \frac{C}{\mu } \bigg(1(x \leq -z) \Big( \frac{\mu}{\abs{b(x)}} \wedge 1 \Big) \notag \\
& \qquad + \P(Y_1\leq -z) \Big( 1(-z < x < 0)   e^{-\int_{0}^{x} r(u)du} + 1(x \geq 0)\Big( \frac{\mu}{\abs{b(x)}} \wedge 1 \Big) \Big) \bigg), \label{eq:mdleft1} \\
\abs{f_z''(x)} \leq&\ \frac{C}{\mu } \bigg(1(x \leq -z)  + \P(Y_1\leq -z) \Big( 1(-z < x < 0)(1 + \abs{x})   e^{-\int_{0}^{x} r(u)du} + 1(x \geq 0) \Big) \bigg)\magenta{.} \label{eq:mdleft2}
\end{align}
\end{lemma}
\startproof{ Proof of Lemma \ref{fEC3} }
We recall from \eqref{eq:stddenf} that the density of $Y_1$ is 
\begin{align*}
 \kappa\frac{2}{a(x)} \exp\Big({\int_0^x \frac{2b(y)}{a(y)}dy}\Big), \quad x \in \R,
\end{align*}
where $\kappa$ is the normalizing constant. 
 One may verify that the solution to the Poisson equation $\frac{1}{2}a(x)f_z''(x)+b(x)f_z'(x)=\P(Y_1 \leq -z) - 1(x\leq -z)$   satisfies
\ben{\label{f18}
f_z'(x)=
\begin{cases}
-\P(Y_1\geq -z) e^{-\int_0^x r(u)du} \frac{1}{\kappa} \P(Y_1\leq x), & x\leq -z, \\
-\P(Y_1\leq -z) e^{-\int_0^x r(u)du} \frac{1}{\kappa} \P(Y_1\geq x), & x\geq -z\magenta{.}
\end{cases}
}
Rearranging the Poisson equation and using \eq{f18} yields
\besn{\label{f19}
f_z''(x)=&-r(x) f_z'(x) +\frac{2(\P(Y_1\leq -z) - 1(x\leq -z))}{a(x)}\\
=&
\begin{cases}
-\P(Y_1\geq -z) \left[ \frac{2}{a(x)} -\frac{r(x)}{\kappa} e^{-\int_0^x r(u) du} \P(Y_1\leq x) \right], & x\leq -z, \\
-\P(Y_1\leq -z) \left[ -\frac{2}{a(x)} -\frac{r(x)}{\kappa} e^{-\int_0^x r(u) du} \P(Y_1\geq x)\right], & x>-z.
\end{cases}
}
We will start by proving the bound on $f_z''(x)$.  The form of the density of $Y_1$ implies
\begin{align*}
\frac{1}{\kappa} \P(Y_1\leq x) = \int_{-\infty}^{x} \frac{2}{ a(y)} e^{\int_{0}^{y}r(u)du} dy.
\end{align*}
Since $b(x)$ is nonincreasing and $b(x) > 0$ when $x < 0$, then for $x < 0$, 
\begin{align}
\frac{1}{\kappa} \P(Y_1\leq x) = \int_{-\infty}^{x} \frac{2}{ a(y)} e^{\int_{0}^{y}r(u)du} dy \leq \frac{1}{b(x)} \int_{-\infty}^{x} \frac{2 b(y)}{ a(y)} e^{\int_{0}^{y}r(u)du} dy = \frac{1}{b(x)} e^{\int_{0}^{x}r(u)du}, \label{eq:main1}
\end{align}
and  since $b(x) < 0$ for $x > 0$, then for $x \magenta{>} 0$, 
\begin{align}
\frac{1}{\kappa} \P(Y_1\geq x) \leq \frac{1}{b(x)} \int_{x}^{\infty} \frac{2 b(y)}{ a(y)} e^{\int_{0}^{y}r(u)du} dy = \frac{-1}{b(x)} e^{\int_{0}^{x}r(u)du}= \frac{ 1}{\abs{b(x)}} e^{\int_{0}^{x}r(u)du}. \label{eq:main2}
\end{align}
Applying \eqref{eq:main1} and \eqref{eq:main2} to the form of $f_z''(x)$ above, we get the desired upper bound  when $x\leq -z$ or $x \geq 0$. When $-z < x < 0$, we use the fact that $\abs{r(x)} \leq C \abs{x}$ to get
\begin{align*}
\abs{f_z''(x)} \leq&\  \magenta{\frac{2}{a(x)}\P(Y_1\leq -z)+C} \P(Y_1\leq -z)\magenta{\abs{x}}   e^{-\int_0^x r(u) du} \int_{x}^{\infty} \frac{2}{ a(y)} e^{\int_{0}^{y} r(u)du} dy \\
\leq&\  \frac{C}{\mu } \P(Y_1\leq -z)(1 + \abs{x})   e^{-\int_0^x r(u) du}.
\end{align*}
 The last inequality follows from the facts that $a(x) \geq \mu$, and that $\int_{x}^{\infty}  e^{\int_{0}^{y} r(u)du} dy$ can be bounded by a constant that depends only on $\beta$, which is evident from the form of $r(x)$. This establishes the bound on $f_z''(x)$. Repeating the same procedure with $f_z'(x)$ gives us the bound 
\begin{align*}
\abs{f_z'(x)} \leq&\ C \bigg(1(x \leq -z) \frac{1}{\abs{b(x)}} + \P(Y_1\leq -z) \Big( 1(-z < x < 0) \magenta{\frac{1}{\mu}}  e^{-\int_{0}^{x} r(u)du} + 1(x \geq 0)\frac{1}{\abs{b(x)}} \Big) \bigg).
\end{align*}
To conclude the proof, we require the following two inequalities from Lemma B.8 of \cite{Brav2017}:
\begin{align}
&  e^{-\int_{0}^{x} r(u)du}\int_{-\infty}^{x} \frac{2}{ a(y)} e^{\int_{0}^{y}r(u)du} dy\leq 
\begin{cases}
\frac{3}{ \mu }, \quad x \leq 0,\\
\frac{1}{\mu}e^{\beta^2} (3 + \beta), \quad x \in [0,\beta],
\end{cases} \label{DS:fbound1}\\
& e^{-\int_{0}^{x}r(u)du}\int_{x}^{\infty} \frac{2}{ a(y)} e^{\int_{0}^{y} r(u)du} dy \leq 
\begin{cases}
\frac{1}{\mu } \Big( 2 + \frac{1}{\beta} \Big), \quad x \in [0,\beta], \\
\frac{1}{\mu \beta}, \quad x \geq \beta.
\end{cases} \label{DS:fbound2}
\end{align}
By repeating the bounding procedure discussed above, but using \eqref{DS:fbound1} and \eqref{DS:fbound2} in place of \eqref{eq:main1} and \eqref{eq:main2}, we arrive at 
\begin{align*}
\abs{f_z'(x)} \leq&\ \frac{C}{\mu } \bigg(1(x \leq -z)   + \P(Y_1\leq -z) \Big( 1(-z < x < 0)   e^{-\int_{0}^{x} r(u)du} + 1(x \geq 0) \Big) \bigg),
\end{align*}
which establishes the desired bound on $f_z'(x)$. 
\finishproof
Lastly, we will need the following lemmas, which are proved in Section~\ref{fap5}.
\begin{lemma}\label{fl2}
There exist constant $c_1,C_1 > 0$ such that for any integer  $k \geq 0$, some $C(k) > 0$, and any $R\geq C_1$ and $0<z\leq c_1 R^{1/4}$,
\ben{\label{11}
\E|1(W\leq  -z) W^k|\leq C(k) (z\vee 1)^{k+1} \P(Y_1\leq -z),
}
The constants $c_1,C_1,C(k)$  depend on  $\beta$. 
\end{lemma}

\begin{lemma}\label{fl3}
There exist constant $c_1,C_1 > 0$ such that for any integer  $k \geq 0$, some $C(k) > 0$, and any $R\geq C_1$ and $0<z\leq c_1 R^{1/4}$,
\ben{\label{14}
\E|1(-z\leq W\leq 0)W^k e^{-\int_0^W r(u)du}|\leq C(k) (z\vee 1)^{k+1},
}
The constants $c_1,C_1,C(k)$  depend on  $\beta$. 
\end{lemma}
We are now ready to prove Lemma~\ref{l4}. 
\startproof{Proof of Lemma~\ref{l4}}
Again, in the following, we C to denote a constant  whose value may change from line to line but only depends on $\beta$. We take expected values on both sides of \eqref{26} and bound the terms on the right-hand side one line at a time. For the first line, we need to bound
\bes{
\left| \E \bigg[ \int_{-\delta}^\delta  K_W(y)  r(W+y) y f_z''(W)     dy \bigg] \right|
\leq &\left| \E \bigg[  r(W) f_z''(W) \int_{-\delta}^\delta  y K_W(y)      dy \bigg]\right|\\
&+\left| \E \bigg[ f_z''(W) \int_{-\delta}^\delta y K_W(y) (r(W+y)-r(W))       dy \bigg] \right|.
}
Using \eqref{MD:k3}, it follows that 
\begin{align*}
\left| \E \bigg[  r(W) f_z''(W) \int_{-\delta}^\delta  y K_W(y)      dy \bigg]\right| =&\ \left| \E \bigg[  r(W) f_z''(W) \frac{ \delta^2 b(W)}{6}\bigg]\right| \leq C \mu \delta^2 \E \abs{W^2 f_{z}''(W)},
\end{align*}
where we used $\abs{r(x)} \leq C \abs{x}$ and $\abs{b(x)}\leq \mu \abs{x}$ in the inequality. Furthermore, since $\abs{a(x)} \leq C \mu$ and recalling from \eqref{eq:rform} that $\abs{r'(x)} \leq 4$, we have  
\begin{align*}
 \left| \E \bigg[ f_z''(W) \int_{-\delta}^\delta y K_W(y) (r(W+y)-r(W))       dy \bigg] \right| \leq&\ \magenta{ \E \bigg[ C |f_z''(W)|\delta^2 \int_{-\delta}^\delta  K_W(y)         dy \bigg] }   \\
=&\  \magenta{ \E \bigg[ C |f_z''(W)|\delta^2 \frac{a(W)}{2} \bigg] }  \\
\leq &\ C \mu \delta^2 \E \abs{f_z''(W)}.
\end{align*}
Therefore, 
\begin{align*}
\left| \E \bigg[ \int_{-\delta}^\delta  K_W(y)  r(W+y) y f_z''(W)     dy \bigg] \right| \leq C \mu \delta^2 \E \abs{(1 + W^2) f_z''(W)}. 
\end{align*}
Applying the bound on $\abs{f_z''(x)}$ from Lemma~\ref{fEC3},  we see that the right-hand side above is bounded by
\begin{align*}
  & \magenta{C\delta^2} \E \bigg((1 + W^2)1(W \leq -z)  \\
  &\qquad + \P(Y_1\leq -z)(1 + W^2) \Big( 1(-z < W < 0)(1 + \abs{W})   e^{-\int_{0}^{W} r(u)du} + 1(W \geq 0) \Big) \bigg)\\
  \leq&\ C \P(Y_1\leq -z) \magenta{\delta^2} (z\vee 1)^4,
\end{align*} 
where the inequality is due to  Lemmas~\ref{fl2} and \ref{fl3} and the fact that  $\E W^2 \leq C$, which was proved in Lemma A.1 of \cite{Brav2017}.
Following the same argument, the second line 
\begin{align*}
&\abs{\int_{-\delta}^\delta  K_W(y)  r(W+y)\int_0^y \int_0^s r(W+u)f_z''(W+u) du ds     dy}\\
&\magenta{\leq C\mu\delta^2\E (1+W^2) \sup_{-\delta\leq u\leq \delta}|f_z''(W+u)|} \\
&\magenta{\leq C\delta^2 \E \bigg( (1+W^2)1(W\leq -z+\delta)} \\
&\magenta{\quad + \P(Y_1\leq -z) (1+W^2) \Big( 1(-z-\delta<W<\delta)(1+|W|)e^{\sup_{-\delta\leq u\leq \delta}(-\int_0^{W+u} r(v)dv)}+1(W\geq -\delta)       \Big) \bigg),}
\end{align*} 
\magenta{which} is also bounded by $C \P(Y_1\leq -z) \delta^2 (z\vee 1)^4$  from simple modifications of Lemmas~\ref{fl2} and \ref{fl3}. For the third line, we use the representation $f_z'(W+u) - f_z'(W) = \int_{0}^{u} f_z''(W+v) dv$ and the triangle inequality to get
\begin{align*}
&\left| \E \bigg[\int_{-\delta}^\delta  K_W(y)    r(W+y)\int_0^y \int_0^s r'(W+u)f_z'(W+u) du ds     dy \bigg]\right| \\
\leq&\ \left| \E \bigg[\int_{-\delta}^\delta  K_W(y)    r(W+y)\int_0^y \int_0^s r'(W+u)\int_{0}^{u} f_z''(W+v) dvdu ds     dy \bigg]\right| \\
&+   \left| \E \bigg[\int_{-\delta}^\delta  K_W(y) f_z'(W) (r(W+y)-r(W))\int_0^y \int_0^s r'(W+u) du ds     dy \bigg]\right| \\
&+  \left| \E \bigg[\int_{-\delta}^\delta  K_W(y) f_z'(W) r(W)  \int_0^y \int_0^s r'(W+u) du ds     dy \bigg]\right| \\
\leq&\magenta{\ C \E \bigg[ (\abs{W} + \delta)\int_{-\delta}^\delta  K_W(y)    \int_0^y \int_0^s \int_{0}^{u} |f_z''(W+v)| dvdu ds     dy \bigg]} \\
&\magenta{+    C \delta^3 \E \bigg[\abs{f_z'(W)} \int_{-\delta}^\delta  K_W(y)dy  \bigg] +   C\delta^2 \E \bigg[\abs{r(W) f_z'(W)}\int_{-\delta}^\delta  K_W(y)  dy  \bigg]}\\
\leq&\magenta{\ C \E \bigg[ (\abs{W} + \delta)\int_{-\delta}^\delta  K_W(y)    \int_0^y \int_0^s \int_{0}^{u} |f_z''(W+v)| dvdu ds     dy \bigg]} \\
&\magenta{+    C \mu \delta^3 \E  \abs{f_z'(W)}   +   C\delta^2 \E  \abs{b(W) f_z'(W)}  .}
\end{align*}
In the second inequality we used $\abs{r'(x)}\leq 4$  and the last inequality is due to \eqref{MD:k3}, that $r(x) = 2b(x)/a(x)$ and the fact that $\abs{a(x)} \leq C \mu$. The right-hand side can be bounded by $C\P(Y_1\leq -z) \delta^2 (z\vee 1)^4$ as follows.
The first term on the right-hand side above can be bounded by repeating the procedure used to bound lines one and two of \eqref{26}. The second and third terms can be bounded by combining the bound on $f_z'(x)$ from Lemma~\ref{fEC3} with Lemmas~\ref{fl2} and \ref{fl3}.  
For the fourth line, repeating the arguments from \eq{f22}, \eqref{eq:interm30} and \eq{f23}, we have
\bes{
&\left| \E \bigg[\int_{-\delta}^\delta  K_W(y) r(W+y)\int_0^y \bigg[\frac{ 1(W+s\leq -z)}{a(W+s)/2}-\frac{ 1(W\leq -z)}{a(W)/2} \bigg]ds      dy \bigg] \right|\\
\leq & \left| \E \bigg[ \int_{-\delta}^\delta  K_W(y) r(W+y)\int_0^y \bigg[\frac{2}{a(W+s)}-\frac{2}{a(W)} \bigg]1(W\leq -z)ds      dy \bigg] \right|\\
&+\left| \E \bigg[ \int_{-\delta}^\delta  K_W(y) r(W+y)\int_0^y \bigg[\frac{2(1(W+s\leq -z)-1(W\leq -z))}{a(W+s)} \bigg]ds      dy \bigg] \right|\\
\leq & C\delta^3 \magenta{\E(1+|W|) 1(W\leq -z)} + C\delta \magenta{(z\vee 1)} [\P(W=-\tilde{z})+\P(W=-\tilde{z}-\delta)]\\
\leq & \magenta{C\P(Y_1\leq -z) \delta^3 (z\vee 1)^2+ C\delta^2 (z\vee 1)^2 \P(W\leq -z)}\\
\leq & C\P(Y_1\leq -z) \delta^2 (z\vee 1)^3,
}
where in the last \magenta{two inequalities} we used Lemma~\ref{fl2}. The fifth line is bounded as follows. Using \eqref{MD:k3}, \eqref{eq:alipsch} and the fact that $\abs{r(x)} \leq C\abs{x}$, we get
\bes{
&\left| \E \bigg[ \int_{-\delta}^\delta  K_W(y) r(W+y) \P(Y_1\leq -z) \int_0^y \bigg[\frac{2}{a(W+s)}-\frac{2}{a(W)}\bigg] ds      dy \bigg] \right|\\
\leq & C\P(Y_1\leq -z) \delta^3 \magenta{\E(1+|W|)}\leq C\P(Y_1\leq -z) \delta^3.
}
The sixth line is bounded as follows. Using $\abs{r'(x)} \leq 4$, \eqref{MD:k3}, and $\abs{1(W=-1/\delta) f_z'(W)} \leq 1(W=-1/\delta) C/\mu$, which follows from Lemma~\ref{fEC3} and the fact that $z < \sqrt{R}-2 < \sqrt{R} = 1/\delta$ from \eqref{fEC4}, we have 
\bes{
  \left|\E \bigg[ 1(W=-1/\delta) f_z'(W) \int_0^\delta K_W(y) \int_0^y r'(W+s) ds dy\bigg]\right|  \leq & C \delta \P(W=-1/\delta),
}
and now we argue that 
\begin{align}
C \delta \P(W=-1/\delta)  \leq C \delta^2  (z \vee 1)\P(Y_1 \leq -z). \label{eq:toargue}
\end{align}
For any integer $0 < k < R$, \magenta{from \eq{f26}}
\begin{align*}
\P(W=-1/\delta) = \P(W = -\sqrt{R}) \leq  \frac{1}{k} \sum_{i=0}^{k} \P(W = \delta(i-R)) = \frac{1}{k} \P(W \leq \delta (k-R)).
\end{align*}
Choose $k$ such that $\delta (k-R)$ is the largest element in the support of $W$ that is less than or equal to $-z$, and use the fact that $z \leq c_1 R^{1/4}$ to conclude that 
\begin{align*}
\P(W=-1/\delta) \leq  C \delta \P(W \leq -z).
\end{align*}
Using Lemma~\ref{fl2} implies \eqref{eq:toargue}.  The seventh line is bounded as follows. Using $\abs{r'(x)} \leq 4$, \eqref{MD:k3}, and $\abs{1(W=\beta) f_z'(W)} \leq 1(W=\beta) C/\mu$, which follows from Lemma~\ref{fEC3}, we have
\bes{
& \left| \E \bigg[1(W=\beta) f_z'(W) \int_{-\delta}^0 K_W(y) \int_0^y r'(W+s) ds dy\bigg]\right| \\
\leq & C \delta \P(W=\beta)\P(Y_1\leq -z)\leq C\delta^2 \P(Y_1\leq -z),
}
where we used $\P(W=\beta)\leq C\delta$, which was argued at the end of the proof of Proposition~\ref{prop:lowerboundec}. The eighth line is bounded as follows. Using $\abs{r''(x)} \leq C \delta$ and  Lemma~\ref{fEC3}, 
\bes{
& \left| \E \bigg[1(W\in [-1/\delta+\delta, \beta-\delta]) f_z'(W) \int_{-\delta}^\delta K_W(y) \int_0^y \int_0^s r''(W+u) duds dy \bigg] \right| \\
\leq & C \delta^3 \P(W\leq -z)  +C \delta^3 \P(Y_1\leq -z) \E|1(-z\leq W\leq 0) e^{-\int_0^W r(u)du}| +C \delta^3 \P(Y_1\leq -z)\\
\leq & C \delta^3 (z\vee 1) \P(Y_1\leq -z),
}
where in the last inequality we used \magenta{Lemmas~\ref{fl2} and \ref{fl3}}. The ninth line is bounded similarly. Using $\abs{r'(x)} \leq 4$ and the bound on $\abs{f_z'(x) b(x)} $ from  Lemma~\ref{fEC3}, we have
\bes{
& \left| E\{1(W\in [-1/\delta+\delta, \beta-\delta]) f_z'(W) r'(W) \frac{\delta^2 b(W)}{6}\right| \\
 \leq&\ C \delta^2 \E \abs{W 1(W\leq -z)}  +C \delta^2 \P(Y_1\leq -z) \E|1(-z\leq W\leq 0) W e^{-\int_0^W r(u)du}| +C \delta^2 \P(Y_1\leq -z) \\
    \leq&\  C \delta^2 (z\vee 1)^2 \P(W\leq -z) \magenta{+C\delta^2 \P(Y_1\leq -z)} \leq   C \delta^2 (z\vee 1)^3 \P(Y_1\leq -z).
}
The second inequality is due to \magenta{Lemmas~\ref{fl2} and \ref{fl3}. }
 The last inequality is due to Lemma~\ref{fl2} with $k = 0$ there. Combining the bounds proves Lemma~\ref{l4}.
\finishproof

\subsubsection{Proof of Lemma~\ref{lem:error_expansion}.} \label{sec:errorexpproof}
\startproof{Proof of Lemma~\ref{lem:error_expansion}}
 Assume for now that for all $y \in (-\delta, \delta)$, 
\begin{align}
&r(W+y)f_z'(W+y) - r(W)f_z'(W) \notag \\
=&\ y r(W+y)f_z''(W) -  r(W+y)\int_{0}^{y}\int_{0}^{s} \Big(r(W+u) f_z''(W+u) + r'(W+u)f_z'(W+u) \Big)duds  \notag  \\
&- r(W+y)\int_{0}^{y}\Big(\frac{2}{a(W+s)} 1(W+s \leq -z) - \frac{2}{a(W)} 1(W \leq -z) \Big)ds \notag \\
&+ \Prob(Y_1 \leq -z)r(W+y)\int_{0}^{y}\Big(\frac{2}{a(W+s)} - \frac{2}{a(W)} \Big)ds + f_z'(W) \int_{0}^{y}r'(W+s)ds. \label{eq:interm}
\end{align}
We postpone verifying \eqref{eq:interm} to the end of this proof, but \eqref{eq:interm} implies
\begin{align*}
&\int_{-\delta}^{\delta} \Big(r(W+y)f_z'(W+y) - r(W)f_z'(W) \Big) K_W(y) dy  \notag \\ 
=&\ \int_{-\delta}^{\delta} K_W(y) y r(W+y)f_z''(W) dy  \\
&- \int_{-\delta}^{\delta}  K_W(y)  r(W+y)\int_{0}^{y}\int_{0}^{s} r(W+u) f_z''(W+u)  duds  dy \\
&-   \int_{-\delta}^{\delta}  K_W(y)  r(W+y)\int_{0}^{y}\int_{0}^{s} r'(W+u)f_z'(W+u)  duds  dy\\  
&-   \int_{-\delta}^{\delta}  K_W(y)  r(W+y) \int_{0}^{y}\Big(\frac{2}{a(W+s)} 1(W+s \leq -z) - \frac{2}{a(W)} 1(W \leq -z) \Big)ds  dy\\  
&+ \Prob(Y_1\magenta{\leq -z}) \int_{-\delta}^{\delta}K_W(y)r(W+y)\int_{0}^{y}\Big(\frac{2}{a(W+s)} - \frac{2}{a(W)} \Big)ds  dy \\
&+ f_z'(W)\int_{-\delta}^{\delta}K_W(y) \int_{0}^{y}r'(W+s)dsdy.
\end{align*}
We are almost done, but the last term on the right-hand side above requires some additional manipulations. Since $r'(x) = 0$ for $x > \beta$ and $K_W(y)=0$ for $W = -1/\delta$ and $y \in [-\delta, 0]$,
\begin{align*}
&\int_{-\delta}^{\delta} K_W(y)\int_{0}^{y}r'(W+s)ds  dy \\
 =&\ 1(W = -1/\delta)\int_{0}^{\delta} K_W(y)\int_{0}^{y}r'(W+s)ds  dy\\
&+ 1(W = \beta) \int_{-\delta}^{0} K_W(y)\int_{0}^{y}r'(W+s)ds  dy \\
&+ 1(W \in [-1/\delta + \delta, \beta - \delta])  \int_{-\delta}^{\delta} K_W(y)\int_{0}^{y}r'(W+s)ds  dy,
\end{align*}
and for $W \in [-1/\delta + \delta, \beta - \delta]$, 
\begin{align*}
&\int_{-\delta}^{\delta} K_W(y)\int_{0}^{y}r'(W+s)ds  dy \\
=&\  \int_{-\delta}^{\delta} K_W(y)\int_{0}^{y} \big(r'(W+s) - r'(W)\big)ds  dy +r'(W) \int_{-\delta}^{\delta} yK_W(y)dy \\
=&\  \int_{-\delta}^{\delta} K_W(y)\int_{0}^{y} \int_{0}^{s}r''(W+u)duds  dy + r'(W) \frac{1}{6}\delta^2b(W).
\end{align*}
To conclude the proof, we verify \eqref{eq:interm}:
\begin{align*}
&r(W+y)f_z'(W+y) - r(W)f_z'(W)  \notag \\
=&\ r(W+y)f_z'(W) + r(W+y) \int_{0}^{y} f_z''(W+s) ds  -  r(W)f_z'(W) \notag \\
=&\ r(W+y) \int_{0}^{y} f_z''(W+s) ds \magenta{+} f_z'(W) \int_{0}^{y} r'(W+s) ds.
\end{align*}
Now
\begin{align*}
\int_{0}^{y}  f_z''(W+s)ds =&\  yf_z''(W) + \int_{0}^{y}  \big(f_z''(W+s) - f_z''(W) \big)ds \\
=&\ y f_z''(W) -  \int_{0}^{y}\Big(r(W+s)f_z'(W+s) - r(W)f_z'(W) \Big)ds   \\
&- \int_{0}^{y}\Big(\frac{2}{a(W+s)} 1(W+s \leq -z) - \frac{2}{a(W)} 1(W \leq -z) \Big)ds\\
&+ \Prob(Y_1 \leq -z)\int_{0}^{y}\Big(\frac{2}{a(W+s)} - \frac{2}{a(W)} \Big)ds,
\end{align*}
and the fundamental theorem of calculus tells us that 
\begin{align*}
&r(W+s)f_z'(W+s) - r(W)f_z'(W)\\
 =&\ \int_{0}^{s} \Big(r(W+u) f_z''(W+u) + r'(W+u)f_z'(W+u) \Big)du,
\end{align*}
which proves \eqref{eq:interm}.
\finishproof

\subsubsection{Proving Lemmas~\ref{fl2} and \ref{fl3}.}\label{fap5}
We first prove an upper bound on $E e^{-tW}$ for $0\leq t\leq c_1 R^{1/4}$. Note that $\E e^{-t W} < \infty$ for $ t \geq 0$ because  $W\geq -\sqrt{R}$ (cf. \eq{f29}).
\begin{lemma}\label{l1}
There exist $c_1,C_1>0$ such that if  $R\geq C_1$ and $0\leq t\leq c_1 R^{1/4}$, then
\ben{\label{04}
\E e^{-tW}\leq C e^{\frac{t^2}{2}-\frac{\delta t^3}{6}}.
}
\end{lemma}
\startproof{Proof of Lemma \ref{l1}}
For $0\leq s\leq t$, define $h(s)=\E e^{-sW}$, 
so $h'(s)= -\E \big(W e^{-sW}\big)$. We will shortly prove that 
\ben{\label{05}
\Big(1+\frac{\delta s}{2}+ \frac{\delta^2 s^2}{6}\Big) h'(s)\leq (s+C\delta^2 s^3) h(s).
}
from which we have 
\begin{align*}
h'(s)\leq  \Big(\frac{s}{1+\frac{\delta s}{2}+ \frac{\delta^2 s^2}{6}}+\frac{C \delta^2 s^3}{1+\frac{\delta s}{2}+ \frac{\delta^2 s^2}{6}}\Big) h(s) \leq&\  \Big(\frac{s}{1+\frac{\delta s}{2}}+C \delta^2 t^3\Big) h(s) \\
\leq&\ \Big(s\big(1-\delta s/2 +  (\delta s/2)^2 \big) +C\delta^2 t^3\Big)h(s) \\
\leq&\  \Big(s-\delta s^2/2   +C\delta^2 t^3\Big)h(s).
\end{align*}
Above we have used $0 \leq s \leq t$. The \magenta{third} inequality follows from the inequality $1/(1+x) \leq (1-x+x^2)$ for $x \geq 0$. 
Since $h(0) = 1$, 
\be{
\log (h(t)) = \int_{0}^{t} \frac{h'(s)}{h(s)} ds  \leq \int_0^t (s-\frac{\delta s^2}{2}+C\delta^2 t^3) ds\leq \frac{t^2}{2}-\frac{\delta t^3}{6}+C\delta^2 t^4\leq \frac{t^2}{2}-\frac{\delta t^3}{6}+\magenta{C}c_1^4,
}
which implies \eq{04}. We are left to prove \eq{05}. Recall from \eqref{eq:gxbarf} that $\E G_{\magenta{\tilde X}} f(W) = 0$ provided  $\E \abs{f(W)} < \infty$. Choose $f(x) = -e^{-sx}/s$, so that $f'(x) = e^{-sx}$ and $f''(x) = -s e^{-sx}$. The form of $G_{\magenta{\tilde X}} f(x)$ in \eqref{eq:steinx} implies  
\begin{align}
\E \bigg[ \frac{a(W)}{2}f''(W)+b(W)f'(W)+\int_{-\delta}^\delta (f''(W+y)-f''(W)) K_W(y) dy \bigg] =0 \label{eq:identstein}
\end{align}
Using the form of $b(x)$ from \eqref{eq:adeff}, we have
\besn{\label{06}
\E[b(W)f'(W)] =&\ \E[1(W\leq \beta) (-\mu W) e^{-sW}+1(W>\beta)(-\mu \beta) e^{-sW}]\\
=&\mu h'(s)+\E[1(W>\beta)\mu(W-\beta)e^{-sW}]\\
\geq & \mu h'(s).
}
Similarly, using the form of $a(x)$ from \eqref{eq:adeff} and the fact that $W \geq -1/\delta$, we have 
\besn{\label{07}
 -\E\Big[\frac{a(W)}{2}f''(W)\Big] =&\ \mu \E\Big[1(-\frac{1}{\delta}\leq W\leq \beta) (1+\frac{\delta W}{2})s e^{-sW}+1(W>\beta)(1+\frac{\delta \beta}{2})s e^{-sW}\Big]\\
=&\ \mu \E\Big[(1+\frac{\delta W}{2})s e^{-sW}\Big]\magenta{-}\mu \E\Big[1(W>\beta) \frac{\delta (W-\beta)}{2}s e^{-sW}\Big]\\
\leq&\ \mu \E\Big[(1+\frac{\delta W}{2})s e^{-sW}\Big]=\mu sh(s)-\mu \frac{\delta s}{2}h'(s).
}
Lastly,
\begin{align*}
-\E \bigg[ \int_{-\delta}^\delta (f''(W+y)-f''(W)) K_W(y) dy\bigg] =&\ \E \bigg[\int_{-\delta}^\delta s (e^{-s(W+y)}-e^{-sW}) K_W(y)dy\bigg] \nonumber\\
\leq&\ \E \bigg[se^{-sW} \int_{-\delta}^\delta (-sy+s^2 y^2e^{\abs{sy}} )K_W(y)dy\bigg] \\
\leq&\ \E \bigg[se^{-sW} \int_{-\delta}^\delta  (-sy+ C s^2 y^2   )K_W(y)dy\bigg].
\end{align*}
The first inequality is due to $e^{-x}-1 \leq -x + x^2 e^{\abs{x}}$, and the second inequality is due to $|sy|\leq c_1 R^{1/4}\frac{1}{\sqrt{R}}\leq C$. Recall from \eqref{MD:k3} that $\int_{-\delta}^\delta K_W(y)dy = \frac{1}{2} a(W) $ and that $\int_{-\delta}^\delta y K_W(y)dy = \frac{\delta^2 b(W)}{6}$, and from \eqref{fEC6} that $\abs{a(x)} \leq C \mu$. Also note from \eqref{eq:adeff} that $b(x) = -\mu x - 1(x > \beta)(\mu \beta - \mu x)$. Therefore, the right-hand side above is bounded by
\begin{align}
& \E \bigg[s^2e^{-sW} \frac{-\delta^2 b(W)}{6}\bigg]+ C\mu \delta^2 s^3 Ee^{-sW} \nonumber\\
=&\ \frac{\delta^2 s^2}{6} \E \big[\mu W e^{-sW}\big] +\frac{\delta^2 s^2}{6} \E \Big[1(W>\beta) (\mu \beta - \mu W) e^{-sW}\Big]+ C\mu \delta^2 s^3 \E e^{-sW}\nonumber\\
=&\ -\mu \frac{\delta^2 s^2}{6} h'(s) +\frac{\delta^2 s^2}{6} \E \Big[1(W>\beta) (\mu \beta - \mu W) e^{-sW}\Big]+C\mu \delta^2 s^3 h(s)\nonumber\\
\leq&\ -\mu \frac{\delta^2 s^2}{6} h'(s)+C\mu \delta^2 s^3 h(s).\label{08}
\end{align}
Combining \eqref{06}--\eqref{08} with \eqref{eq:identstein} concludes the proof. 
\finishproof
Now we are ready to prove Lemmas~\ref{fl2} and \ref{fl3}. 
\proof{Proof of Lemma~\ref{fl2}}
Suppose $k \geq 1$. We prove the lemma by showing that 
\begin{align}
\E|1(W\leq  -z) W^k| \leq&\ C(k) (z \vee 1)^{k} e^{-\frac{z^2}{2}-\frac{\delta z^3}{6}},  \notag \\
\text{ and } \quad C \frac{1}{z} e^{-\frac{z^2}{2}-\frac{\delta z^3}{6}} \leq&\ \Prob(Y_1 \leq -z). \label{eq:l1part1}
\end{align}
Let us start with the first inequality. 
\magenta{For $0<z<1$, the first inequality holds because of Lemma \ref{l1}. Therefore, we only need to consider the case $z\geq 1$.}
Integration by parts yields
\be{
\E|1(W\leq  -z) W^k|=z^k \P(W\leq -z)+\int_{-\infty}^{-z} k (-y)^{k-1} \P(W\leq y) dy.
}
For the first term, note that
\ben{\label{12}
z^k\P(W\leq -z) \leq  z^{k} \E\Big( e^{z(-z -W)}   1(W \leq -z) \Big)  \leq z^k e^{-z^2} \E e^{-zW}\leq Cz^k e^{-\frac{z^2}{2}-\frac{\delta z^3}{6}},
}
where the last inequality is due to Lemma~\ref{l1}. For the second term, we have
\bes{
 \int_{-\infty}^{-z} k (-y)^{k-1} \P(W\leq y) dy \magenta{\leq}&\ \int_{-\infty}^{-z} k (-y)^{k-1} \E\Big( e^{z(y-W)}   1(W \leq y) \Big) dy\\
 \leq&\ \int_{-\infty}^{-z} k (-y)^{k-1} e^{zy} \E e^{-zW} dy\\
\leq&\ C e^{\frac{z^2}{2}-\frac{\delta z^3}{6}} \int_{-\infty}^{-z} k (-y)^{k-1}e^{zy} dy\\
\leq&\ C(k) z^{k-2} e^{-\frac{z^2}{2}-\frac{\delta z^3}{6}}.
}
The second-last inequality is due to Lemma~\ref{l1} and in the last inequality we used  $\int_{-\infty}^{-z} (-y)^{k-1}e^{zy}dy\leq C(k)z^{k-2}e^{-z^2}$ for $k\geq 1$ \magenta{and $z\geq 1$}. This proves the first inequality in \eqref{eq:l1part1}, and we now argue the second one. Recall that  $r(x) = 2b(x)/a(x)$ and that the density of $Y_1$ in \eqref{eq:stddenf} implies 
\bes{
 \P(Y_1\leq -z)=  \int_{-\infty}^{-z} \frac{2\kappa}{a(y)} e^{\int_{0}^{y} r(u) du} dy \geq  \int_{-z-1}^{-z} \frac{2\kappa}{a(y)} e^{\int_{0}^{y} r(u) du} dy =  &\  \int_{-z-1}^{-z} \frac{2\kappa}{a(y)}  e^{\int_{0}^{y} \frac{-2u}{2 + \delta u} du} dy.
} 
The last equality follows from the assumption  $z \leq \sqrt{R} - 2 $ in \eqref{fEC4} so that $-z-1 \geq -\sqrt{R}  $, meaning $a(x) = 2\mu - \delta b(x)$ for $x \in [-z-1, \magenta{0}]$. We now argue that $2\kappa/a(x) \geq C$ for some constant $C> 0 $ that depends only on $\beta$. The density of $Y_1$ is given by \eqref{eq:stddenf}, so we know that the normalizing constant $\kappa$ satisfies  
\begin{align*}
\frac{1}{\kappa} =&\  \int_{-\infty}^{\infty} \frac{2}{a(x)} e^{\int_{0}^{x} 2b(y)/a(y) dy} dx \\
 =&\ \int_{-\infty}^{0} \frac{2}{a(x)} e^{\int_{x}^{0} -2\abs{b(y)}/a(y) dy} dx  +  \int_{0}^{\infty} \frac{2}{a(x)} e^{\int_{0}^{x} -2\abs{b(y)}/a(y) dy} dx \\
\leq &\ \frac{C}{\mu} \int_{-\infty}^{0}  e^{\int_{x}^{0} -\frac{2\abs{b(y)}} {\mu (2+\delta \beta)} dy} dx  + \frac{C}{\mu} \int_{0}^{\infty}   e^{\int_{0}^{x} -\frac{2\abs{b(y)}} {\mu (2+\delta \beta)} dy} dx.
\end{align*}
The second equality is true because $b(x)$ in \eqref{eq:adeff} is nonincreasing and $b(0) = 0$, meaning $b(x) = -\abs{b(x)} 1(x \geq 0) + \abs{b(x)} 1(x < 0)$. The inequality is true because $a(x) \geq \mu$ and is \magenta{nondecreasing} with $  a(x) \leq a(\beta) = \mu(2 + \delta \beta)$. Since $b(x)/\mu$ is a function that depends only on $\beta$, the right-hand side is a quantity that increases in $\delta$. In \eqref{fEC4} we assumed $\delta < 1/2$, so $ 1/\kappa \leq \sup_{ \delta \in (0,1/2)} 1/\kappa \leq C/\mu$. Combining this with the fact that   $a(x) \leq  C\mu$ from  \eqref{fEC6} we get $2\kappa/a(x) \geq C$. Therefore, 
\bes{
\P(Y_1\leq -z) \geq &\  C \int_{-z-1}^{-z}   e^{\int_{0}^{y} \frac{-2u}{2 + \delta u} du} dy =   C\int_{-z-1}^{-z}   e^{ \frac{4\log(\delta y+2) -\magenta{2} \delta y}{\delta^2} -\frac{4\log(2)}{\delta^2}} dy .
}
Taylor expansion tells us that 
\begin{align}
\frac{4\log(\delta y+2) - \magenta{2}\delta y}{\delta^2} -\frac{4\log(2)}{\delta^2} = \frac{-y^2}{2} + \frac{\delta y^3}{6} - \frac{\delta^2 \xi^4}{16} \label{eq:logtaylor}
\end{align}
for some $\xi$ between $0$ and $y$. Recall that $\delta = 1/\sqrt{R}$ and $z\leq c_1 R^{1/4}$, meaning $\delta^2 \abs{\xi}^4 \leq C$ when $y \in [-z-1,-z]$, so
\bes{
 \P(Y_1\leq -z)  \geq C\int_{-z-1}^{-z} e^{-\frac{y^2}{2}+\frac{\delta y^3}{6}}dy =&\ C\int_{-1}^0 e^{-\frac{(-z+y)^2}{2}+\frac{\delta(-z+y)^3}{6}}dy.
}
Now for $y \in [-1,0]$, and $z\leq c_1 R^{1/4}$ we use the fact that $\delta z^2 \leq C$ to get
\begin{align*}
-\frac{(-z+y)^2}{2}+\frac{\delta(-z+y)^3}{6} =&\ -\frac{z^2}{2}+zy - \frac{y^2}{2}+ \delta\frac{z^2y}{2} - \delta \frac{z y^2}{2} + \delta \frac{y^3}{6} -\frac{\delta z^3}{6}\\
\geq&\ -\frac{z^2}{2}+zy - \frac{y}{2}(y - \delta z^2) -\frac{\delta z^3}{6} - C\\
\geq&\ -\frac{z^2}{2}+zy -\frac{\delta z^3}{6} - C,
\end{align*}
so 
\besn{\label{f30}
 \P(Y_1\leq -z)  \geq  C\int_{-1}^0 e^{-\frac{z^2}{2}+zy-\frac{\delta z^3}{6}}dy \geq  C \frac{1}{z}e^{-\frac{z^2}{2}-\frac{\delta z^3}{6}}.
}
This proves \eqref{eq:l1part1} for $k\geq 1$. For $k=0$,  the lemma follows from \eq{12} and \eq{f30}.  
\finishproof
\startproof{Proof of Lemma~\ref{fl3}}
\magenta{The bound \eq{14} trivially holds if $0<z<1$. Therefore, we assume $z\geq 1$ in the following.}
 Using integration by parts, we have for $k\geq 1$,
\besn{\label{15}
&\E|1(-z\leq W\leq 0)W^k e^{-\int_0^W r(u)du}|\\
=&\ -z^k e^{-\int_0^{\magenta{-z}} r(u)du} \P(W\leq -z)+\int_{-z}^0 \big(k(-y)^{k-1}+(-y)^kr(y)\big)e^{-\int_0^y r(u)du}\P(W\leq y) dy\\
\leq&\ C(k)\int_{-z}^0 \big((-y)^{k-1}+(-y)^{k+1}\big)e^{-\int_0^y r(u)du}\P(W\leq y) dy.
}
The last inequality is due to $r(y) = 2b(y)/a(y) = 2\mu(-y)/a(y) \leq C (-y)$ when $y \leq 0$, because $a(y) \geq \mu$. The same argument tells us that for $k=0$,
\ben{\label{16}
\E[1(-z\leq W\leq 0) e^{-\int_0^W r(u)du}]\leq \magenta{1+}C\int_{-z}^0 (-y)e^{-\int_0^y r(u)du}\P(W\leq y) dy.
}
In what follows we prove  that for $k\geq 0$,
\ben{\label{13}
\int_{-z}^0 (-y)^k e^{-\int_0^y r(u)du}\P(W\leq y) dy\leq Cz^{k\vee 1}.
}
 Lemma~\ref{fl3} follows from \eq{15}, \eq{16} and \eq{13}.
Without loss of generality assume that $z$ is an integer. If not, increase it to the nearest integer.
The taylor expansion in \eqref{eq:logtaylor} implies that for $0\leq -y\leq z \leq c_1 R^{1/4}$,
\be{
-\int_0^y r(u)du = -\int_0^y \frac{2b(u)}{a(u)} du =   -\frac{4\log(\delta y+2) - \magenta{2}\delta y}{\delta^2} + \frac{4\log(2)}{\delta^2} \leq \frac{y^2}{2}-\frac{\delta y^3}{6}+C.
}
Thus,
\bes{
&\int_{-z}^0 (-y)^k e^{-\int_0^y r(u)du}\P(W\leq y) dy\\
\leq &C\sum_{j=-z}^{-1} |j|^k \int_j^{j+1} e^{\frac{y^2}{2}-\frac{\delta y^3}{6}} e^{|j|y}e^{-|j|y} \P(W\leq y) dy\\
\leq &C \sum_{j=-z}^{-1} |j|^k \sup_{j\leq y\leq j+1} [e^{\frac{y^2}{2}+|j|y}] \sup_{j\leq y\leq j+1} [e^{-\frac{\delta y^3}{6}}]
\int_{j}^{j+1} e^{-|j|y} \P(W\leq y) dy\\
\leq & C \sum_{j=-z}^{-1} |j|^k e^{-\frac{j^2}{2}} \sup_{j\leq y\leq j+1} [e^{-\frac{\delta y^3}{6}}] \int_{-\infty}^{\infty} e^{-|j|y} \P(W\leq y) dy\\
=&C \sum_{j=-z}^{-1} |j|^k e^{-\frac{j^2}{2}} \sup_{j\leq y\leq j+1} [e^{-\frac{\delta y^3}{6}}] \frac{1}{|j|} E e^{-|j| W}.
}
We used $\sup_{j\leq y\leq j+1} [e^{\frac{y^2}{2}+|j|y}]=e^{-\frac{j^2}{2}+\frac{1}{2}}$ in the last inequality and integration by parts in the last equality. Invoking Lemma~\ref{l1}, we have
\bes{
\int_{-z}^0 (-y)^k e^{-\int_0^y r(u)du}\P(W\leq y) dy \leq&\  C \sum_{j=-z}^{-1} |j|^{k-1} e^{-\frac{j^2}{2}} \sup_{j\leq y\leq j+1} [e^{-\frac{\delta y^3}{6}}]  e^{\frac{j^2}{2}} e^{-\frac{\delta |j|^3}{6}}\\
\leq&\ C \sum_{j=-z}^{-1} |j|^{k-1} \leq C z^{k\vee 1},
}
where we used
\be{
\sup_{j\leq y\leq j+1} [-\frac{\delta y^3}{6} -\frac{\delta |j|^3}{6}]\leq C \delta j^2\leq C, 
\ \text{for}\  1\leq -j\leq z \leq c_1 R^{1/4}.
}
This proves \eq{13}.
\finishproof

\subsubsection{Proof of \eq{f13ec}.}
\startproof{Proof of \eq{f13ec}.} 
Inequality \eqref{f13ec} contains an upper bound on $\left|\frac{\Prob(Y_1\geq z)}{\Prob(W\geq z)}-1\right|$. 
\magenta{A} similar upper bound on $\left|\frac{\Prob(W\geq z)}{\Prob(Y_1\geq z)}-1\right|$ is a consequence of Theorem~4.1 of \cite{Brav2017}. The following simple modification of the argument in \cite{Brav2017} implies \eq{f13ec}. It follows from (4.8), (4.9), (4.11), and (4.12) of \cite{Brav2017} that there exist some $c_1,C_1 > 0$ such that for $R\geq C_1$ and $0<z\leq c_1 R$,
\bes{
&|P(Y_1\geq z)-P(W\geq z)|\\
\leq& C \delta^2 \P(Y_1\geq z)+C\delta^2 \P(W\geq z)+C\delta^2 \min\{(z\vee 1), \frac{1}{\delta^2}\}\P(Y_1\geq z)+C\delta \P(W\geq z).
}
Dividing both sides by $\P(W\geq z)$ we have 
\be{
\left|\frac{\P(Y_1\geq z)}{\P(W\geq z)}-1 \right|\leq C\delta^2 (z\vee 1) \frac{\P(Y_1\geq z)}{\P(W\geq z)}+C\delta.
}
Note that $C\delta^2 (z\vee 1) \leq C(\delta^2 z + \delta^2) = C(z/R + 1/R) \leq C(c_1 + 1/C_1) $, and choose $c_1$ small enough and $C_1$ large enough so that $C(c_1 + 1/C_1) < 1/2$. Then, repeating the argument used below \eqref{fEC5} implies \magenta{\eq{f13ec} for $R\geq C_1$, and the argument above \eq{fEC4} implies \eq{f13ec} for $R<C_1$.}
\finishproof

\subsection{Proof of Theorem~\ref{thm:md-std}}\label{fap3}
We first recall Theorem~\ref{thm:md-std}.
\begin{theorem} 
\label{thm:md-stdec}
Assume $n = R + \beta \sqrt{R}$ for some fixed $\beta > 0$. \blue{There exist positive constants $c_0$ and $C$ depending only on $\beta$ such that 
  \begin{align}
  & \left|\frac{\Prob(Y_0\geq z)}{\Prob(W\geq z)}-1\right|\leq \frac{C}{\sqrt{R}}\left(1+z\right) \quad \text{ for } 0<z\leq c_0 R^{1/2}\ \text{and} \label{f12ec} \\
&\left|\frac{\Prob(Y_0\leq -z)}{\Prob(W\leq -z)}-1\right|\leq \frac{C}{\sqrt{R}}\left(1+z^3\right),\ \text{ for } 0<z\leq \min\{c_0 R^{1/6}, R^{1/2}\}. \label{f15ec}
  \end{align}
}
\end{theorem}
Theorem~\ref{thm:md-stdec} follows from a similar and simpler proof than that of Theorem~\ref{thm:md-highec}. 
In particular, the bound \eq{f15ec} can be proved by a simple adaption of the arguments in \cite{ChFaSh13a} and therefore its proof is omitted.
The proof of \eq{f12ec} below may be useful for other exponential approximation problems.
We \magenta{will use $c_0,\ C,\ C_0$ to denote positive constants} that may differ from line to line, but will only depend on $\beta$. We require the following two lemmas. The first one is proved in Section~\ref{fap7}, and the second in Section~\ref{fap6}.
\begin{lemma}\label{l2}
There exist    $C,C_0 > 0$ that depend  only on $\beta$, such that for any $R \geq C_0$ and any $z > \beta$, 
\ben{\label{5}
\P(W\geq \magenta{z+\delta})\leq C e^{-(\beta- 3\beta^2 \delta)z}.
}
\end{lemma}
To state the second lemma we let $f_z(x)$ solve the Poisson equation
\ben{\label{8}
b(x) f_z'(x) + \mu f_z''(x)= 1(x\geq z) - \P(Y_0\geq z).
}
\begin{lemma}\label{lem:eq9}
There exist $C_0,C > 0$ depending on $\beta$ such that for all $R \geq C_0$ \magenta{and any $z>\beta$}, 
\begin{align}
&\E \Big( \int_{-\delta}^{\delta} \big( b(W+y)f_z'(W+y) - b(W) f_z'(W)\big) \frac{K_W(y)}{\mu } dy \Big) \notag \\ \leq&\ C\delta  \P(Y_0\geq z) \E \Big( \sup_{|s|\leq \delta}\big(1+e^{\beta (W+s)}1(\beta \leq W+s \leq z)\magenta{\big)\Big)}, \label{eq:lem91}
\end{align}
where $K_W(y)$ is defined in \eqref{eq:kdef}. Furthermore, 
\begin{align}
 & - \frac{\delta}{2\mu } \E \big( b^2(W) f_z'(W)\big) \leq  C\delta \Big(\P(W\geq z)+\P(Y_0\geq z)\big(1+\E e^{\beta W}1(\beta\leq W\leq z) + |W|\big) \Big), \label{eq:lem92}\\
&- \frac{\delta}{2\mu } \E \big( b(W)\big( \Prob(Y_0 \geq z)   -1(W \geq \delta + z) \big)\big) \leq  C\delta  \P(W\geq z+\delta)\magenta{.} \label{eq:lem93}
\end{align}  
\end{lemma} 
\startproof{Proof of \eq{f12ec}}
Note that if $R$ is bounded, then, for fixed $\beta$, $n$ is also bounded and $P(W\geq z)$ is bounded away from 0 for $z$ in the bounded range $0<z\leq c R^{1/2}$. By choosing a sufficiently large $C$, \eq{f12ec} trivially holds. Therefore, in the following, we assume $R\geq C_0$ for a sufficiently large $C_0$. Additionally, the result for finite range $z \in (0, \beta+\delta]$ follows from the Berry-Esseen bound in Theorem 3 of \cite{BravDaiFeng2016}, so we assume $z > \beta + \delta$. 

Recall that the density of $Y_0$ is given in \eqref{eq:stddenf}, and that $b(x) =-(\mu x \wedge \mu \beta)$. 
Just like in \eqref{f18}, \magenta{one may verify that the solution to the Poisson equation \eq{8} satisfies}
\ben{\label{f33}
f_z'(x)=
\begin{cases}
-\P(Y_0\geq z)e^{-\int_{0}^{x} \frac{b(u)}{\mu}  du} \int_{-\infty}^{x} \frac{1}{\mu } e^{\int_{0}^{y} \frac{b(u)}{\mu}  du} dy , & x<z, \\
-\P(Y_0\leq z)e^{-\int_{0}^{x} \frac{b(u)}{\mu}  du} \int_{x}^{\infty}  \frac{1}{\mu } e^{\int_{0}^{y} \frac{b(u)}{\mu}  du} dy, & x\geq z,
\end{cases}
}
\magenta{so} \eq{eq:gxbarf} implies  $\E G_{\magenta{\tilde X}} f_z(W) = 0$. Taking expected values in \eqref{8} then gives us 
\begin{align*}
 &\P(W\geq z) - \P(Y_0\geq z)  \\
 =&\ \E\Big( b(W) f_z'(W) + \mu f_z''(W)\Big) - \E G_{\magenta{\tilde X}} f_z(W) \\
=&\ \E\bigg(\mu f_z''(W)- \int_{-\delta}^{\delta}f_{\magenta{z}}''(W+y)  K_W(y) dy \bigg) \\
=&\ \E\big( 1(W\geq z) - \P(Y_0\geq z)  -b(W)f_z'(W)\big)\\
&-\E\bigg(\int_{-\delta}^{\delta} \big( 1(W+y\geq z) - \P(Y_0\geq z) - b(W+y) f_z'(W+y)\big) \frac{1}{\mu}K_W(y) dy\bigg),
\end{align*} 
where $K_W(y)$ is defined in \eqref{eq:kdef}, the second equality is due to \eqref{eq:steinx}, and the final equality follows from  $\mu f_z''(x) = - b(x) f_z'(x) + 1(x \geq z) - \Prob(Y_0 \geq z)$. Using $K_W(y) \geq 0$ and recalling from \eqref{MD:k3} that $\int_{-\delta}^{\delta} K_W(y) dy =  \mu - \frac{\delta}{2}b(W)$, we have 
\begin{align*}
&b(W)f_z'(W) = \int_{-\delta}^{\delta}b(W) f_z'(W) \frac{K_W(y)}{\mu } dy + \frac{\delta}{2\mu } b^2(W) f_z'(W),  \\
&-\E\bigg(\int_{-\delta}^{\delta}1(W+y\geq z)  \frac{1}{\mu}K_W(y) dy\bigg) \leq -\E\Big(1(W \geq z+\delta) \big(1 - \frac{\delta}{2 \mu} b(W)\big) \Big),\\
&\E\bigg(\int_{-\delta}^{\delta}  \P(Y_0\geq z) \frac{1}{\mu}K_W(y) dy\bigg)  = \P(Y_0\geq z) \E\Big(  1 - \frac{\delta}{2 \mu} b(W) \Big),
\end{align*}
and therefore, 
\begin{align}
&\P(W\geq z) - \P(Y_0\geq z)  \notag \\
\leq&\ \P(z \leq W < z+\delta) + \E \Big( \int_{-\delta}^{\delta} \big( b(W+y)f_z'(W+y) - b(W) f_z'(W)\big) \frac{K_W(y)}{\mu } dy \Big) \notag \\
&- \frac{\delta}{2\mu } \E \big( b^2(W) f_z'(W)\big) - \frac{\delta}{2\mu } \E \big( b(W)\big( \Prob(Y_0 \geq z)   -1(W \geq \delta + z) \big)\big)\magenta{.} \label{10}
\end{align}
Subtracting $\P(z \leq W < z+\delta)$ from both sides and using Lemma~\ref{lem:eq9} to bound the remaining three terms,  we get 
\begin{align*}
&\P(W\geq z+\delta) - \P(Y_0\geq z) \notag \\
\leq&\ C\delta\sup_{|s|\leq \delta}\Big( \P(\magenta{W}\geq z)+  \P(Y_0\geq z)\big(1+\E|\magenta{W}|+\E e^{\beta(W+s)}1(\beta\leq W+s\leq z)\big) \Big).
\end{align*}
We now argue that the right-hand side above can be bounded by $ C\delta \P(Y_0\geq z) (z\vee 1)$. Note that $\E \abs{W} \leq C$ due to Lemma~\ref{lem:higherders}. Next, since $b(x)/\mu = -(x \wedge \beta)$ only depends on $\beta$ and $z \geq \beta$, then $\Prob(Y_0 \geq z) =  \frac{\int_{z}^{\infty} e^{-\beta y} dy}{\int_{-\infty}^{\infty} e^{\int_{0}^{y} b(u)/\mu du }}  = C e^{-\beta z}$ for some $C>0$ depending only on $\beta$. Thus, Lemma~\ref{l2}  tells us that  for \magenta{$R\geq C_0$ and $\beta+\delta<z\leq c_0 R^{1/2}$},
\be{
\P(\magenta{W}\geq z) = \frac{ \P(\magenta{W}\geq z)}{\P(Y_0\geq z)}\P(Y_0\geq z)\leq C\P(Y_0\geq z).
}
Moreover, \magenta{for $|s|\leq \delta$,}
\bes{
&\E \big(e^{\beta(W+s)}1(\beta\leq W+s\leq z)\big)\\
&\magenta{= e^{\beta^2}\P(\beta\leq W+s\leq z) +\int_\beta^z \beta e^{\beta y} \P(y<W+s\leq z) dy } \\
&\leq C+\beta\int_{\beta}^z e^{\beta u} \P(W+s\geq u)du\\
&\leq C+ \int_{\beta}^z C du\\
&\leq C(z\vee 1).
}
The \magenta{first} equality follows from integration by parts, and the second-last inequality is due to Lemma~\ref{l2}. Thus we have shown that
\be{
\P(W\geq z+\delta)-\P(Y_0\geq z)\leq C\delta \P(Y_0\geq z) (z\vee 1).
} 
Note that \magenta{$|\P(Y_0\geq z)/\P(Y_0\geq z+\delta)-1|=|e^{\beta \delta}-1|\leq \delta C$}. Therefore, 
\be{
\P(W\geq z+\delta)-\P(Y_0\geq z+\delta)(1+\delta C)\leq C\delta \P(Y_0\geq z+\delta)(1+\delta C)(z\vee 1).
}
Dividing both sides of the above inequality by $\P(W\geq z+\delta)$, we obtain
\be{
1-\frac{\P(Y_0\geq z+\delta)}{\P(W\geq z+\delta)}(1+\delta C)\leq C\delta \frac{\P(Y_0\geq z+\delta)}{\P(W\geq z+\delta)}(1+\delta C)(z\vee 1).
}
Choose $C_0$ large enough and $c_0$ small enough and \magenta{since} $R\geq C_0$ and $0<z\leq c_0 R^{1/2}$, the constant in front of  $\frac{\P(Y_0\geq z+\delta)}{\P(W\geq z+\delta)}$ on the right-hand side can be made less than \magenta{1/2, implying }
\be{
1-\frac{\P(Y_0\geq z+\delta)}{\P(W\geq z+\delta)} \leq C\delta  (z\vee 1).
}
A similar argument can be repeated to show $\frac{\P(Y_0\geq z+\delta)}{\P(W\geq z+\delta)} -1 \leq  C\delta  (z\vee 1)$, implying \eqref{f12ec}.
\finishproof

\subsubsection{Proof of Lemma~\ref{l2}.}\label{fap7}
We require the following auxiliary result, whose proof is provided after the proof of Lemma~\ref{l2}.
\begin{lemma}\label{l3}
There exist constants $C,C_0 > 0$ depending only on $\beta$ such that for $R \geq C_0$,  
\ben{\label{6}
\E e^{(\beta-3\beta^2 \delta)W}\leq C/\delta.
}
\end{lemma}

\startproof{Proof of Lemma~\ref{l2}}
Since $\delta = 1/\sqrt{R}$, we can choose $C_0$ large enough   so that $(\beta-3\beta^2 \delta)>\beta/2$. For notational convenience, we set  $\nu = 3\beta^2 \delta$. Define
\be{
f''(x)=
\begin{cases}
(\beta-\nu)e^{(\beta-\nu)x}, & x\leq z\\
\text{linear interpolation}, & z<x\leq z+\Delta\\
0, & x>z+\Delta,
\end{cases}
}
with $0<\Delta\leq \delta$ to be chosen,
$f'(x) = 1 + \int_{0}^{x} f''(y) dy$ and $f(x) = \int_{0}^{x} f'(y) dy$. Note that $f(x)$ grows  linearly in $x$  when $x \geq z + \Delta$ because $f'(x) = f'(z+\Delta)$ for $x \geq z + \Delta$, so $\E \abs{f(W)} < \infty$ because $W$ is bounded from below and $\E \abs{W} < \infty$. Therefore,  $\E G_{\magenta{\tilde X}} f(W) = 0$ due to \eqref{eq:gxbarf}, implying $ \E\big(-b(W)f'(W)\big) = \E \Big(\int_{-\delta}^{\delta} f''(W+y) K_{W}(y) dy\Big)$ if we use the form of $G_{\magenta{\tilde X}} f(x)$ in \eqref{eq:steinx}. Note also that $f'(0) = 1$, $f'(x) \geq 0$, and  $f'(x) \geq e^{(\beta-\nu)z}$ for $x \geq z + \Delta$. Using $b(x) = -(\mu x \wedge \mu \beta)$ \magenta{and the assumption that $z>\beta$} we therefore have
\bes{
 \E\big(-b(W)f'(W)\big) \geq&\ \E\big(\mu We^{(\beta-\nu)W}1(W<\beta)\big) \\
&\quad + \E\big(\mu \beta e^{(\beta-\nu)W}1(\beta\leq W\leq z)\big) +\E\big(\mu \beta e^{(\beta-\nu)z}1(W> z+\Delta)\big)\\
 \geq&\ -\mu C+ \E\big(\mu \beta e^{(\beta-\nu)W}1(\beta\leq W\leq z)\big) +\E\big(\mu \beta e^{(\beta-\nu)z}1(W> z+\Delta)\big),
}
where the second inequality is due to $-\abs{x} e^{-(\beta-\nu)\abs{x}} 1(x < \beta) \geq -C$. Recalling from \eqref{MD:k3} that $\int_{-\delta}^{\delta} K_{W}(y)dy = \mu-\delta b(W)/2 $, we have 
\bes{
\E \Big(\int_{-\delta}^{\delta} f''(W+y) K_{W}(y) dy\Big) &\leq\E\Big(  \sup_{|s|\leq \delta} f''(W+s) \int_{-\delta}^{\delta} K_{W}(y)dy \Big)\\
&=\E \Big( \mu \sup_{|s|\leq \delta} f''(W+s)  \big(1-\frac{\delta}{2\mu}b(W)\big)1(W\leq z+\Delta+\delta)\Big)\\
&\leq \E\Big( \mu (\beta-\nu) e^{(\beta-\nu)(W+\delta)}\big(1-\frac{\delta}{2\mu}b(W)\big)1(W\leq z+\Delta+\delta)\Big)\\
&=\E\Big( \mu (\beta-\nu) e^{(\beta-\nu)(W+\delta)}\big(1+\frac{\delta}{2}\beta\big)1(\beta\leq W\leq z+\Delta+\delta)\Big)\\
&\quad +\E\Big( \mu (\beta-\nu) e^{(\beta-\nu)(W+\delta)}\big(1+\frac{\delta}{2}W\big)1(W<\beta )\Big)\\
&\leq \E\Big( \mu (\beta+C\delta) e^{(\beta-\nu)W}1(\beta\leq W\leq z)\Big)\\
&\quad +\E \Big( \mu (\beta+C\delta)e^{(\beta-\nu)z}1(z\leq W\leq z+\Delta+\delta)\Big)+ \mu C.
}
Combining the inequalities above, we have
\bes{
&\E\big(\beta e^{(\beta-\nu)z}1(W> z+\Delta+\delta)\big)\\
&\leq C+C\delta \E \big( e^{(\beta-\nu)W}1(\beta\leq W\leq z+\Delta+\delta)\big)+\beta e^{(\beta-\nu)z} 
\P(z\leq W\leq z+\Delta).
}
Without loss of generality we assume $z$ does not belong to the support of $W$ and let $\Delta\to 0$, and observe that $\P(z\leq W\leq z+\Delta) \to 0$. 
Therefore, we have
\be{
\magenta{\beta e^{(\beta-\nu)z}}\P(W>z+\magenta{\delta})\leq C+ C\delta \E e^{(\beta-\nu)W}\leq C
}
where we \magenta{have used Lemma~\ref{l3}.}  
\finishproof

\startproof{Proof of Lemma~\ref{l3}} 
Since $\delta = 1/\sqrt{R}$, we can choose $C_0$ large enough   so that $(\beta-3\beta^2 \delta)>\beta/2$. For notational convenience, we set  $\nu = 3\beta^2 \delta$. 
Fix $M >  \beta $ and let  $f(x) = \int_{0}^{x}e^{(\beta-\nu)(y  \wedge M)} dy$. We recall from \eqref{eq:steinx} that $G_{\magenta{\tilde X}} f(x) = b(W)f'(W) +\int_{-\delta}^{\delta} f''(W+y) K_{W}(y) dy$ where $K_W(y)$ is defined in \eqref{eq:kdef}. Since $\beta - \nu > \beta/2$ by assumption, the function $f(x)$ grows linearly for $x \geq M$, so $\E \abs{ f(W)} < \infty$ because $\E \abs{W} < \infty$. Therefore,  $\E G_{\magenta{\tilde X}} f(W) = 0$, or $ \E\big(-b(W)f'(W)\big) = \E \Big(\int_{-\delta}^{\delta} f''(W+y) K_{W}(y) dy\Big)$, due to \eqref{eq:gxbarf}. Now
\besn{\label{f40}
 \E\big(-b(W)f'(W)\big)=&\  \E\big(\mu \beta e^{(\beta-\nu)(W \wedge M)} 1(W\geq \beta) \big)+\E\big(\mu W e^{(\beta-\nu)W} 1(W< \beta) \big)\\
=&\ \E\big(\mu \beta e^{(\beta-\nu)(W \wedge M)}  \big)+\E\big(\mu (W-\beta) e^{(\beta-\nu)W} 1(W< \beta)\big)\\
\geq&\ \E\big(\mu \beta e^{(\beta-\nu)(W \wedge M)}  \big)-\mu C
}
where in the last inequality we used the fact that $\abs{(x-\beta) e^{(\beta-\nu)x}} 1(x<\beta) \leq C$  if $(\beta-\nu)>\beta/2$. Furthermore, since $f''(x) = (\beta-\nu) e^{(\beta-\nu)(x \wedge M)} 1(x < M)$ and $\int_{-\delta}^{\delta} K_{W}(y)dy = \mu -\delta b(W)/2 $, we have
\besn{\label{f41}
\E \bigg(\int_{-\delta}^{\delta} f''(W+y) K_{W}(y) dy\bigg) \leq&\ \E \bigg(\int_{-\delta}^{\delta}(\beta-\nu)e^{(\beta-\nu) ((W+y)\wedge M ) }K_W(y) dy\bigg)\\
\leq&\   (\beta-\nu) e^{(\beta-\nu)\delta}\E\Big(e^{(\beta-\nu)(W \wedge M)}(\mu-\frac{\delta}{2}b(W))\Big)\\
=&\  \mu (\beta-\nu) e^{(\beta-\nu)\delta}\E\Big(e^{(\beta-\nu)(W \wedge M)}(1+\frac{\delta}{2}W)1(W \leq \beta)\Big) \\
&+ \mu (\beta-\nu) e^{(\beta-\nu)\delta}\E\Big(e^{(\beta-\nu)(W \wedge M)}(1+\frac{\delta}{2}\beta)1(W \geq \beta)\Big)\\
\leq&\  \mu C+ \mu (\beta-\nu)e^{(\beta-\nu)\delta}(1+\frac{\delta}{2}\beta)\E\big(e^{(\beta-\nu)(W \wedge M)}\big).\\
}
Divide both sides of \eq{f40} and \eq{f41} by \magenta{$\mu \delta$} and combine these two inequalities, and also substitute $3 \beta^2 \delta$ for $\nu$, to get
\be{
\frac{\Big(\beta-(\beta-3\beta^2\delta)e^{(\beta-3\beta^2\delta)\delta} (1+\frac{\delta}{2}\beta) \Big)}{\delta} \E \big( e^{(\beta-3 \beta^2 \delta)(W\wedge M)}\big)\leq C/\delta.
}
Since the coefficient in front of the expected value on the left-hand side converges to a positive constant as $\delta \to 0$,  for sufficiently small $\delta$ (or sufficiently large $C_0$), we have 
\be{
\E \big(e^{(\beta-3\beta^2\delta)(W\wedge M)}\big)\leq C/\delta.
}
We conclude by letting $M \to \infty$.
\finishproof

\subsubsection{Proof of Lemma~\ref{lem:eq9}.}\label{fap6}
We begin by proving \eqref{eq:lem91}. 
\begin{align*}
&\E \Big( \int_{-\delta}^{\delta} \big( b(W+y)f_z'(W+y) - b(W) f_z'(W)\big) \frac{K_W(y)}{\mu } dy \Big) \\
=&\ \E \Big( \int_{-\delta}^{\delta}\frac{K_W(y)}{\mu } \int_{0}^{y}\big( b(x)f_z'(x)\big)'\big|_{x=W+s} ds dy \Big) \\
\leq&\ \delta \E \Big( \sup_{|s|\leq \delta} \abs{\big( b(x)f_z'(x)\big)'\big|_{x=W+s}} \big( 1 - \frac{\delta}{2\mu}b(W) \big) \Big)\\
\leq&\ C \delta \E \Big( \sup_{|s|\leq \delta} \abs{\big( b(x)f_z'(x)\big)'\big|_{x=W+s}}  \Big).
\end{align*}
The first inequality is due to $K_W(y) \geq 0$ and $\int_{-\delta}^{\delta} K_{W}(y)/\mu dy =  1-\delta b(W)/(2\mu)$ from \eqref{MD:k3}, and the last inequality  is true because $\big| \frac{\delta}{2\mu}b(W)\big| = \big|\frac{\delta}{2}(W \wedge \beta) \big| \leq C$ since $W\geq -1/\delta$. To bound the right-hand side we note that \magenta{(cf. \eq{f33})}
\be{
-b(x)f_z'(x)=
\begin{cases}
\P(Y_0\geq z)b(x)  e^{-\int_{0}^{x} \frac{b(u)}{\mu}  du} \int_{-\infty}^{x}  \frac{1}{\mu } e^{\int_{0}^{y} \frac{b(u)}{\mu}  du} dy , & x<z, \\
\P(Y_0\leq z)b(x) e^{-\int_{0}^{x} \frac{b(u)}{\mu}  du} \int_{x}^{\infty} \frac{1}{\mu }  e^{\int_{0}^{y} \frac{b(u)}{\mu}  du} dy, & x\geq z,
\end{cases}
}
so for \magenta{$x> z>\beta$},
\bes{
(-b(x)f_z'(x))' =&\ -\beta \P(Y_0\leq z)\Big(-1+  \beta e^{-\int_{0}^{x} \frac{b(u)}{\mu}  du}\int_{x}^{\infty}    e^{ \int_{0}^{y} \frac{b(u)}{\mu}  du} dy \Big)\\
 =&\ -\beta \P(Y_0\leq z)\Big(-1+  \beta e^{- \frac{\beta^2}{2} + \beta x}\int_{x}^{\infty}    e^{ \frac{\beta^2}{2} - \beta y} dy \Big)\\
 =&\ 0.
}
For $\beta<x<z$,  
\bes{
\abs{(b(x)f_z'(x))'} =&\ \beta \P(Y_0\geq z) \Big(1+ \beta e^{-\frac{\beta^2}{2}+\beta x} \int_{-\infty}^{x}  e^{\int_{0}^{y} \frac{b(u)}{\mu}  du} dy   \Big) \leq \beta \P(Y_0\geq z) (1 + C e^{\beta x}).
}
In the inequality above we used the fact that $\int_{-\infty}^{x}  e^{\int_{0}^{y} \frac{b(u)}{\mu}  du} dy \leq C$  because $b(x)/\mu = -(x \wedge \beta)$ depends only on $\beta$. Lastly, for $x<\beta$,
\bes{
\abs{(b(x)f_z'(x))'} =&\  \P(Y_0\geq z)\Big|x+(e^{\frac{x^2}{2}}+x^2 e^{\frac{x^2}{2}})\int_{-\infty}^{x}  e^{- \frac{y^2}{2}} dy \Big|.
}
When $ -1 \leq x\magenta{ < \beta}$, the right-hand side is bounded by $C \P(Y_0\geq z)$ and when $x < -1$, we use the bound $\frac{1}{-x-1} e^{-\frac{x^2}{2}}\leq \int_{-\infty}^{x}  e^{- \frac{y^2}{2}} dy \leq  \frac{1}{-x} e^{-\frac{x^2}{2}}$ to conclude that 
\begin{align*}
\P(Y_0\geq z)\Big|x+(e^{\frac{x^2}{2}}+x^2 e^{\frac{x^2}{2}})\int_{-\infty}^{x}  e^{- \frac{y^2}{2}} dy \Big| \leq \P(Y_0\geq z)\Big|x+C-x  \Big| \leq C \P(Y_0\geq z).
\end{align*}
Combining the three cases yields \eqref{eq:lem91}. We now prove the bound on $- \frac{\delta}{2\mu } \E \big( b^2(W) f_z'(W)\big)$ in \eqref{eq:lem92}. From the form of $b(x) f_z'(x)$ above,  we have for $x>z\magenta{>\beta}$,
\bes{
\magenta{-}\frac{1}{\mu}b^2(x)f_z'(x)&= \magenta{\beta} \P(Y_0\leq z),
}
for $\beta\leq x\leq z$,
\bes{
\magenta{-}\frac{1}{\mu}b^2(x)f_z'(x)&=  \P(Y_0\geq z) \beta e^{-\frac{\beta^2}{2}+\beta x} \int_{-\infty}^{x}  e^{\magenta{\int_{0}^{y}} \frac{b(u)}{\mu}  du} dy 
\leq C  e^{ \beta x}\P(Y_0\geq z),
}
and for $x<\beta$,
\be{
\magenta{-}\frac{1}{\mu}b^2(x)f(x) =  \P(Y_0\geq z)  x^2 e^{\frac{x^2}{2}}\int_{-\infty}^{x}  e^{- \frac{y^2}{2}} dy   \leq  \P(Y_0\geq z) \magenta{(1+|x|)}.
}
Combining the three cases implies \eqref{eq:lem92}. Lastly we prove \eqref{eq:lem93}. In the proof of Lemma~\ref{lem:momentequiv} we showed that  $\E b(W) = 0$. Furthermore,    $-\mu\beta\leq b(x)\leq 0$ for $x\geq 0$, so
\begin{align*}
 -\frac{\delta}{2\mu} \E \big( b(W) \big(\P(Y_0\geq z)-1(W\geq z+\delta)\big) \big)=  \frac{\delta}{2\mu} \E \big( b(W)1(W\geq z+\delta) \big) \leq  C\delta  \P(W\geq z+\delta).
\end{align*} 
\finishproof

\section{Companion for the Hospital Model}
\label{app:hospital_proofs}
In this portion of the electronic companion, we motivate the $v_3$ approximation for the hospital model presented in Section~\ref{sec:hospcompare}, where we suggested  using 
\begin{align}
  \label{eq:hosv3ec}
   v_3(x)=\max\Big\{ \delta +\frac{1}{2} \Big(\delta^{2} 1(x<0)- \delta^2 ( x^{-} - \beta )-\delta^2-2\delta^2\beta\Big),  \delta/2  \Big\}.
\end{align}
We recall that $\tilde X = \{ \tilde X(n) = \delta (X(n) - N)\}$, where $X(n)$ is the customer count at the end of time unit $n$, that  $W$ and $W'$ have the distributions of $\tilde X(0)$ and $\tilde X(1)$ when $\tilde X(0)$ is initialized according to the stationary distribution of $\tilde X$, and that $\Delta = W' - W$.  We use $\epsilon(x)$ and $\epsilon_i(x)$ to denote generic functions satisfying
\begin{align}
\abs{\epsilon(x)} \leq C (1 + \abs{x})^{5}. \label{eq:epsdefec}
\end{align}
The following lemma gives us the conditional moments of $\Delta$. It is proved in Section~\ref{sec:hosmecproof}.
\begin{lemma}\label{lem:hosmec}
  For the hospital model with $N$ servers,  $\Lambda=\sqrt{N}-\beta$ and $\mu=\delta=1/\sqrt{N}$,
  \begin{align}
 b(x) = \E(\Delta| W= x) =&\  \delta( x^{-} - \beta ), \label{eq:1ec}\\ 
   \E(\Delta^2|W=x)  =&\ 2\delta + \Big(b^{2}(x)-\delta b(x) -\delta^2 - 2\delta^2\beta\Big) + \delta^3 x^{-}  \label{eq:2ec} \\ 
   \E(\Delta^3|W=x) =&\ 6\delta b(x)  + \delta^{3} \epsilon(x) , \label{eq:3ec}\\ 
    \E(\Delta^4|W=x) =&\  12\delta^2  + \delta^{3} \epsilon(x), \label{eq:4ec}\\
   \E(\Delta^5|W=x)=&\  \delta^{3} \epsilon(x). \label{eq:5ec}
  \end{align}
\end{lemma}
To derive $v_3(x)$, we begin with the Taylor expansion in \eqref{eq:taylorgeneric} with $n = 4$, which says that for sufficiently smooth $f(x)$, 
\begin{align}
- \E \Delta f'(W) =  \E \bigg[\sum_{i=2}^{4} \frac{1}{i!}  \Delta^i f^{(i)}(W) + \frac{1}{5!} \Delta^{5} f^{(5)}(\xi_1)\bigg], \label{eq:txec}
\end{align}
where $\xi_i$ denote numbers lying between $W$ and $W'$. By combining \eqref{eq:txec} with Lemma~\ref{lem:hosmec}, we will show that
\begin{align}
&- \E b(W) f'(W) - \frac{1}{2} \E\Big(2 \delta + \delta^{2} 1(W<0)- \delta^2( W^{-} - \beta ) -\delta^2-2\delta^2\beta\Big) f''(W)\notag \\
=&\    \delta^3 \Big(\frac{1}{2}\E \epsilon_0(W)  f''(W)+ \frac{1}{6}  \E   \epsilon_3(W) f'''(W) + \frac{1}{24}  \E \epsilon_4(W)  f^{(4)}(W) + \frac{1}{120}   \E \epsilon_5(W) f^{(5)}(\xi_1)\Big) \notag \\ 
&+ \frac{1}{2}\delta^3 \E \Big(\epsilon_1(W) f^{(4)}(W) +  \epsilon_2(W) f^{(5)}(\xi)\Big) \notag \\
&+\frac{1}{2} \delta^3 \E \Big(\epsilon_6(W)f''(W)+\epsilon_7(W)f'''(W)  +  \epsilon_2(W) \Big(\frac{d^2}{d x^2} \big(   (x^{-}-\beta) f''(x) \big)\big|_{x = \xi_3}\Big)\Big). \label{eq:toprovehospec}
\end{align}
Truncating the term in front of $f''(W)$ on the left-hand side from below by $\delta/2$ gives us $v_3(x)$ in \eqref{eq:hosv3ec}.  The truncation level $\delta/ 2$ is chosen because the support of $W$ is in $[-\delta N, \infty)$ and the term in front of $f''(W)$ on the left-hand side of \eqref{eq:toprovehospec} equals $\delta \big(\frac{1}{2} - \frac{1}{2}\delta  \beta \big) \approx\frac{\delta}{2}$ when evaluated at the point $W = -\delta N$. We could have chosen $\delta \big(\frac{1}{2} - \frac{1}{2}\delta  \beta \big)$ instead (when this quantity is positive), but in practice this does not make a big difference.  Let us now prove \eqref{eq:toprovehospec}. Combining \eqref{eq:txec} with Lemma~\ref{lem:hosmec} yields
\begin{align*}
- \E b(W) f'(W) =&\ \frac{1}{2} \E\Big(2\delta + b^{2}(W)-\delta b(W) -\delta^2 - 2\delta^2\beta + \delta^3 W^{-}\Big) f''(W) \notag \\
&+ \frac{1}{6} \E \big(6\delta b(W)  + \delta^{3} \epsilon_3(W)\big)f'''(W) + \frac{1}{24} \E\big(12\delta^2 + \delta^3 \epsilon_4(W) \big) f^{(4)}(W)  \notag \\
&+ \frac{1}{120} \delta^3 \E \epsilon_5(W) f^{(5)}(\xi_1).
\end{align*}
Let us write $W^{-}f''(W)$ as $\epsilon_{0}(W) f''(W)$  and rearrange the right-hand side above into the more convenient form:
\begin{align}
&- \E b(W) f'(W) - \frac{1}{2} \E\Big(2\delta + b^{2}(W)-\delta b(W) -\delta^2 - 2\delta^2\beta  \Big) f''(W) \notag \\
=&\     \frac{1}{2} \delta \E   b(W)  f'''(W) + \frac{1}{2}  \delta \big(\E   b(W)  f'''(W) + \delta  \E  f^{(4)}(W)\big)  \notag \\
&+\delta^3 \Big( \frac{1}{2} \E \epsilon_0(W)  f''(W)+ \frac{1}{6}  \E   \epsilon_3(W) f'''(W) + \frac{1}{24}  \E \epsilon_4(W)  f^{(4)}(W) + \frac{1}{120}   \E \epsilon_5(W) f^{(5)}(\xi_1)\Big). \label{1}
\end{align}
The last row is considered as error because of the $\delta^3$ there. We wish to transform the first row on the right-hand side into an expression involving $f''(x)$ plus error. To this end we require the following lemma, which is proved at the end of this section.
\begin{lemma}
\label{lem:hospaux}
Suppose that $g \in C^{3}(\R)$ is such that $\E g(W') - \E g(W) = 0$. Then
\begin{align*}
 \E b(W) g'(W) + \delta \E g''(W) =&\     \delta^2 \E \Big(\epsilon_1(W) g''(W) +  \epsilon_2(W) g'''(\xi)\Big),  
\end{align*} 
where $\epsilon_i(x)$ are generic functions satisfying \eqref{eq:epsdefec} and $\xi$ lies between $W$ and $W'$.
\end{lemma}
We apply Lemma~\ref{lem:hospaux} with $g(x) = \frac{1}{2} \delta f''(x)$ to get
\begin{align}
& \frac{1}{2}\delta \Big(\E b(W) f'''(W) + \delta \E f^{(4)}(W)\Big) = \frac{1}{2}\delta^3 \E \Big(\epsilon_1(W) f^{(4)}(W) +  \epsilon_2(W) f^{(5)}(\xi)\Big). \label{eq:gfdelta}
\end{align}
The left-hand side above coincides with one of the terms in the second row of \eqref{1}. 
Next we choose $g(x) = \int_{0}^{x} b(y) f''(y) dy$ and note that $g''(x) = b'(x)f''(x) + b(x)f'''(x)$. Applying Lemma~\ref{lem:hospaux} with our new choice of $g(x)$, we get 
\begin{align*}
& \E b^2(W) f''(W) + \delta \E \big(b'(W)f''(W) + b(W)f'''(W)\big) \\
 =&\  \delta^2 \E \Big(\epsilon_1(W) \big(b'(W)f''(W) + b(W)f'''(W)\big) +  \epsilon_2(W) \Big(\frac{d^2}{d x^2} \big(   b(x) f''(x) \big)\big|_{x = \xi_3}\Big)\Big)\\
 =&\  \delta^3 \E \Big(\epsilon_6(W)f''(W)+\epsilon_7(W)f'''(W)  +  \epsilon_2(W) \Big(\frac{d^2}{d x^2} \big(   (x^{-}-\beta) f''(x) \big)\big|_{x = \xi_3}\Big)\Big).
\end{align*}
The last equality follows from the fact that $b(x) = \delta(x^{-}-\beta)$. Multiplying both sides by $1/2$ and rearranging terms, we get
\begin{align}
& \frac{1}{2}\delta \E b(W) f'''(W) \notag \\
=&\ - \frac{1}{2}\E\big( b^2(W) + \delta b'(W)\big) f''(W)  \notag  \\
&+\frac{1}{2} \delta^3 \E \Big(\epsilon_6(W)f''(W)+\epsilon_7(W)f'''(W)  +  \epsilon_2(W) \Big(\frac{d^2}{d x^2} \big(   (x^{-}-\beta) f''(x) \big)\big|_{x = \xi_3}\Big)\Big). \label{f25}
\end{align}
Plugging  \eqref{eq:gfdelta} and \eqref{f25} into \eqref{1}, we conclude that 
\begin{align*}
&- \E b(W) f'(W) - \frac{1}{2} \E\Big(2\delta + b^{2}(W)-\delta b(W) -\delta^2 - 2\delta^2\beta  \Big) f''(W)\notag \\
=&\    - \frac{1}{2}\E\big( b^2(W) + \delta b'(W)\big) f''(W) \\
&+ \delta^3 \Big(\frac{1}{2}\E \epsilon_0(W)  f''(W)+ \frac{1}{6}  \E   \epsilon_3(W) f'''(W) + \frac{1}{24}  \E \epsilon_4(W)  f^{(4)}(W) + \frac{1}{120}   \E \epsilon_5(W) f^{(5)}(\xi_1)\Big) \\ 
&+ \frac{1}{2}\delta^3 \E \Big(\epsilon_1(W) f^{(4)}(W) +  \epsilon_2(W) f^{(5)}(\xi)\Big) \\
&+\frac{1}{2} \delta^3 \E \Big(\epsilon_6(W)f''(W)+\epsilon_7(W)f'''(W)  +  \epsilon_2(W) \Big(\frac{d^2}{d x^2} \big(   (x^{-}-\beta) f''(x) \big)\big|_{x = \xi_3}\Big)\Big).
\end{align*}
To conclude \eqref{eq:toprovehospec}, we move $- \frac{1}{2}\E\big( b^2(W) + \delta b'(W)\big) f''(W)$   to the left-hand side and note that 
\begin{align*}
& 2\delta + b^{2}(W)-\delta b(W) -\delta^2 - 2\delta^2\beta  - b^2(W) - \delta b'(W) \\
=&\ 2\delta + \delta^2( W^{-} - \beta )^2 -\delta^2( W^{-} - \beta ) -\delta^2 - 2\delta^2\beta  - \delta^2( W^{-} - \beta )^2 + \delta^2 1(W<0) \\
=&\ 2\delta - \delta^2(W^{-}-\beta)- \delta^2 - 2\delta^2\beta + \delta^2 1(W < 0),
\end{align*}
which coincides with the term in front of $f''(W)$ on the left-hand side of \eqref{eq:toprovehospec}.

\startproof{Proof of Lemma~\ref{lem:hospaux}}
Since   $\E g(W') - \E g(W) = 0$, performing a third-order Taylor expansion gives us
\begin{align*}
0 =&\  \E \Delta g'(W) +  \frac{1}{2}  \E\Delta^2 g''(W) + \frac{1}{6} \E\Delta^{3} g'''(\xi) \\
=&\ \E b(W) g'(W) + \frac{1}{2} \E\Big(2\delta + b^{2}(W)-\delta b(W) -\delta^2 - 2\delta^2\beta + \delta^3 W^{-}\Big) g''(W) \notag \\
&+ \frac{1}{6} \E \big(6\delta b(W)  + \delta^{3} \epsilon_0(W)\big)g'''(\xi).
\end{align*}
The second equality is due to Lemma~\ref{lem:hosmec}.
We set  
\begin{align*}
\epsilon_1(x) =&\  -\frac{b^{2}(x)-\delta b(x) -\delta^2 - 2\delta^2\beta + \delta^3 x^{-}}{2\delta^2}  \quad \text{ and } \quad \epsilon_2(x) =  - \frac{  6\delta b(x)  + \delta^{3} \epsilon_0(x) }{6\delta^2}
\end{align*}
and note that both $\epsilon_1(x)$ and $\epsilon_2(x)$ satisfy \eqref{eq:epsdefec} because  $b(x) = \delta(x^{-}-\beta)$.
\finishproof

\subsection{Proof of Lemma~\ref{lem:hosmec} }
\label{sec:hosmecproof}
By the definition of the hospital model the change in customers $\Delta = W'-W$ satisfies 
  \begin{align*}
    \Delta = \delta(A - D),
  \end{align*}
  where $A$ is a mean $\Lambda$ Poisson random variable, and
  conditioned on $W=x=\delta(k-N)$, $D\sim \text{Binomial}(k\wedge N, \mu)$;
    see \cite{DaiShi2017}. To prove the lemma, we utilize the following Stein
  identities for a  mean $\Lambda$ Poisson random variable $X$  and
  a Binomial($k, \mu$) random varabile $Y$:
  \begin{align}
    &   \E X f(X) =\Lambda \E f(X+1) \quad \text{ for each } f:\Z_+\to \R \text{ with } \E \abs{ X f(X)} <\infty, \label{eq:steinpoisson} \\
    &   \E Yf(Y) = \mu  k \E f(Y+1) - \mu \E\Big[Y\Big(f(Y+1)-f(Y)\Big)\Big] \text{ for each } f:\Z_+\to \R. \label{eq:steinbinom}
  \end{align}
See, for example, Lectures VII and VIII of \cite{Stei1986}.

\startproof{Proof of Lemma~\ref{lem:hosmec}}
We prove \eqref{eq:1ec}--\eqref{eq:5ec} in sequence. Using the facts that  $\Lambda=\sqrt{N}-\beta$ and $\mu=\delta=1/\sqrt{N}$, for $x = \delta(k - N)$ we have 
\begin{align}
\E(A-D|W=x)=\Lambda - (k\wedge N) \mu =  \mu \big(\Lambda/\mu - N   -   (k- N \wedge 0)\big)  =&\  -\beta  -  (x \wedge 0) = x^{-} - \beta. \label{eq:6}
\end{align}
For the remainder of the proof we adopt the convention that all expectations are conditioned on $W=x$. We now prove \eqref{eq:2ec}:
  \begin{align}
   \E \Delta^2 &=\delta^2\Big(\E\big[A(A-D)\big]-\E\big[D(A-D)\big]\Big)\nonumber\\
        &=\delta^2\Big(\Lambda\E (A-D+1)-\E\big[D(A-D)\big] \Big) \nonumber\\
        &=\delta^2\Big(\Lambda\E (A-D+1)-\mu(n\wedge N) \E (A-D-1) - \mu  \E D \Big) \nonumber\\
        & =\delta^2\Big(2\Lambda +  ( \Lambda - \mu(n \wedge N) )\E (A-D-1)  + \mu(x^{-} - \sqrt{N})\Big) \nonumber\\
        &= \delta^2\Big(2\Lambda + (x^{-} - \beta)^2 - (x^{-} - \beta) + \mu(x^{-} - \sqrt{N})\Big).\label{eq:8}
  \end{align}
We used \eqref{eq:steinpoisson} in the second equality and \eqref{eq:steinbinom} in the third equality. The fourth and fifth equalities are due to  
\begin{align*}
\E D =&\ \E(D-A)+\E A =-(x^{-} - \beta)+\Lambda = -x^- +  \sqrt{N} \\
&\ \text{ and } \Lambda-\mu(n\wedge N) = \E(A-D) = x^- - \beta, \text{ respectively}.
\end{align*}
Lastly, bringing the $\delta^2$ term in   \eqref{eq:8} inside the parenthesis and recalling that $\delta^2 \Lambda = \delta - \delta^2 \beta$, $\delta (x^- - \beta) = b(x)$, and $\mu = \delta = 1/\sqrt{N}$ yields (\ref{eq:2ec}). We now prove \eqref{eq:3ec}. Note that $\E\Delta^3 = \delta^3 \E (A-D)^3$ equals
  \begin{align*}
& \delta^3\Big(\E\big[A(A-D)^2\big]-\E\big[D(A-D)^2\big]\Big)\\
    =&\delta^3\Big(\Lambda \E(A-D+1)^2-\mu (n\wedge N)\E(A-D-1)^2 + \mu \E D\big[(A-D-1)^2-(A-D)^2\big] \Big) \\
     =&\delta^3\Big(4 \Lambda \E  (A-D) + \big(\Lambda-(n\wedge N)\mu\big) \E (A-D-1)^2  + \mu \E D\big[-2(A-D)+1\big] \Big).
  \end{align*}
  The first equality is due to \eqref{eq:steinpoisson} and \eqref{eq:steinbinom} and to get the second equality we use $\E (A-D+1)^2  = \E (A-D- 1 + 2)^2  = \E\big[(A-D- 1)^2 + 4(A-D)\big]$. Let us analyze the terms above one by one.  First, we have 
  \begin{align*}
\delta^3 4 \Lambda \E  (A-D) =&\ 4 (x^- - \beta) \delta^3 \Lambda = 4 (x^- - \beta) (\delta^2 - \delta^3 \beta) = 4\delta^2(x^- - \beta)  + \delta^3\epsilon(x).
  \end{align*}
 Second, since $\delta(A-D) = \Delta$, we have
 \begin{align*}
 \delta^3 \big(\Lambda-(n\wedge N)\mu\big) \E (A-D-1)^2 =&\ \delta (x^- - \beta) \E \big(\Delta^2 - 2\delta \E \Delta + \delta^2 \big)\\
 =&\ \delta (x^- - \beta) \E \big(\Delta^2 - 2\delta^2(x^- - \beta) + \delta^2 \big) \\
 =&\ \delta (x^- - \beta) \E \Delta^2  + \delta^3 \epsilon(x)\\
 =&\ 2\delta^2(x^- - \beta) + \delta^3 \epsilon(x).
 \end{align*}
 The last equality is due to \eqref{eq:2ec}, which says that $\E \Delta^2 = 2\delta + \delta^2 \epsilon(x)$. Lastly, 
 \begin{align*}
 \delta^3 \mu \E D\big[-2(A-D)+1\big] =&\  \delta^4 \big( 2 \E (A-D)^2  - 2\E A(A-D) + \E D\big) \\
 =&\  \delta^2 \big( 2 \E \Delta^2 - 2\delta^2\Lambda\E (A-D+1) + \delta^2(- x^- + \sqrt{N})\big) \\
 =&\ \delta^3 \epsilon(x).
 \end{align*}
The second equality is due to \eqref{eq:steinpoisson} and the last equality follows from $\E \Delta^2 = \delta \epsilon(x)$ and $\delta^2 \Lambda\E(A-D +1) = \delta (1-\delta \beta)(x^- - \beta+1) = \delta \epsilon(x)$. Putting the pieces together yields $\E \Delta^3 = 6\delta^2 (x^- - \beta) + \delta^3 \epsilon(x)$, which proves \eqref{eq:3ec}. We now prove (\ref{eq:4ec}):
  \begin{align*}
    \E \Delta^4 =&\delta^4\Big(\E\big[A(A-D)^3\big]-\E\big[D(A-D)^3\big]\Big)\\
         &= \delta^4\Big(\Lambda \E\big[(A-D+1)^3\big]-\mu(n\wedge N) \E(A-D-1)^3 \Big)\\
    & \quad + \mu \E\big[D(A-D-1)^3-D(A-D)^3]\Big)\\
    =&\delta^4\Big(\Lambda \E\big[(A-D+1)^3-(A-D-1)^3\big] + \big(\Lambda-(n\wedge N)\mu\big) \E(A-D-1)^3 \\
         &+ \mu \E D\big[(A-D-1)^3-(A-D)^3\big] \Big) \\
    =&\delta^4\Big(\Lambda \E\big[6(A-D)^2+2\big] + (x^{-} - \beta) \E(A-D-1)^3 \\
         &+ \mu \E D\big[-3(A-D)^2+3(A-D)-1\big] \Big).
    \end{align*}
Let us analyze the terms above one by one. First, 
\begin{align*}
\delta^4 \Lambda \E\big[6(A-D)^2+2\big] = \delta^2 \Lambda \E\big[6\Delta^2+2\delta^2\big] = 12\delta^2 + \delta^3 \epsilon(x).
\end{align*}
Second, 
\begin{align*}
\delta^4 (x^{-} - \beta) \E (A-D-1)^3 = \delta (x^{-} - \beta) \big( \E \Delta^3 - 3 \delta \E \Delta^2 + 3 \delta^2 \E \Delta - 1\big) = \delta^3 \epsilon(x),
\end{align*}
and third, 
\begin{align*}
& \delta^4  \mu \E D\big[-3(A-D)^2+3(A-D)-1\big] \\
=&\ \delta^5 \E \big[ 3 (A-D)^3 - 3(A-D)^2 + D  \big] + \delta^5 \E A\big[-3(A-D)^2+3(A-D)\big]\\
=&\ \delta^2 \E \big[ 3 \Delta^3 - 3\delta \Delta^2 + \delta^3 D  \big] + \delta^5 \Lambda \E \big[-3(A-D+1)^2+3(A-D+1)\big]\\
=&\ \delta^3 \epsilon(x).
\end{align*}
Putting the pieces together yields $\E \Delta^4 = 12\delta^2 + \delta^3 \epsilon(x)$. The proof of (\ref{eq:5ec}) is analogous to the proof of (\ref{eq:4ec}) and is omitted. 
\finishproof
  
\section{Companion for the AR(1) Model} \label{app:ar1proof}
In this section we derive the $v_3$ approximation for the AR(1) model and then prove Lemma~\ref{lem:ar1} in Section~\ref{sec:ar1pf}. We recall that $W = \delta(D_{\infty}-R)$, $W' = e^{-\alpha Z}W + \delta\big(X +  R(e^{-\alpha Z} - 1)\big)$, and $\Delta = W' - W$,  where $\delta = \sqrt{\alpha}$, $R = 1/\alpha$,  $X$ and $Z$ are independent unit-mean exponentially distributed random variables that are also independent of $W$, and  $D_\infty > 0$  has the limiting distribution of the AR(1) model defined by \eqref{eq:defar1}. The asymptotic regime we consider is $\alpha \to 0$, so we assume that $\alpha \in (0,1)$.  We recall Lemma~\ref{lem:ar1}:
 \begin{lemma}\label{lem:ar1ec}
 Recall that $\delta = \sqrt{\alpha}$. For any $k \geq 1$, 
 \begin{align*}
 \E(\Delta^{k} | D_{\infty} = d) =  \delta^{k} k! \bigg(1 + \sum_{i=1}^{k} (-1)^{i}d^{i} \prod_{j=1}^{i} \frac{ \alpha }{1 + j \alpha}  \bigg), \quad d > 0.
\end{align*}
\end{lemma} 
We also recall that for $x \geq -1/\sqrt{\alpha}$, 
\begin{align*}
  \E(\Delta^{k} | W = x) =  \E(\Delta^{k} | D_{\infty} = x/\delta + R) =&\  \delta^{k} k! \bigg(1 + \sum_{i=1}^{k} (-1)^{i}\Big( x\sqrt{\alpha}  +  1 \Big)^{i} \prod_{j=1}^{i} \frac{1}{1 + j \alpha}  \bigg),
\end{align*}
and observe that $\E(\Delta^{k} | W = x)  = \delta^k p_{k}(x)$ for some  degree-$k$ polynomial $p_{k}(x)$. To derive $v_3(x)$, we start with the Taylor expansion in \eqref{eq:taylorgeneric} with $n = 4$; i.e.,   for any function $f: \R \to \R$ satisfying $\E f(W') - \E f(W) = 0$,  
\begin{align}
& \delta \E p_{1}(W) f'(W)+\frac{1}{2}\delta^2 \E p_{2}(W) f''(W) + \frac{1}{6} \delta^3 \E  p_{3}(W) f'''(W)+\frac{1}{24}\delta^4 \E  p_{4}(W) f^{(4)}(W) \notag \\
=&\ -\frac{1}{120} \delta^5 \E p_{5}(W) f^{(5)}(\xi). \label{eq:taylorfourth}
\end{align}   
Since $\sup_{\alpha \in (0,1)} \abs{p_{k}(x)} < \infty$ for each $x  \in \R$, the right-hand side is of order  $\delta^5$. When deriving $v_3$, we want it to account for all terms of order $\delta, \ldots, \delta^4$  and treat terms of order   $\delta^5$  as error.   The following lemma is the basis for the $v_3$ approximation. It converts the third and fourth derivative terms in \eqref{eq:taylorfourth} into expressions involving $f''(x)$ plus error. Its proof is  similar to the $v_2$ derivation in Section~\ref{sec:v2def},   so we postpone it until the end of this section. 
\begin{lemma}
\label{lem:ar1taylor}
Define
\begin{align*}
\bar p_3(x) =&\   \frac{1}{6}  \Big( p_3(x) -   \frac{p_1(x)p_4(x)}{ 2 p_2(x)} - \frac{1}{4} \delta  p_2(x)\Big(\frac{p_4(x)}{p_2(x)}\Big)'  \Big),\\
\underline p_2(x) =&\   \Big(\frac{p_2(x)}{2}-\frac{p_1(x)p_3(x)}{3p_2(x)}-\frac{p_2(x)}{6}\Big(\frac{p_3(x)}{p_2(x)}\Big)'\Big).
\end{align*}
Let $W$ and $W'$ be as in Section~\ref{fse5}. If $f \in C^{5}(\R)$ is such that $\E f(W') - \E f(W) = 0$, then 
\begin{align}
&\delta \E p_{1}(W)    f'(W)+\delta^2 \E \Big(\frac{p_{2}(W)}{2}-\frac{p_{1}(W)\bar p_3(W)}{ \underline{p}_2(W)}-\delta \underline{p}_2(W)\Big(\frac{\bar p_3(W)}{\underline{p}_2(W)}\Big)'\Big)f''(W)  \notag \\
=&\  \ -\frac{1}{120} \delta^5 \E p_5(W) f^{(5)}(\xi_1)+ \frac{1}{72} \delta^5 \E p_3(W)  \Big( \frac{  p_4(x) }{p_2(x) }f'''(x) \Big)''\Big|_{x = \xi_2 }  \notag \\
& \hspace{1.2cm}  + \frac{1}{24}\delta^5 \E p_4(W)  \Big( \frac{\bar p_3(x) }{\underline{p}_2(x) }f''(x) \Big)'''\Big|_{x = \xi_3 } - \frac{1}{18} \delta^5 \E p_3(W)  \Big(\frac{p_3(x)}{p_2(x)} \frac{\bar p_3(x) }{\underline{p}_2(x) }f''(x) \Big)'' \Big|_{x = \xi_4 }. \label{eq:arlem}
\end{align}  
\end{lemma}
We choose  
\begin{align*}
\underline{v}_3(x) = \delta^2 \Big(\frac{p_{2}(x)}{2}-\frac{p_{1}(x)\bar p_3(x)}{ \underline{p}_2(x)}-\delta \underline{p}_2(x)\Big(\frac{\bar p_3(x)}{\underline{p}_2(x)}\Big)'\Big)
\end{align*} 
based on the term in front of $f''(W)$ on the left-hand side of \eqref{eq:arlem}, which is the basis for the $v_3$ approximation in Section~\ref{fse5}. Our choice is based on the presumption that all the terms on the right-hand side of \eqref{eq:arlem} are of order $\delta^5$ when $\delta$ is close to zero. We do not prove this claim rigorously in this paper. Nevertheless, our $v_3$ approximation performs quite well numerically.

\startproof{Proof of Lemma~\ref{lem:ar1taylor}}
In Section~\ref{sec:v2def}, we used $a(x), \ldots, d(x)$ to represent $\E(\Delta^{k} | W = x)$; i.e.,  
\begin{align}
b(x) =&\ \E(\Delta | W = x), \quad a(x) = \E(\Delta^{2} | W = x), \quad c(x) = \E(\Delta^{3} | W = x), \notag \\
d(x) =&\ \E(\Delta^{4} | W = x), \quad e(x) = \E(\Delta^{5} | W = x).  \label{eq:abcd}
\end{align}
Since this proof relies heavily on Section~\ref{sec:v2def}, we use this notation  and then convert to use   $\delta^k p_k(x) = \E(\Delta^{k} | W = x)$ at the end. Our starting point is equation \eqref{eq:taylorfourth}, which we recall for convenience:
\begin{align}
&  \E b(W) f'(W)+\frac{1}{2} \E a(W) f''(W) + \frac{1}{6}   \E  c(W) f'''(W)+\frac{1}{24}  \E  d(W) f^{(4)}(W) \notag \\
=&\ -\frac{1}{120}  \E e(W) f^{(5)}(\xi). \label{eq:taylorfourthinpf}
\end{align} 
Our proof relies on several key equations from Section~\ref{sec:v2def}, which we recall as we go. Let $g_{1}(x) = \int_{0}^{x}\frac{d(y)}{a(y)} f'''(y) dy $ and use \eqref{eq:taylorsecond}, or 
\begin{align}
\E b(W) f'(W)+\frac{1}{2}\E a(W)f''(W)    = -\frac{1}{6} \E c(W) f'''(\xi_2), \label{eq:taylorsecondec}
\end{align}
 with $g_{1}(x)$ in place of $f(x)$ there to get
\begin{align*}
\E \frac{b(W)d(W)}{a(W)}f'''(W)+\E \frac{a(W)}{2}\Big(\frac{d(W)}{a(W)}\Big)'f'''(W)+\frac{1}{2}\E d(W)f^{(4)}(W) =  -\frac{1}{6} \E c(W) g_{1}'''(\xi_2).
\end{align*} 
 We multiply both sides  by $1/12$ and subtract the result from \eqref{eq:taylorfourthinpf} to get
\begin{align}
&\E b(W) f'(W)+\frac{1}{2}\E a(W)f''(W) + \E \bar c(W) f'''(W) \notag \\
=&\ \frac{1}{72} \E c(W) g_{1}'''(\xi_2) -\frac{1}{120} \E e(W) f^{(5)}(\xi_1), \label{eq:arintermv3}
\end{align} 
where $\bar c(x) =   \frac{1}{6}  c(x) - \frac{b(x)d(x)}{ 12 a(x)} - \frac{a(x)}{24}\Big(\frac{d(x)}{a(x)}\Big)'   $. Note that $\bar c(x) = \delta^3 \bar p_3(x)$. Next, let $g_2(x) = \int_{0}^{x}\frac{\bar c(y) }{\underline{v}_2(y) }f''(y) dy $, where 
\begin{align*}
\underline{v}_2(x) = \frac{a(x)}{2}-\frac{b(x)c(x)}{3a(x)}-\frac{a(x)}{6}\Big(\frac{c(x)}{a(x)}\Big)', \quad x \in \R
\end{align*}
 is identical to $\underline{v}_2(x)$ defined in \eqref{eq:underlinev2}, and note that $\underline{v}_2(x) = \delta^2 \underline{p}_2(x)$.
We use \eqref{f4}, or 
\begin{align}
& \E b(W)f'(W)+\E\Big(\frac{a(W)}{2}-\frac{b(W)c(W)}{3a(W)}-\frac{a(W)}{6}\Big(\frac{c(W)}{a(W)}\Big)'\Big)f''(W) \notag \\
=&\  \frac{1}{18} \E c(W) g'''(\xi_2) -\frac{1}{24} \E d(W) f^{(4)}(\xi_1), \label{f4ec}
\end{align}
 with $g_2(x)$ in place of $f(x)$ there to get 
\begin{align*}
& \E b(W)\frac{\bar c(W) }{\underline{v}_2(W) }f''(W)+\E \underline{v}_2(W) \Big(\frac{\bar c(W) }{\underline{v}_2(W) }\Big)' f''(W) + \E \underline{v}_2(W) \frac{\bar c(W) }{\underline{v}_2(W) }f'''(W) \notag \\
=&\  \frac{1}{18} \E c(W)   \Big(\frac{c(x)}{a(x)} g_{2}''(x) \Big)'' \Big|_{x = \xi_4 } -\frac{1}{24} \E d(W) g_2^{(4)}(\xi_3).  
\end{align*}  
Subtracting the equation above from  \eqref{eq:arintermv3}, we conclude that  
\begin{align*}
&\E b(W) f'(W)+\E \Big(\frac{a(W)}{2}-\frac{b(W)\bar c(W)}{ \underline{v}_2(W)}-\underline{v}_2(W)\Big(\frac{\bar c(W)}{\underline{v}_2(W)}\Big)'\Big)f''(W)  \notag \\
=&\  \ - \frac{1}{18} \E c(W)   \Big(\frac{c(x)}{a(x)} g_{2}''(x) \Big)'' \Big|_{x = \xi_4 } + \frac{1}{24} \E d(W) g_2^{(4)}(\xi_3) \\
&+ \frac{1}{72} \E c(W) g_{1}'''(\xi_2) -\frac{1}{120} \E e(W) f^{(5)}(\xi_1).
\end{align*}  
To conclude, we  note that $g_2^{(4)}(\xi_3) =\Big( \frac{\bar c(x) }{\underline{v}_2(x) }f''(x) \Big)'''\Big|_{x = \xi_3 }$ and $g_{1}'''(\xi_2) =  \Big( \frac{  d(x) }{a(x) }f'''(x) \Big)''\Big|_{x = \xi_2 }$  and then substitute $\delta p_1(x)$ for $b(x)$, $\delta^2 p_2(x)$ for $a(x)$, etc., where they appear above.
\finishproof

\subsection{Proof of Lemma~\ref{lem:ar1}}
\label{sec:ar1pf}
\startproof{}
Recall that $\Delta = W' - W = \delta \big(D_{\infty}(e^{-\alpha Z} - 1) + X\big)$, so
\begin{align*}
\E(\Delta^{k} | D_{\infty} = d) =&\ \delta^k \E \big(d(e^{-\alpha Z} - 1) + X \big)^{k} =  \delta^{k} \sum_{i=0}^{k} {k \choose i} \E \Big(X^{i}\big( e^{-\alpha Z} - 1\big)^{k-i}\Big) d^{k-i}.
\end{align*}
Since $X$ and $Z$ are independent and exponentially distributed with mean $1$,  we have  $\E X^i = i!$ and
\begin{align*}
{k \choose i} \E \Big(X^{i}\big( e^{-\alpha Z} - 1\big)^{k-i}\Big) =&\ \frac{k!}{(k-i)!} \E \big( e^{-\alpha Z} - 1\big)^{k-i} = \frac{k!}{(k-i)!}  \int_{0}^{\infty} \big(e^{-\alpha z} - 1\big)^{k-i} e^{-z} dz.
\end{align*}
Using integration by parts, 
\begin{align*}
&\int_{0}^{\infty} \big(e^{-\alpha z} - 1\big)^{k-i} e^{-z} dz \\
=&\ (k-i) (-\alpha) \int_{0}^{\infty} \big(e^{-\alpha z} - 1\big)^{k-i-1} e^{-(1+\alpha)z} dz\\
 =&\ (k-i)(k-i-1) (-\alpha)^2 \frac{1}{1 + \alpha} \int_{0}^{\infty} \big(e^{-\alpha z} - 1\big)^{k-i-2} e^{-(1+2\alpha)z} dz\\
 & \ldots \\
 =&\ (k-i)! (-\alpha)^{k-i} \frac{1}{1 + \alpha}\frac{1}{1 + 2\alpha} \ldots \frac{1}{1+(k-i-1)\alpha}  \int_{0}^{\infty}  e^{-(1+(k-i)\alpha)z} dz\\
 =&\ (k-i)! (-\alpha)^{k-i} \frac{1}{1 + \alpha}\frac{1}{1 + 2\alpha} \ldots \frac{1}{1+(k-i)\alpha}.
\end{align*}
\finishproof


\ACKNOWLEDGMENT{We thank Zhuosong Zhang for proving Lemma~\ref{lem:hosmec}.  We thank Yige Hong and Zhuoyang Liu for producing
some figures of this paper.  Xiao Fang is partially supported by Hong
Kong RGC grants 24301617, 14302418 and 14304917, a CUHK direct grant
and a CUHK start-up grant. J. G. Dai is partially supported by NSF grant CMMI-1537795.}


\bibliographystyle{informs2014} 
\bibliography{dai20170728,fang20190831} 



%
%

\end{document}